\documentclass[a4paper,10pt]{article}
\usepackage{epsf}  
\usepackage{amsfonts, amssymb, amsthm, amsmath}  

\usepackage{graphicx}
\usepackage{psfrag}
\newtheorem{thm}{Theorem}  
\newtheorem{conj}[thm]{Conjecture}
\newtheorem{cor}[thm]{Corollary}  
\newtheorem{lemma}[thm]{Lemma}  
\newtheorem{remark}[thm]{Remark}  
\newtheorem{defn}[thm]{Definition}  
\newtheorem{prop}[thm]{Proposition}  
  
\newtheorem{example}[thm]{Example}  
\numberwithin{thm}{section}  
\def\pf{\noindent\emph{Proof: }}  
\def\stop{\hfill$\square$}  
\providecommand{\tot}[1]{\ensuremath{[ #1]}}
\providecommand{\totl}[1]{\ensuremath{\lceil #1\rceil }}
\providecommand{\totb}[1]{\ensuremath{\lfloor #1\rfloor}}
\newcommand{\fun}{\totl}

\newcommand{\ex}{\mathfrak}
\providecommand {\e}[1]{\mathfrak t^{#1}}
\providecommand{\tC}[2]{\ensuremath \mathcal E^\times\left(#2\right)}
\providecommand{\tCp}[2]{\ensuremath{ {}^{+}\mathcal E^{\times}\left(#2 \right)}}

\newcommand{\T}{\mathbf{T}}  
\newcommand{\TC}{{}^{\text{log}}_{\phantom{t}\mathbb T}C}

\newcommand{\dlog}{{}^{\text{log}}d}
\newcommand{\cuspT}{{}^cT}
\newcommand{\logT}{{}^\text{log} T}

\providecommand{\fp}[2]{{}_{\hspace{3pt}#1\hspace{-2pt}}\times_{#2}}
  
\DeclareMathOperator{\dist}{dist}  

\DeclareMathOperator{\id}{id}
\DeclareMathOperator{\expl}{Expl}
\newcommand{\dbar}{\bar{\partial}}

\newcommand{\exploded}{exploded $\mathbb T$ }
\newcommand{\fib}{\bar M}
  
\providecommand{\tl}[2]{\mathbb T^{#1}_{#2}}
\providecommand{\et}[2]{\ensuremath{\mathfrak T^{#1}_{#2}}}
\providecommand{\lrb}[1]{\ensuremath{\left(#1\right)}}
\providecommand{\abs}[1]{\left\lvert #1\right\rvert}  
\providecommand{\norm}[1]{\left\lVert #1\right\rVert}
  
\author{Brett Parker\\  \texttt{parker@math.mit.edu} }  
  
\title{Holomorphic curves in Exploded Torus Fibrations: Compactness}

\begin{document}

\maketitle

\begin{abstract}
The category of exploded torus fibrations is an extension of the category of smooth manifolds in which some adiabatic limits look smooth. (For example, the type of limits considered in tropical geometry appear smooth.) In this paper, we prove a compactness theorem for (pseudo)-holomorphic curves in exploded torus fibrations. In the case of smooth manifolds, this is just a version of Gromov's compactness theorem in a topology strong enough for gluing analysis. 

\end{abstract}

\tableofcontents

\section{Introduction}

As the subject of this paper is holomorphic curves in a new category, the reader desiring rigorous statements will be disappointed in this introduction, and should wait until terms have been defined. 

Holomorphic curves have been a powerful tool in symplectic topology since Gromov's original paper on the subject, \cite{gromov}. 
The basic idea is that given a symplectic manifold, there is a contractible choice of almost complex structure so that the symplectic form is positive on holomorphic planes. The `moduli space' of (pseudo)-holomorphic curves with a choice of one of these almost complex structures is in some sense `smooth', `compact' and finite dimensional, and topological invariants of this `moduli space' give invariants of the underlying symplectic manifold.

 The number of quotation marks required for the above statements is perhaps one indication that the category of smooth manifolds is not the ideal setting for working with holomorphic curves. In particular, in order for our moduli space to be `compact', we have to include  connected `smooth' families of holomorphic curves which exhibit bubbling behavior, and change topology. This  certainly does not fit the usual definition of a smooth family (although it is more natural from the perspective of algebraic geometry). The topology in which the moduli space is compact is also so unnatural to the smooth category that it is difficult to state simply.
 
 A second reason to leave the smooth category is that holomorphic curves can be very difficult to find directly there. A common technique used to find holomorphic curve invariants is to consider a family of almost complex structures which degenerates somehow to a limit in which holomorphic curves become simpler, and then reconstruct holomorphic curve invariants from limiting information. In the category of exploded torus fibrations, some of these degenerations can be viewed as smooth families, and the limiting problem is a problem of the same type. 
 
One example of such a degeneration is that considered in the algebraic case by \cite{Li}. This is a degeneration represented by a flat family of projective schemes so that the fiber at $0$ has two components intersecting transversely over a smooth divisor, and all other fibers are smooth. This corresponds in the symplectic setting to the degeneration used to represent a symplectic sum, used in  \cite{IP} to prove the symplectic sum formula for Gromov-Witten invariants. Similar degenerations have been used to study holomorphic curves in \cite{ruan} and \cite{sft}. More generally, one might consider the case of a flat family in which the central fiber is reduced with simple normal crossing singularities, or the analogous case in the symplectic setting. This type of degeneration can be represented by a smooth family in the category of  exploded torus fibrations. In the exploded category, an object $\ex B$ corresponding to this central fiber will have a smooth part $\totl{\ex B}$ which is equal to the central fiber, and a tropical part $\totb{\ex B}$ which is a stratified integral affine or `tropical' space which has a vertex for every component of the central fiber $\totl{\ex B}$, an edge for every intersection, and an $n$ dimensional face for each complex codimension $n$ intersection of some number of components. Holomorphic curves in $\ex B$ look like usual holomorphic curves projected to $\totl {\ex B}$, and tropical curves projected to $\totb{\ex B}$. 

This tropical part $\totb{\ex B}$ is related to the large complex structure limit studied in some approaches to mirror symmetry such as found in \cite{gross}, \cite{Fukaya} or \cite{kontsevich}. Studying the tropical curves in $\totb{\ex B}$ can sometimes give complete information about holomorphic curve invariants. One example of this is in \cite{mikhalkin}.

    The goal of this paper is to provide the relevant compactness theorems for holomorphic curves in exploded torus fibrations analogous to Gromov's compactness theorem ( stated precisely on page \pageref{completeness theorem}).  This is one step towards having a good holomorphic curve theory for exploded torus fibrations. It is anticipated that the moduli stack of holomorphic curves in exploded fibrations also admits a natural type of Kuranishi structure which together with these compactness results will allow the definition of holomorphic curve invariants such as Gromov Witten invariants.

\section{Definitions}

 \subsection{Tropical semiring and exploded semiring}\label{tropical semiring}
 
  Functions on exploded fibrations will sometimes take values in the following semirings:
 
 \ 
  
  The tropical semiring is a semiring $\e {\mathbb R}$ which is equal to $\mathbb R$ with `multiplication' being the operation of usual addition and `addition' being the operation of taking a minimum. We will write elements of  $\e {{\mathbb R}}$ as $\e x$ where $x\in\mathbb R$. Then we can write the operations as follows
  \[\e x\e y:=\e{x+y}\]
  \[\e x+\e y:=\e{\min \{x,y\}}\]
  $\e 1$ can be thought of as something which is infinitesimally small. With this in mind, we have the following order on $\e{\mathbb R}$:
  
  \[\e x>\e y\text{ when } x<y\]
 
  Given a ring $R$, the exploded semiring $R\e{ {\mathbb R}}$ consists of elements 
  $c\e x$ with $c\in R$ and $x\in\mathbb R$. Addition and multiplication are as follows:
  \[c_1\e xc_2 \e y=c_1c_2\e{x+y}\]
  \[c_1\e x+c_2\e y=\begin{cases}&c_1\e x\text{ if }x<y
  \\&(c_1+c_2)\e x\text{ if }x=y
  \\&c_2\e y\text{ if }x>y\end{cases}\]
 
  It is easily checked that addition and multiplication are associative and obey the usual distributive rule. We will mainly be interested in $\mathbb C\e {\mathbb R}$ and $\mathbb R\e {\mathbb R}$.
 
 There are semiring homomorphisms
 \[ R\xrightarrow{\iota} R\e{\mathbb R}\xrightarrow{\totb{\cdot}}\e{\mathbb R}\]
defined by 
\[\iota(c):=c\e 0\]
\[\totb{c\e x}:=\e x\]
The homomorphism $\totb\cdot$ is especially important. We shall call $\totb{c\e x}=\e x$ the tropical part of $c\e x$. There will be an analogous `tropical part' of `basic' exploded fibrations which can be thought of as the large scale.

Define the positive exploded semiring $R\e{\mathbb R^+}$ to be the sub semiring of $R\e{\mathbb R}$ consisting of elements of the form $c\e x$ where $x\geq 0$. There is a semiring homomorphism 
\[\fun\cdot:R\e{\mathbb R^+}\longrightarrow R\] 
given by 
\[\fun{c\e x}:=\begin{cases}& c \text{ if }x=0
\\&0\text{ if }x>0\end{cases}\]
Note that the above homomorphism can be thought of as setting $\e 1=0$, which is intuitive when $\e 1$ is thought of as infinitesimally small.

\
 
We shall use the following order on $(0,\infty)\e{\mathbb R}$, which again is intuitive if $\e 1$ is thought of as being infinitesimally small and positive:
\[x_1\e {y_1}<x_2\e{y_2}\text{ whenever }y_1>y_2\text{ or } y_1=y_2\text{ and }x_1<x_2\]

\subsection{Exploded $\mathbb T$ fibrations}

The following definition of abstract exploded fibrations is far too general, but it allows us to talk about local models for exploded torus fibrations as abstract exploded fibrations without giving too many definitions beforehand. 

\begin{defn}
An abstract exploded fibration  $\ex B$ consists the following:
\begin{enumerate}
\item A set of points 
\[ \{p\rightarrow \ex B\}\]
\item A Hausdorff topological space $\tot{\ex B}$ with a map
\[ \{p\rightarrow  \ex B\}\xrightarrow{[\cdot]} \tot{\ex B}\]

\item A sheaf of Abelian groups on $\tot{\ex B}$, $\mathcal E^\times(\ex B)$ called the sheaf of exploded functions on $\ex B$ so that:
\begin{enumerate}

\item Each element  $f\in \mathcal E^\times(U)$ is a 
function defined on the set of points over $U$,
\[f:[\cdot]^{-1}(U)\longrightarrow  R^*\e {\mathbb R}\]
where $R^*$ denotes the multiplicatively invertible elements of some ring $R$. (For this paper, we shall use $R=\mathbb C$.) 
\item $\mathcal E^\times(U)$ includes the constant functions if $U\neq\emptyset$.
\item  Multiplication is given by pointwise multiplication in $R^*\e{\mathbb R}$.
\item Restriction maps are given by restriction of functions.

\end{enumerate}
\end{enumerate}
\end{defn}

We shall use the terminology of an open subset $U$ of an abstract exploded fibration $\ex B$ to mean an open subset $U\subset \tot{\ex B}$ along with the set of points in $[\cdot]^{-1}(U)$ and the sheaf $\mathcal E^\times (U)$. The strange notation for the set of points in $\ex B$, $\{p\rightarrow \ex B\}$ is used because it is equal to the set of morphisms of a point into $\ex B$.

\begin{defn}
A morphism $f:\ex B\longrightarrow\ex C $ of abstract exploded fibrations is a map on points
and a continuous map of topological spaces so that the following diagram commutes

\[\begin{array}{ccc} \{p\rightarrow \ex B\}&\xrightarrow{f}& \{p\rightarrow \ex C\}
\\ \downarrow[\cdot]& &\downarrow[\cdot]
\\ \tot{\ex B}&\xrightarrow{\tot f} &\tot{\ex C} 
\end{array}\]
and so that $f$ preserves $\mathcal E^{\times}$ in the sense that if  $g\in\mathcal E^\times(U)$, then  $f\circ g$ is in $\mathcal E^\times\left(\tot f^{-1}(U)\right)$
\[f^*g:=f\circ g\in \mathcal E^\times\left(\tot f^{-1}(U)\right)\] 

\end{defn}

\

Note that the reader should not attempt to work directly with the above definition as it is too general. It simply enables us to talk about examples without giving a complete definition beforehand.
 The next sequence of examples will give us local models for exploded $\mathbb T$ fibrations.

\

\begin{example} \end{example}

Any smooth manifold $M$ is a smooth  exploded $\mathbb T$ fibration. In this case, the underlying toplogical space $[M]$ is just $M$ with the usual topology, the set of points $ \{p\rightarrow  M\}$ the usual set of points in $M$, and $[\cdot]:\{p\rightarrow M\}\longrightarrow\tot M$ the identity map. 	The sheaf $\mathcal E^\times(M)$ consists of all functions of the form $f\e a$ where $f\in C^\infty(M,\mathbb C^*)$ and $a$ is a locally constant $\mathbb R$ valued function.

 The reader should convince themselves that this is just a strange way of encoding the usual data of a smooth manifold, and that a exploded morphism between smooth manifolds is simply a smooth map. This makes the category of smooth manifolds a full subcategory of the category of exploded $\mathbb T$ fibrations.

\

The above example should be considered as a `completely smooth' object. At the other extreme, we have the following `completely tropical' object. 

\begin{example}\end{example}
 The exploded fibration $\ex T^n$ is best described using coordinates $(\tilde z_1,\dotsc, \tilde z_n)$ where each $\tilde z_i\in\mathcal E^\times (\ex T^n)$. The set of points in $\ex T^{n}$,  $\{p\rightarrow \ex T^n\}$ are then identified with $\lrb{\mathbb C^*\e {\mathbb R}}^n$ by prescribing the values of $(\tilde z_1,\dotsc,\tilde z_n)$. The underlying topological space 
 $\tot{\ex T^n}$ is $\e{\mathbb R^n}$  (given the usual topology on $\mathbb R^n$), and
 \[\tot{(\tilde z_1,\dotsc,\tilde z_n)}=(\totb{\tilde z_1},\dotsc,\totb{\tilde z_n})\text{ or }[(c_1\e {a_1},\dotsc, c_n\e {a_n})]=(\e{a_1},\dotsc,\e{a_n})\]
 Note that there is a $(\mathbb C^*)^n$ worth of points $p\longrightarrow \ex T^n$ over every topological point $[p]\in\tot{\ex T^n}$.
 
  Exploded functions in $\tC \infty{\ex T^n}$ can be written in these coordinates as 
  \[c\e y\tilde z^\alpha:=c\e y\prod \tilde z_i^{\alpha_i}\]
   where $c\in\mathbb C^*$, $y\in\mathbb R$ and $\alpha\in \mathbb Z^n$ are all locally constant functions. 
   
   \
  
  The reader should now be able to verify the following observations:
  
\begin{enumerate}
   \item A morphism from $\ex T^n$ to a smooth manifold is simply given by a map $\tot{\ex T^n}\longrightarrow M$ which has as its image a single point in $M$. 
 
\item A morphism from a connected smooth manifold $M$ to $\ex T^n$  has the information of a map $\tot f$ from $M$ to a point  $[p]\in\tot{\ex T^n}$ and a smooth map $f$ from $M$ to the $(\mathbb C^*)^n$ worth of points over $[p]$.  The map $f$ is determined by specifying the $n$ exploded functions \[f^*(\tilde z_i)\in \tC\infty M\]
\end{enumerate}

Note in particular that $\tC\infty M$ is equal to the sheaf of smooth morphisms of $M$ to $\ex T$. 
This will be true in general. A smooth morphism $f:\ex B\longrightarrow \ex T^n$ is equivalent to the choice of $n$ exploded functions in $\tC\infty {\ex B}$ corresponding to $f^*(\tilde z_i)$.

\

Now, a hybrid object, part `smooth', part `tropical'. 

\begin{example}\label{et 11}
\end{example}
 
 The exploded fibration $\et 1{[0,\infty)}:=\et 11$ is more complicated. We shall describe this using the  coordinate $\tilde z\in\mathcal E^\times(\et 11)$. The set points  $\{p\rightarrow \et 11\}$ is identified with $\mathbb C^*\e{\mathbb R^+}$ by specifying the value of $\tilde z(p)\in\mathbb C^*\e{\mathbb R^+}$.
 
 The underlying topological space is given as follows:
 \[\begin{split}\tot {\et 11}&:=\frac{\mathbb C\coprod \e{\mathbb R^+}}{\{0\in\mathbb C\}=\{\e 0\in\e{\mathbb R^+}\}}
 \\&:=\{(z,\e a)\in\mathbb C\times\e{\mathbb R^+} \text{ so that }za=0\}\end{split}\]

\psfrag{ETF1math}{$\totb{\et 11}:=\e{\mathbb R^{+}}$}
\psfrag{ETF1mathC}{$\totl{\et 11}:=\mathbb C$}
\psfrag{ETF1mathtot}{A picture of $\tot{\et 11}$}
\psfrag{ETF1Cs}{$\mathbb C^{*}$}
\includegraphics{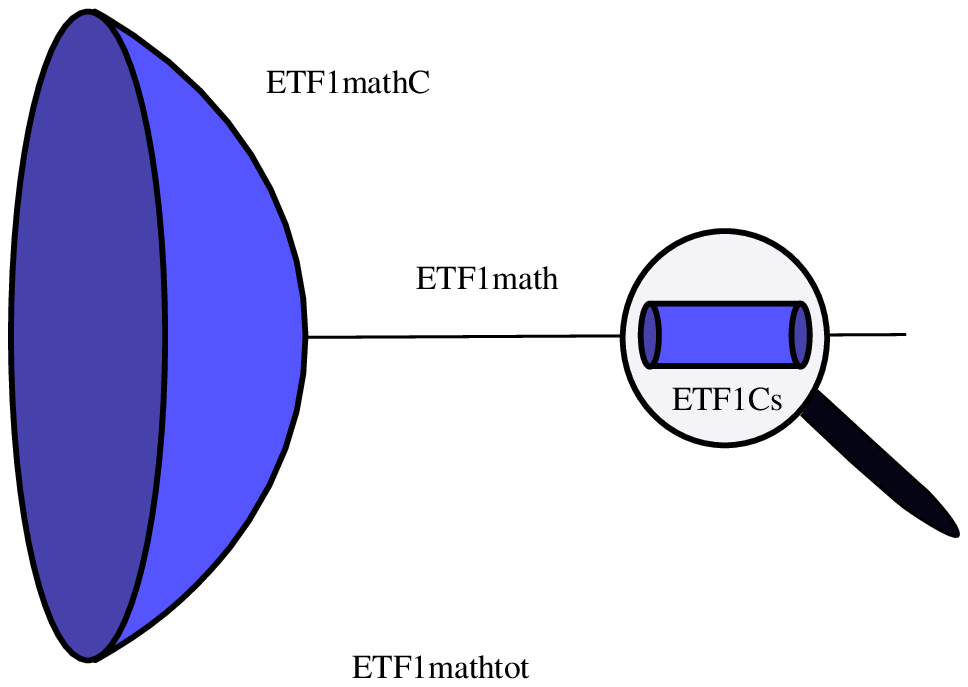}

\

Our map from points to $\tot{\et 11}$ is given by 
\[[\tilde z]:=\left(\fun {\tilde z},\totb{\tilde z}\right)\]
The above expression uses the maps defined in section \ref{tropical semiring}, $\totb{c\e a}=\e a$, and $\fun {c\e a}$ is equal to $c$ if $a=0$ and $0$ if $a>0$. Note that there is one point $p\rightarrow\et 11$ over points of the form $(c,\e 0)\in\tot{\et 11}$, a $\mathbb C^*$ worth of points over $(0,\e a)\in\tot{\et 11}$ for $a>0$ and no points over $(0,\e 0)\in\tot{\et 11} $.

 This object should be thought of as a copy of $\mathbb C$ which we puncture at $0$, and consider this puncture to be an asymptoticaly cylindrical end. Each copy of $\mathbb C^*$ over $(0,\e a)$ should be thought of as some `cylinder at $\infty$'. Note that even though there is a $(0,\infty)$ worth of two dimensional cylinders involved in this object, it should still be thought of being two dimensional. This strange feature is essential for the exploded category to have a good holomorphic curve theory. (Actually, we could just have easily had a $\mathbb Q^{+}$ worth of cylinders at infinity, and worked over the semiring $\mathbb C\e{\mathbb Q}$ instead of $\mathbb C\e {\mathbb R}$, but this author prefers the real numbers.)

We shall use the notation
\[\totl{\et 11}:=\frac{\tot{\et 11}}{\e{\mathbb R^+}=0}=\mathbb C\]
\[\totb{\et 11}:=\frac{\tot{\et 11}}{\mathbb C=0}=\e{\mathbb R^+}\]
We shall call $\totl{\et 11}$ the smooth or topological part of $\et 11$ and $\totb{\et 11}$ the tropical part of $\et 11$.

 We can write any exploded function $h\in\tC \infty {\et 11}$ as 
 \[h(\tilde z)=f(\totl{\tilde z})\e y\tilde z^\alpha\text{ for }f\in C^\infty(\mathbb C,\mathbb C^*)\text{, and } y\in\mathbb R,\alpha \in\mathbb Z\text{ locally constant.}\]

We can now see that $\et 11$ is a kind of hybrid of the last two examples. Restricting to an open set contained inside $\{\totb{\tilde z}=\e 0\}\subset\tot{\et 11}$ we get part of a smooth manifold. Restricting to an open set contained inside $\{\totl{\tilde z}=0\}\subset\tot{\et 11}$, we get part of $\ex T$.

\

We can define  $\tCp\infty{\ex B}$ to consist of all functions in $\mathcal E^\times (\ex B)$ which take values only in $\mathbb C^*\e{\mathbb R^+}$. 
Note that a smooth morphism $f:\ex B\longrightarrow\et 11$ from any  exploded $\mathbb T$ fibration $\ex B$ that we have described in examples so far is given by a choice of exploded function 
$f^*(\tilde z)\in \tCp\infty{\ex B}$. This will be true for exploded $\mathbb T$ fibrations in general.

\
We can build up `exploded curves' using coordinate charts modeled on open subsets of $\et 11$. For an example of this see the picture on page \pageref{curve picture}. For an example of a map of this to $\ex T^{n}$ see example \ref{curve in tropical plane} on page \pageref{curve in tropical plane}. Those interested in tropical geometry may like the following example \ref{tropical curve}. For a local example of a family of such objects, see example \ref{local family model} on page \pageref{local family model}. 

\

 We shall use the special notation $\et mn$ to denote $(\et 11)^n\times\ex T^{(m-n)}$. (These spaces will be described rigorously below.)
Before giving the most general local model for exploded $\mathbb T$ fibrations, we give the following definition which can be thought of as a local model for the tropical part of our exploded fibrations. 

\begin{defn}
   An integral cone $A\subset \e{\mathbb R^m}$ is a subset of $\e{\mathbb R^m}$ given by a collection of inequalities 
   \[A:=\left\{\e a:=(\e{a_1},\dotsc,\e{a_m})\in \e{\mathbb R^m} \text { so that }\e {a\cdot\alpha^i}:=\prod_j (\e{a_j})^{\alpha_j^i}\in\e{\mathbb R^+}\right\}\] 
   \[\text{In other words, }A:=\{\e a\text{ so that }a\cdot\alpha^i\geq 0\}\]
   where $\alpha^i$ are vectors in $\mathbb Z^m$. 
   
   The dual cone $A^*$ to $A$ is defined to be 
   \[A^*:=\{\alpha\in \mathbb Z^m: \alpha\cdot a\geq 0, \forall \e a\in A\} \]
    \end{defn}

\includegraphics{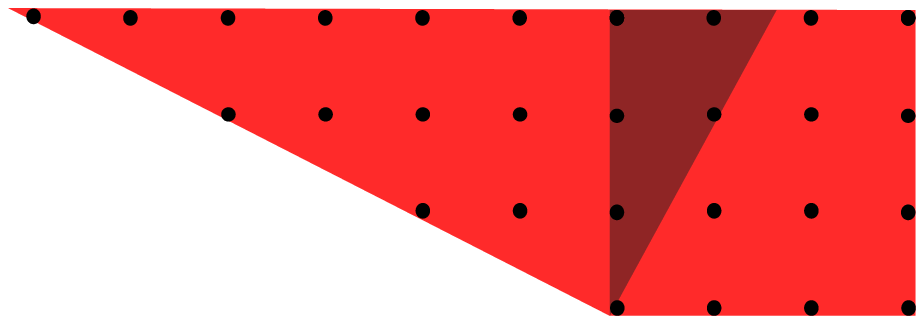}

\

The following is the most general local model for exploded $\mathbb T$ fibrations.

\begin{example}\label{smooth local model}\end{example}

Given any integral cone $A\subset\e{ \mathbb R^m}$ which generates $\e{\mathbb R^m}$ as a group, we will now describe the exploded fibration $\mathbb R^n\times\et mA$.  We shall use coordinates  

 \[(x_1,\dotsc,x_n,\tilde z_1,\dotsc,\tilde z_m)\] 
The set of points $ \{p\rightarrow \mathbb R^n\times\et mA\}$ is then given by 
\[\left\{(x_1,\dotsc,x_m,c_1\e{a_1},\dotsc,c_m\e{a_m})\in\mathbb R^n\times\left(\mathbb C^*\e{\mathbb R}\right)^m ;(  \e {a_1},\dotsc,\e {a_m})\in A\right\}\]

 We will describe the underlying topological space $\tot{\mathbb R^n\times\et mA}$ using some extra coordinate choices. Let $\{\alpha^1,\dotsc,\alpha^k\}$  generate the dual cone $A^*\subset\mathbb Z^m$ as a semigroup with the relations 
 \[\sum_i r_{i,j}\alpha^i=0\] 
Label coordinates for  $\mathbb C^k$ as $(z^{\alpha^1},\dotsc,z^{\alpha^k})$, and use $\e a$ to denote points in $A$. We can describe the underlying topological space, $\tot{\mathbb R^n\times\et mA}$ as a subspace 
\[\tot{\mathbb R^n\times\et mA}=\left\{(x,z^{\alpha^i},\e a);\ \prod_i (z^{\alpha^i})^{r_{i,j}}=1, (\alpha^i\cdot a)z^{\alpha^i}=0\right\}\subset \mathbb R^n\times\mathbb C^k\times A\]
Again we shall use the notation 
\[\totb{\mathbb R^n\times\et mA}:= A,\text{ the tropical part  }\]
\[\totl{\mathbb R^n\times\et mA}:=\left\{\prod_i(z^{\alpha^i})^{r_{i,j}}=1\right\}\subset\mathbb R^n\times\mathbb C^k,\text{ the smooth part}\]
(This smooth part is equal to $\mathbb R^{n}$ times the toric variety corresponding to the cone $A$.)
Using the shorthand $\tilde z^\alpha:=\prod_i\tilde z_i^{\alpha_i}$, and recalling that $\fun {c\e a}$ is equal to $c$ if $a=0$, and equal to $0$ if $a>0$, we can write the map from the set of points down to the smooth part as
\[\totl{(x,\tilde z)}:=(x,\fun {\tilde z^{\alpha^1}},\dotsc,\fun {\tilde z^{\alpha^k}})\]
The map from the set of points down to the tropical part  is given by 
\[\totb{(x, \tilde z)}:=(\totb{\tilde z_1},\dotsc,\totb{\tilde z_m})\text{ or  }\totb{(x,c_1\e {a_1},\dotsc,c_m\e{a_m})}=(\e {a_1},\dotsc,\e{a_m})\]
The map from the set of points down to the underlying topological space is then defined by 

\[\tot{(x,\tilde z)}:=\left(\totl{(x,\tilde z)},\totb{(x,\tilde z)}\right)\]
(Note that we can also describe $\tot{\mathbb R^{n}\times \et mA}$ as the closure of the image of the above map inside $\mathbb R^{n}\times\mathbb C^{k}\times A$.)
Exploded functions  $h\in\tC \infty{\mathbb R^n\times\et mA}$  are all functions that can be written as 

\[h(x,\tilde z):=f\left(\totl{(x,\tilde z)}\right)\tilde z^\alpha\e y\]
 where $f\in C^\infty(\mathbb R^n\times \mathbb C^k,\mathbb C^*)$, and  $\alpha\in\mathbb Z^m$ and $y\in \mathbb R$ are locally constant. It is important to note that two exploded functions are equal if they have the same values on points (which means that we can represent the same exploded function using different smooth functions $f$.) Note also that this definition is independent of the choice  of basis $\{\alpha^i\}$ for $A^*$. 

\

The earlier examples $\et 11=\et 1 {\e{[0,\infty)}}$, and $\mathfrak T^n=\et n{\e{\mathbb R^n}}$, and \[\et mn=(\et 11)^{n}\times\ex T^{(m-n)}:=\et m{\e{(\mathbb R^{+})^{n}\times\mathbb R^{(m-n)}}}\]
These examples are much easier to understand than the most general model, so if the above definition was confusing, it is worthwhile reading through it again with these examples in mind. In this case, we can choose the obvious basis for $A^{*}$ given by the first $n$ standard basis vectors, so $\tilde z^{\alpha^{i}}:=\tilde z_{i}$. The smooth part of $\et mn$, is then  $\totl{\et mn}:=\mathbb C^{n}$ with coordinates $(\totl{\tilde z_{1}},\dotsc,\totl{\tilde z_{n}})$. An important example of exploded fibrations locally modeled on open subsets of $\et mn$ is given in example \ref{expl} on page \pageref{expl}. 
We need the more confusing models $\et mA$ because they arise naturally in the intersection theory of exploded fibrations. They also occur in moduli spaces of holomorphic curves. We need to consider families of holomorphic curves with total spaces locally modeled on these more complicated $\et mA$.

  \begin{defn}
   A smooth exploded $\mathbb T$ fibration $\ex B$ is an abstract exploded fibration locally modeled on open subsets of $\mathbb R^n\times \et mA$. 
    
    In other words, each point in $\tot {\ex B}$ has a neighborhood $U$ which is isomorphic as an abstract exploded fibration to some open set in $\mathbb R^n\times \et mA$. The model $\mathbb R^n\times \et mA$ used may depend on the open set.
  \end{defn}

 There are some extra sheaves of functions which are defined on any smooth exploded $\mathbb T$ fibration $\ex B$ which will come in useful:
 
 \begin{defn}
 The sheaf of smooth functions $C^\infty(\ex B)$ is the sheaf of smooth morphisms of $\ex B$ to $\mathbb R$ considered as a smooth \exploded fibration. 
 \end{defn}
 
 \begin{defn}
 The sheaf of integral affine functions $\mathcal A^\times(\ex B)$ consists of functions of the form $\totb{f}$ where $f\in\mathcal E^\times(\ex B)$. 
 \end{defn}
 
 Recall that $\totb{c\e a}:=\e a$. We shall also use the terminology that the order of $c\e a$ is $a$. The sheaf $\mathcal A^{\times}(\et mA)$ corresponds to  functions on $A$ of the form 
 $\e a\mapsto\e{y+\alpha\cdot a}$ where  $\alpha\in\mathbb Z^m$ and $y\in\mathbb R$. 
 
 \begin{defn}
 The sheaf of exploded tropical functions $\mathcal E(\ex B)$ is the sheaf of $\mathbb C\e {\mathbb R}$ valued functions which are locally equal to a finite sum of exploded functions $\mathcal E^\times$. (`Sum' means sum using pointwise addition in $\mathbb C\e {\mathbb R}$.) The sheaf of tropical functions $\mathcal A(\ex B)$ consists of functions of the form $\totb{f}$ for some $f\in\mathcal E(\ex B)$.
 \end{defn}
 
 For example $\mathcal A(\ex T^n)$ consists of functions of the form $\e {f}$ where $f:\mathbb R^n\longrightarrow \mathbb R$ is some continuous, concave piecewise integral affine function. Functions of this form are usually called tropical functions. We shall usually not use $\mathcal E(\ex B)$ because the addition is strange, and $\mathcal E^{\times}(\ex B)$ is naturally defined as the sheaf of morphisms to $\ex T$. The only place that $\mathcal E(\ex B)$ will be used in this paper will be in the definition of the tangent space of $\ex B$. The operation of addition is needed here to state the usual property of being a derivation. The other reason that addition was mentioned in this paper was to emphasize the links with tropical geometry.
 
 \
    
   The inclusion $\iota:\mathbb C\longrightarrow \mathbb C\e{\mathbb R}$ defined by $\iota (c)=c\e 0$ induces an inclusion of functions
    \[\iota: C^\infty(\ex B)\longrightarrow \mathcal E(\ex B)\]
    \[\iota(f)(p):=\iota(f(p))\] 
    
     \begin{defn}
      Using the notation of example \ref{smooth local model}, a $C^k$ function on $\mathbb R^n\times\et mA$, $f\in C^k(\mathbb R^n\times \et mA)$ is a function  on $\tot{\mathbb R^{n}\times\et mA}$ which is equal to a $C^k$ function applied to the variables  $ x_i$ and $z^{\alpha^i}$. A $C^k$ exploded function $f\in \mathcal E^{k\times}(\mathbb R^n\times \et mA)$ is a function of the form $f\tilde z^\alpha\e y$ where $f$ is a $C^{k}$ function taking values in $\mathbb C^{*}$. We can define $C^{k}$ morphisms of exploded fibrations to be morphisms of abstract exploded fibrations using the sheaf $\mathcal E^{k\times}$ instead of $\mathcal E^{\times}$
      \end{defn}
      
      We can consider an arbitrary Hausdorff topological space to be a $C^{0}$ exploded fibration. A continuous morphism of $\ex B$ to a Hausdorff topological space $X$ is given by a map $f:\ex B\longrightarrow \tot{\ex B}\longrightarrow X$ so that  the pullback of continuous real valued functions on $X$ are $C^{0}$ functions on $\ex B$. Note that in general, the map  $\tot\cdot:\ex B\longrightarrow \tot{\ex B}$ is not a continuous morphism in this sense.

     \begin{defn}
     The topological or smooth part of $\ex B$, $\totl{\ex B}$ is the unique Hausdorff topological space $\totl{\ex B}$ with the following universal property: There is a $C^{0}$ morphism 
     \[\ex B\xrightarrow{\totl{\cdot}}\totl{\ex B}\]
     so that any other $C^{0}$ morphism from $\ex B$ to a Hausdorff topological space factors uniquely through $\totl{\cdot}$:
      
     \[\text{Given }f:\ex B\longrightarrow X\]
     \[\text{there exists a unique }h:\totl{\ex B}\longrightarrow X\]
     \[\text{so that }\ \ \ \ \ f=h\circ\totl\cdot\]

     \end{defn}
     
     We can construct $\totl {\ex B}$ as the image of $\ex B$ in $\mathbb R^{C^{0}(\ex B)}$. The universal property of $\totl{\cdot}$ makes it clear that taking the topological part of $\ex B$ is a functor from the exploded category to the category of Hausdorff topological spaces.

\begin{defn}\label{topological convergence}
We shall use the adverb `topologically' to indicate a property of $\totl{\ex B}$. In particular, 
\begin{enumerate}
 \item A sequence of points $p_{i}\longrightarrow \ex B$ converges topologically to $p\longrightarrow \ex B$ if their images converge in $\totl{\ex B}$. In this case we use the notation $\totl{p_i}\rightarrow\totl{p}$ or say these points converge in $\totl{\ex B}$. Note that this is equivalent to $f(p_{i})$ converging to $f(p)$ for every continuous $\mathbb R$ valued function $f\in C(\ex B)$.
\item The exploded fibration $\ex B$ is topologically compact if $\totl{\ex B}$ is compact.
\item A map $f:\ex B\longrightarrow \ex C$ is topologically proper if the induced map $\totl{f}:\totl{\ex B}\longrightarrow\totl {\ex C}$ is proper.
\end{enumerate}
\end{defn}

Note that limits in the above sense are of course not unique.  A manifold is topologically compact when considered as an exploded fibration if and only if it is compact. Another example is any connected open subset of $\ex T^n$. Sometimes the topological part of $\ex B$ gives slightly misleading information about $\ex B$, as the following example is intended to demonstrate. We shall give a more useful analogue of `compactness' in definition \ref{complete} shortly.

  \begin{example}\label{notbasic}\end{example}
  Our local model $\mathbb R^{n}\times \et mA$ has a smooth part $\totl{\mathbb R^{n}\times \et mA}$ which is equal to the product of $\mathbb R^{n}$ with the toric variety corresponding to the integral affine cone $A$.
  The following two examples have `misbehaved' smooth parts: Consider the exploded fibration formed by quotienting $\mathbb R\times \ex T$ by $(x,\tilde z)=(x+1,\tilde z)$ and $(x,\tilde z)=(x+\pi,\e 1\tilde z)$. The smooth part of this is a single point, even though it is made up of pieces which have smooth parts equal to a circle. Another example is given by the subset of $\mathbb R\times \ex T$ which is the union of the set where $\totb {\tilde z}<\e 0$ with the set where $\totb{\tilde z}>\e 1$ and $x>0$. (Note that these sets overlap because $\e 0>\e 1$.)  The smooth part of this is equal to $\mathbb R$.
      
      \
      
    For many of our proofs, we shall restrict to a class of exploded fibrations in which the smooth part and the tropical part of the exploded fibration give better information than the above two examples. To explain this class of `basic' exploded fibrations, we shall first construct some elementary examples.

     \begin{example}\label{polygon}
     \end{example}
   
     Suppose that $P\subset \e{\mathbb R^{m}}$ is a convex integral affine polygon with nonempty interior defined by some set of inequalities:
     \[P:=\{\e a \in \e{\mathbb R^{m}}\text{ so that }\e {c_{i}+a\cdot\alpha^{i}}\leq \e{0},\e {c_{i'}+a\cdot\alpha^{i'}}< \e{0} \}\]
     where $\alpha^{i}\in \mathbb Z^{m}$ and $c_{i}\in\mathbb R$. This can simply be viewed as a convex polyhedron in $\mathbb R^{m}$ which has faces with rational slopes. View $P$ as a stratified space in the standard way, (so the strata of $P$ consist of the vertices of $P$, the edges of $P$ minus the vertices of $P$, the two dimensional faces of $P$ minus all vertices and edges, and so on, up to the  interior of $P$).  Note that as some of the defining inequalities may be  strict, $P$ may be missing some faces or lower dimensional strata.
     
     Then there is an exploded fibration $\et mP$ associated with $P\subset \e{\mathbb R^{m}}$, which is constructed by gluing together coordinate charts modeled on $\et m{A_{k}}$ where $A_{k}$ are the integral affine cones corresponding to the strata  of $P$. 
     
     Explicitly, suppose that the point $\e{a^{k}}\in\e{\mathbb R^{m}}$ is inside the $k$th strata of $P$, and $A_{k }$ is the integral affine cone so that locally near $\e{a^{k}}$, $P$ is equal to $A_{k}$ shifted by $\e{a^{k}}$. 
     \[\text{near }\e{a^{k}}, \ \ \ \  P=\{\e a\in \e{\mathbb R^{m}}\text{ so that }\e{a-a^{k}}\in A_{k}\}\] 
    
 There are coordinates $\tilde z_{1},\dotsc,\tilde z_{m}\in\mathcal E^{\times}\lrb{\et mP}$ so that the set of points in $\et mP$ are given by specifying values $(\tilde z_{1}(p),\dotsc,\tilde z_{m}(p))\in\lrb{\mathbb C^{*}\e{\mathbb R}}^{m}$ so that 
     \[(\totb{\tilde z_{1}(p)},\dotsc,\totb{\tilde z_{m}(p)})\in P\]
   
        Denote by $P_{k}\subset P$ the subset of $P$ which is the union of all strata of $P$ whose closure contains the $k$th strata of $P$.  There is a corresponding subset $\et m{P_{k}}$ of $\et mP$:
        \[\et m {P_{k}}:=\left\{\tilde z\text{ so that } \totb{\tilde z}\in P_{k}\right\}\]
     We put an exploded structure on  $\et m {P_{k}}$ so it is isomorphic to the corresponding subset of $\et m{A_{k}}$ with the standard coordinates given by $(\e{-a^{k}_{1}}\tilde z_{1},\dotsc,\e{-a^{k}_{m}}\tilde z_{m})$. It is an easy exercise to check that the restriction of the exploded structure this gives to $\et m{P_{k}}\cap \et m{P_{k'}}$ is compatible with the restriction of the one coming from $\et m{P_{k'}}$, so this gives a globally well defined exploded structure. 
     
     An alternative way to describe this structure is as follows: Consider the collection of functions $\tilde w:=\e {c}\tilde z^{\alpha}$ so that on $P$, $\totb{\tilde w}\leq \e 0$. Choose some finite generating set $\{\tilde w_{1},\dotsc,\tilde w_{n}\}$ of these functions so that any other function of this type $\tilde w$ can be written as $\tilde w=\e{c}\tilde w_{1}^{\beta_{1}}\dotsb \tilde w_{n}^{\beta_{n}}$  where $\beta_{i}\in\mathbb N$ and $c\in [0,\infty)$. Exploded functions on $\et mP$ can then be described as functions of the form $f(\totl{\tilde w_{1}},\dotsc,\totl{\tilde w_{n}})\e{c}\tilde z^{\alpha}$, where $f:\mathbb C^{n}\longrightarrow \mathbb C^{*}$ is smooth. The underlying topological space $\tot{\et mP}$ is equal to the closure of the set 
     \[\{\totl{\tilde w_{1}},\dotsc,\totl{\tilde w_{n}},\totb{\tilde z_{1}},\dotsc,\totb{\tilde z_{m}}\}\subset\mathbb C^{n}\times P\]
     
     The tropical part of $\et mP$ is given by 
     \[\totb{\et mP}=P=\{\totb{\tilde z_{1}},\dotsc,\totb{\tilde z_{m}}\}\subset \e{\mathbb R^{m}}\]
 The smooth part of $\et mP$ is given by 
 \[\totl{\et mP}:=\{\totl{\tilde w_{1}},\dotsc,\totl{\tilde w_{n}}\}\subset \mathbb C^{n}\]
     
     Simple examples of the above construction are $\et 1{(a,b)}$ which is equal to the subset of $\ex T$ where $\e a>\totb{\tilde z}>\e b$, and $\et 1{[0,b)}$ which is equal to the subset of $\et 11$ where $\totb{\tilde z}>\e b$. (Strictly speaking, we should be writing $\et 1{\e{(a,b)}}$, as the above confuses the interval $(a,b)\subset \mathbb R$ with $\e{(a,b)}\subset\e{\mathbb R}$.) 
     
     \
     
     Similarly, we can construct $\mathbb R^{k}\times \et mP$. This has coordinates $(x,\tilde z)\in\mathbb R^{k}\times \et mP$, has underlying topological space $\tot{\mathbb R^{k}\times \et mP}=\mathbb R^{k}\times \tot{\et mP}$, so that $\tot{(x,\tilde z)}=(x,\tot{\tilde z})$, and has exploded functions equal to functions of the form $f(x,\totl{\tilde w_{1}},\dotsc,\totl{\tilde w_{n}})\e c\tilde z^{\alpha}$, where $f:\mathbb R^{k}\times\mathbb C^{n}\longrightarrow \mathbb C^{*}$ is smooth.  
     
     \begin{example}\label{polygon2}\end{example}
     Given a smooth manifold $M$ and  $m$ complex line bundles on $M$, we can construct the exploded fibration $M\rtimes \et mP$ as follows: Denote by $E$ the corresponding space of  $m$ $\mathbb C^{*}$ bundles over $M$. This has a smooth free $\lrb{\mathbb C^{*}}^{m}$ action. The exploded fibration $\et mP$ also has a $\lrb{\mathbb C^{*}}^{m}$ action given by multiplying the coordinates $\tilde z_{1},\dotsc,\tilde z_{m}$ by the coordinates of $\lrb{\mathbb C^{*}}^{m}$. $M\rtimes \et mP$ is the exploded fibration constructed by taking the quotient of $E\times \et mP$ by the action of $\lrb{\mathbb C^{*}}^{m}$ by $(c,c^{-1})$. As the action of $\lrb{\mathbb C^{*}}^{m}$ is trivial on $\totb{\et mP}$, the tropical part of $M\rtimes\et mP$ is still defined, and is equal to $P$.
     
     Alternately, choose coordinate charts on $E$ equal to $U\times\lrb{\mathbb C^{*}}^{m}\subset \mathbb R^{n}\times \mathbb C^{n}$. The transition maps are of the form $(u,z_{1},\dotsc,z_{m})\mapsto(\phi(u),f_{1}(u)z_{1},\dotsc,f_{m}(u)z_{m})$. Replace these coordinate charts with $U\times \et mP\subset\mathbb R^{n}\times\et mP$, and replace the above transition maps with maps of the form $(u,\tilde z_{1},\dotsc,\tilde z_{m})\mapsto(\phi(u),f_{1}(u)\tilde z_{1},\dotsc,f_{m}(u)\tilde z_{m})$. The map to the tropical part $\totb{M\rtimes\et mP}:=P$ in these coordinates is given by 
     \[\totb{(u,\tilde z_{1},\dotsc,\tilde z_{n})}=(\totb{\tilde z_{1}},\dotsc,\totb{\tilde z_{n}})\in P\]

\begin{defn}\label{basic definition}
  
  An exploded fibration $\ex B$ is basic if

  \begin{enumerate}\item There exists a Hausdorff topological space $\totb{\ex B}$ called the tropical part of $\ex B$, along with a map
  \[\totb{\cdot}: \{p\rightarrow  \ex B\}\longrightarrow\tot{\ex B}\longrightarrow \totb{\ex B}\]
  
  \item The space $\totb{\ex B}$ is a union of strata $\totb{\ex B_i}$ so that 
  \begin{enumerate}
  \item Each strata $\totb{\ex B_i}\subset\totb{\ex B}$ is an integral affine space equal to some open convex integral affine polygon  in $\e{\mathbb R^k}$ of the form considered in example \ref{polygon}. (The dimension $k$ depends on the strata.) 
  \item The closure of $\totb{\ex B_i}\subset\totb{\ex B}$  is a union of strata which is also equal to an integral affine polygon $\overline{\totb{\ex B_{i}}}$ of the form considered in example \ref{polygon} with some strata identified.
  \end{enumerate}
\item Use the notation $\ex B_{i}$ for inverse image of $\totb{\ex B_{i}}$ under the map $\totb\cdot$, and call this a strata of $\ex B$.
\begin{enumerate}
\item Each strata $\ex B_{i}$ has a smooth part $\totl{\ex B_{i}}$ which is a smooth manifold.
\item Each strata $\ex B_{i}$ is equal to an exploded space $\totl{\ex B_{i}}\rtimes\et k{\totb{\ex B_{i}}}$ using the construction of example \ref{polygon2}. The map $\totb\cdot$ restricted to $\ex B_{i}$ is the map $\totb\cdot :\totl{\ex B_{i}}\rtimes\et k{\totb{\ex B_{i}}}\longrightarrow \totb{\ex B_{i}}$ from example \ref{polygon2}.
\item There is a neighborhood of each strata $\ex B_{i}$ which is open in $\totl {\ex B}$ which is isomorphic to some subset of $\totl{\ex B_{i}}\rtimes\et k{\overline{\totb{\ex B_{i}}}}$ which is a topologically open neighborhood of $\totl{\ex B_{i}}\rtimes\et k{\totb{\ex B_{i}}}\subset {\totl{\ex B_{i}}\rtimes\et k{\overline{\totb{\ex B_{i}}}}}$.
\end{enumerate}   

  \end{enumerate}
\end{defn}

It is an easy exercise to check that if $\ex B$ is basic, the smooth part $\totl{\ex B}$ is also a stratified space, with strata equal to the manifolds $\totl{\ex B_{i}}$. The stratification of $\totl{\ex B}$ is in some sense dual to the stratification of $\totb{\ex B}$ in that $\totb{\ex B_{i}}$ is in the closure of $\totb{\ex B_{j}}$ if and only if $\totl{\ex B_{j}}$ is in the closure of $\totl{\ex B_{i}}$.

    All our examples up until now apart from those introduced in example \ref{notbasic} have been basic. Another example of an exploded fibration that is not basic is given by quotienting $\mathbb R\times \ex T$ by $(x,\tilde z)=(x+1,\tilde z^{-1})$. A less pathological example of an exploded fibration that is not basic is  given by the quotient of $\ex T$ by $\tilde z=\e 1\tilde z$. Any exploded fibration is locally basic.

   \
     
     To summarize the topological spaces we have involved in the case where $\ex B$ is basic, we have the following topological spaces. 
     
     \[\begin{split}\{p\rightarrow\ex B\}\longrightarrow &\tot{\ex B}\longrightarrow\totl{\ex B}
     \\ &\downarrow 
     \\ &\totb{\ex B}\end{split}\]
      \begin{enumerate}
      \item The smooth or topological part $\totl {\ex B}$. Any continuos morphism of $\ex B$ to a Hausdorff topological space factors through $\totl{\ex B}$. The smooth part of a smooth manifold $M$ is  $M$. The smooth part of $\ex T^{n}$ is a single point.
     \item The tropical part $\totb {\ex B}$. This is a stratified integral affine space, which can be thought of as what $\ex B$
     looks like on the large scale. The tropical part of a connected smooth manifold is a single point. The tropical part of $\ex T^{n}$ is $\e{\mathbb R^{n}}$.
     \item The underlying topological space $\tot{\ex B}$. This strange topological space is a mixture of $\totl{\ex B}$ and $\totb{\ex B}$. It is the topological space on which $\mathcal E^{\times}$ is a sheaf, and this is the topology in which $\ex B$ is locally modeled on $\mathbb R^{n}\times \et mA$.
     \item In the case that $\ex B$ is not a manifold, the set of points $\{p\rightarrow\ex B\}$ is not Hausdorff in any of the above topologies. We can put a topology on this set by choosing the strongest topology in which all continuous morphisms from manifolds into $\ex B$ are continuous maps to ${\{p\rightarrow\ex B\}}$. For example  $\{p\rightarrow \et 11\}$ is homeomorphic to the disjoint union of $[0,\infty)$ copies of  $\mathbb C^{*}$. (When we have defined tangent spaces and metrics on exploded fibrations, this would be the topology coming from any metric.)
     
      \end{enumerate}
      
       \
       
       The following is a picture  of an `exploded curve' built up using coordinate charts modeled on open subsets of $\et 11$. Note that $\totl{\ex C}$ is what we would see if we mapped $\ex C$ to a smooth manifold. On the other hand, we would see $\totb{\ex C}$ if we mapped $\ex C$ to $\ex T^{n}$. \label{curve picture}
       
       \psfrag{ETF3tot}{$\tot{\ex C}$}
      \psfrag{ETF3totl}{$\totl{\ex C}$}
      \psfrag{ETF3totb}{$\totb{\ex C}$}
      
       \includegraphics{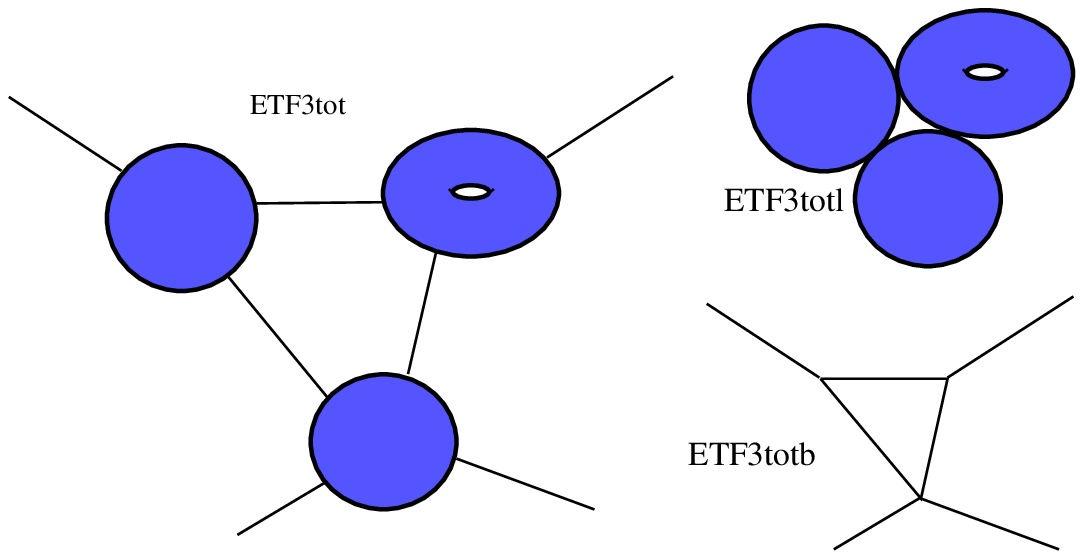}

\begin{defn}\label{complete}
An exploded $\mathbb T$ morphism $f:\ex B\longrightarrow \ex C$ is complete if the following conditions hold
\begin{enumerate}
\item $f$ is topologically proper.
\item Any map \[g:\et 1{(0,l)}\longrightarrow \ex B\]
extends to a continuous map $\et 1{[0,l)}\longrightarrow\ex B $ if and only if the map $f\circ g$ extends to a continuous map $\et 1{[0,l)}\longrightarrow \ex C$.
\end{enumerate}
We say that an exploded fibration $\ex B$ is complete if the map from $\ex B$ to a point is complete.
\end{defn}

The notation $\et 1{(0,l)}$ is that introduced in example \ref{polygon}.
An equivalent way of stating the second condition is that the map $\tot g:(0,l)\longrightarrow \tot{\ex B}$ extends to a continuous map $[0,l)\longrightarrow \tot{\ex B}$ if and only if $\tot{f\circ g}:(0,l)\longrightarrow\tot{\ex C}$ extends to a continuous map $[0,1)\longrightarrow \tot{\ex C}$ 

Complete maps should be thought of as being analogous to proper maps.
A smooth manifold is complete if and only if it is compact. An example of a non-compact exploded fibration which is complete is $\ex T^n$.  In the case that $\ex B$ is basic, completeness means that the topological part $\totl{\ex B}$ is compact, and the tropical part $\totb{\ex B}$ is a complete stratified integral affine space.

 An interesting example of an exploded fibration which is topologically compact but not complete is given by the following: Let $M\times \ex T$ be the product of a compact manifold with $\ex T$. This is complete, but we can construct an interesting non complete example as follows: choose an open subset $U\subset M$. Consider the subset of $M\times \ex T$ given by the union of $U\times \et 1{(-\infty,1)}$ with $M\times\et 1{(0,\infty)}$. The smooth part of this is the compact manifold $M$, but it is not complete. 

An example of an interesting exploded fibration which is complete but not basic is given by the quotient of $\mathbb R\times \ex T$ by $(x,\tilde z)=(x+1,\e 1\tilde z)$.

\

The following lemma is easy to prove:

\begin{lemma}

\begin{enumerate}
\item if $f$ and $g$ are complete, then $f\circ g$ is complete.

\item if $f\circ g$ is complete, then $g$ is complete.
\end{enumerate}
\end{lemma}

  \subsection{Tangent space}

\begin{defn}
A smooth exploded vectorfield $v$ on $\ex B$ is determined by maps
\[v:C^\infty(\ex B)\longrightarrow C^\infty(\ex B)\]  
   
  \[v:\mathcal E(\ex B)\longrightarrow \mathcal E(\ex B)\]
so that   
 \begin{enumerate}
   \item \[v(f+g)=v(f)+v(g)\]
   \item \label{derivation condition}\[v(fg)=v(f)g+fv(g)\]
   \item \[v(c\e yf)=c\e yv(f)\]
   \item The action of $v$ is compatible with the inclusion $\iota:C^{\infty}\longrightarrow \mathcal E$ in the sense that 
   \[v(\iota f)=\iota v(f)\]

   \end{enumerate}
   
      Smooth exploded vectorfields form a sheaf. The action of the restriction of $v$ to $U$ on the restriction of $f$ to $U$ is the restriction to $U$ of the action of $v$ on $f$. 
 
  We can restrict a vector field to a point $p\longrightarrow \ex B$ to obtain a tangent vector at that point. This is determined by  maps 
   \[v:C^\infty(\ex B)\longrightarrow \mathbb R\]
   \[v:\mathcal E(\ex B)\longrightarrow \mathbb C\e{\mathbb R}\]
    satisfying the above conditions  with condition \ref{derivation condition} replaced by
    \[v(fg)=v(f)g(p)+f(p)v(g)\]
Denote by $T_p\ex B$ the vector space of tangent vectors over $p\longrightarrow\ex B$.
  \end{defn}

  We can add smooth exploded vectorfields and multiply them by smooth real valued functions. We shall show below that the sheaf of smooth exploded vectorfields is equal to the sheaf of smooth sections of $T\ex B$, which is a real vector bundle over $\ex B$.

  \begin{lemma}\label{differentiation order}
  Differentiation does not change the order of a function in the sense that given any smooth exploded vectorfield $v$ and exploded function $f\in\mathcal E$, 
  \[\totb {f}=\totb {vf}\]
  \end{lemma}
  
  \pf
  \[v\e 0=v(\iota 1)=0\e 0\]
  Now we can apply $v$ to the equation
   $f=1\e 0f$, so $vf=0\e 0f+1\e 0vf$. Taking $\totb{\cdot}$ of this equation gives 
  \[\totb {vf}=\totb{ f}+\totb {vf}\text{ i.e. }\totb {vf}\geq\totb{f}\]
  (Recall that we use the order $\e x<\e y$ if $x>y$ as we are thinking of $\ex t$ as being tiny. So $\e x+\e y=\e x$ means that $\e x\geq\e y$.)
  Now, 
  \[v\e a=v(\e a\e 0)=\e a0\e 0=0\e a\]
  
   Now suppose that at some point $\totb {vf}>\totb{ f}$. Then we can restrict to a small neighborhood of this point, and choose an $a$ so that $\totb {vf}>\e a>\totb{f}$. But then $f+\e a=\e a$, but $\totb{v(f+\e a)}=\totb{vf}+\totb{\e a}=\totb{ vf}$, which is a contradiction, and the lemma is proved.
  
  \stop
  
  \begin{lemma}
  For any smooth exploded fibration $\ex B$, there exists a smooth exploded fibration $T\ex B$, the tangent space of $\ex B$, along with a canonical smooth projection $\pi: T\ex B\longrightarrow \ex B$ that makes  $T\ex B$ into a real vector bundle over $\ex B$. The sheaf of smooth exploded vectorfields is equal to the sheaf of smooth sections of this vector bundle.
  
  In particular, \[T(\mathbb R^n\times\et mA)=\mathbb R^{2n+2m}\times\et mA\]
  \end{lemma}
  
  \pf
  
  We shall first prove that $T(\mathbb R^n\times\et mA)=\mathbb R^{2n+2m}\times\et mA$. We shall use coordinate functions $x_i$ for $\mathbb R^n$ and $\tilde z_i$ for $\et mA$, as in example \ref{smooth local model} on page \pageref{smooth local model}.  A section of $\mathbb R^{2n+2m}\times\et mA\longrightarrow \mathbb R^n\times\et mA$ is given by $n+2m$ smooth functions on $\mathbb R^n\times \et mA$. To a vectorfield $v$, associate the $n$ smooth functions $v(x_i)$, and the $m$ smooth complex valued functions given by $\fun{ v(\tilde z_i)\tilde z_i^{-1}}$. (Lemma \ref{differentiation order} tells us 
 that $\totb{\tilde z_i^{-1}v(\tilde z_i)}=\e 0$, so this makes sense.) This gives us $n+2m$ smooth real valued functions, and therefore gives us a section to associate with our vectorfield.

Now we must show that the  arbitrary choice of $n+2m$ smooth functions on $\et mA$ to equal $v(x_i)$ and  $\fun{ v(\tilde z_i)\tilde z_i^{-1}}$ will uniquely determine a smooth exploded vectorfield.

First, recall that exploded functions are a sum of functions of the form $f\e y \tilde z^\alpha$, where $f$ is some smooth function of $x$ and $z^{\alpha^j}$. Use the notation
\[z^{\alpha^j}:=e^{t_{\alpha^j}+i\theta_{\alpha^j}}\]
There is the following formula for $v(f)$ which follows from the axioms:
\[v(f)=\sum v(x_i)\frac {\partial f}{\partial x_i}+\sum_{i,j}\Re \left(v(\tilde z_i)\tilde z_i^{-1}\right)\alpha_i^j\frac {\partial f}{\partial t_{\alpha^j}}+\Im\left(v(\tilde z_i)\tilde z_i^{-1}\right)\alpha_i^j\frac {\partial f}{\partial \theta_{\alpha^j}}\]

Note that this is well defined, despite the fact that $\frac {\partial f}{\partial t_{\alpha^j}}$ and $\frac {\partial f}{\partial \theta_{\alpha^j}}$ are not necessarily well defined. Note also that this is a smooth function, and is real valued if $f$ is real valued.  The corresponding formula for an exploded function is

\[v\left(\sum_{\alpha}f_\alpha\e {y_\alpha} \tilde z^\alpha\right):= \sum_{\alpha}\left(v(f_\alpha)+f_\alpha\sum_{i}\alpha_i v(\tilde z_i)\tilde z_i^{-1}\right)\e {y_\alpha}\tilde z^\alpha\]

It can be shown that $v$ satisfying such a formula satisfies all the axioms for being a smooth exploded vectorfield, is well defined, and is zero if and only if $v(x_i)=0$ and $\fun {v(\tilde z_i)\tilde z_i^{-1}}=0$.

 This shows that $T\mathbb R^n\times\et mA=\mathbb R^{2n+2m}\times \et mA$. The lemma then follows from our coordinate free definition of a smooth vectorfield, and the fact that every smooth exploded $\mathbb T$ fibration is locally modeled on exploded fibrations of this type.

\stop

For any smooth exploded $\mathbb T$ fibration $\ex B$, we now have that $T\ex B$ is a real vector bundle. In local coordinates $\{x_i,\tilde z_j\}$, a basis for this vectorbundle is given by $\{\frac \partial {\partial x_i}\}$ and the real and imaginary parts of $\{\tilde z_i\frac\partial{\partial\tilde z_i}\}$. The dual of this bundle $T^*\ex B$ is the cotangent space. A basis for this is locally given by $\{dx_i\}$ and the real and imaginary parts of $\{\tilde z_i^{-1}d\tilde z_i\}$. We can take tensor powers (over smooth real valued functions) of these vector bundles, to define the usual objects found on smooth manifolds. For example, a metric on $\ex B$ is a smooth, symmetric, positive definite section of $T^*\ex B\otimes T^*\ex B$.

\begin{defn}
 An integral vector $v$ at a point $p\longrightarrow\ex B$ is a vector $v\in T_p\ex B$ so that for any exploded function $f\in \mathcal E^\times(\ex B)$, 
 \[v(f)f^{-1}\in\mathbb Z\]
 Use the notation ${}^{\mathbb Z}T_p\ex B\subset T_p\ex B$ to denote the integral vectors at $p\longrightarrow \ex B$.
 
\end{defn}

For example, a basis for ${}^{\mathbb Z}T\ex T^n$ is given by the real parts of $ \tilde z_i\frac \partial{\partial \tilde z_i}$. The only integral vector on a smooth manifold is the zero vector.

\

Given a smooth morphism  $f:\ex B\longrightarrow \ex C$, there is a natural smooth morphism $df:T\ex B\longrightarrow T\ex C$ which is the differential of $f$, defined as usual by \[df(v)g:=v( g\circ f)\] Of course, $df$ takes integral vectors to integral vectors.

\subsection{$C^{k,\delta}$ regularity}\label{regularity}

In this section, we define the topology in which we shall prove our compactness theorem for holomorphic curves. The regularity $C^{\infty,\delta}$ is the level of regularity that the moduli stack of holomorphic curves can be expected to exhibit. A $C^{\infty,\delta}$ function is a generalization of a function on a manifold with a cylindrical end which is smooth on the interior, and which decays exponentially along with all its derivatives on the cylindrical end.  The reader wishing a simple introduction to holomorphic curves in \exploded fibrations may skip this somewhat technical section on first reading.

The following is a definition of a strata of an integral affine cone, similar to the definition of strata given in definition \ref{basic definition} on page \pageref{basic definition}.

\begin{defn}
The closure of a strata  of an integral cone $A:=\{\e a\ ;\ a\cdot\alpha^i\geq0\}$ is a subset of $A$ defined by the vanishing of $a\cdot\alpha$ for some $\alpha\in A^*$. A strata $S$ of $A$ is a subset which is equal to the closure of a strata minus the closure of all strata of smaller dimension. The zero substrata is the substrata on which $a\cdot\alpha=0$ for all $\alpha\in A^*$. 
 
\end{defn}

\begin{defn}
Given any $C^0$ function  $f$ defined on a strata $S\subset A$ in $\mathbb R^n\times \et mA$, define 
\[e_S(f)(x,\tilde z):=f(x,\tilde z\e a)\]
where $\e a:=(\e{a_1},\dotsc,\e{a_m})$ is any point  in $S$, and $\tilde z\e a$ means $(\tilde z_1\e {a_1},\dotsc,\tilde z_m\e {a_m})$.
%

\end{defn}

For example, $\et 22:=\et 2{[0,\infty)^{2}}$ has two one dimensional strata 
\[S_{1}:=\{\totb{\tilde z_{2}}=\e 0,\totb {\tilde z_{1}}\neq \e 0\}\ \ \ \ \ S_{2}:=\left\{\totb{\tilde z_{1}}=\e 0,\totb{\tilde z_{2}}\neq \e 0\right\}\]
 If we have a function $f\in C^0(\et 22)$, then

 \[e_{S_{1}}f(z_{1},z_{2})=f(0,z_{2})\ \ \ \ \ \ \ e_{S_{2}}f(z_{1},z_{2})=f(z_{1},0)\]
 Note that the operations $e_{S_{i}}$ commute and $e_{S_{i}}e_{S_{i}}=e_{S_{i}}$.
\begin{defn}If $I$ denotes any collection of strata $\{S_{1},\dotsc, S_{n}\}$, we shall use the notation 
\[e_{I}f:=e_{S_{1}}\lrb{e_{S_{2}}\lrb{\dotsb e_{S_{n}}f}}\]
\[\Delta_{I}f:=\lrb{\prod_{S_{i}\in I}(\id-e_{S_{i}})}f\]
\end{defn}
For example,

\[\begin{split}\Delta_{S_{1},S_{2}}f(z_{1},z_{2})&:=(1-e_{S_{1}})(1-e_{S_{2}})(f)(z_1,z_2)
\\&:=f(z_1,z_2)-f(0,z_2)-f(z_1,0)+f(0,0)\end{split}\]
Note that if $S\in I$, $e_{S}\Delta_{I}=0$. In the above example, this corresponds to $\Delta_{S_{1},S_{2}}f(z_{1},0)=0$ and $\Delta_{S_{1},S_{2}}f(0,z_{2})=0$.

\

We shall need a weight function $w_{I}$ for every collection of nonzero strata $I$. This will have the property that if $f$ is any smooth function $w_{I}^{-1}\Delta_{I}f$ will be bounded on any topologically compact subset of $\mathbb R^{n}\times \et mA$ (in other words, a subset with a compact image in the topological part $\totl{\mathbb R^{n}\times \et mA}$). Consider the ideal of functions of the form $\tilde z^{\alpha}$ so that $\Delta_{I}\totl{z^{\alpha}}=\totl{z^{\alpha}}$ (In other words, $e_{S}\totl{z^{\alpha}}=0$ for all $S\in I$.) Choose some finite set of generators $\{\tilde z^{\alpha^{i}}\}$ for this ideal. Then define 

\[w_{I}:=\sum\abs{\totl{\tilde z^{\alpha^{i}}}}\]

Continuing the example started above, we can choose  $w_{S_{1}}=\abs{z_{1}}$, $w_{S_{2}}=\abs{z_{2}}$,  
$w_{S_{1},S_{2}}=\abs{z_{1}z_{2}}$, and $w_{S_{1}\cap S_{2}}=\abs {z_{1}}+\abs{z_{2}}$.

Note that given any other choice of generators, the resulting $w'_{I}$ is bounded by a constant times $w_{I}$ on any topologically compact subset of $\mathbb R^{n}\times \et mA$. Note also that $w_{I}w_{I'}$ is bounded by a constant times $w_{I\cup I'}$ on any topologically compact subset of $\mathbb R^{n}\times \et mA$. 

\

We shall now define $C^{k,\delta}$ for any $0<\delta<1$:

\begin{defn}
Define $C^{0,\delta}$ to be the same as $C^0$. A sequence of smooth functions $f_i\in C^\infty(\mathbb R^n\times\et mA)$ converge to a continuous function $f$ in  $C^{k,\delta}(\mathbb R^n\times \et mA)$ if the following conditions hold:
\begin{enumerate}
\item Given any collection of at most $k$ nonzero strata $I$, 
\[\abs{w_{I}^{-\delta}\Delta_{I}(f_i-f)}\]
converges to $0$ uniformly on compact subsets of $\totl{\mathbb R^{n}\times\et mA}$ as $i\rightarrow\infty$. (This includes the case where our collection of strata is empty and $f_{i}\rightarrow f$ uniformly on compact subsets.)
\item For any smooth exploded vectorfield $v$, $v(f_i)$ converges to some function $vf$ in $C^{k-1,\delta}$.
\end{enumerate} 

Define $C^{k,\delta}(\mathbb R^n\times \et mA)$ to be the closure of $C^\infty$ in $C^0$ with this topology. Define $C^{\infty,\delta}$ to be the intersection of $C^{k,\delta}$ for all $k$.
\end{defn}

In particular, $C^k\subset C^{k,\delta}$ for $0<\delta< 1$. Functions in $C^{k,\delta}$ can be thought of as functions which converge a little slower than $C^k$ functions when they approach different strata. Thinking of a single strata as being analogous to a cylindrical end, this is similar to requiring exponential convergence (with exponent $\delta$) on the cylindrical end.  We shall now start showing that we can replace smooth functions with $C^{\infty,\delta}$ functions in the definition of exploded torus fibrations to create a category of $C^{\infty,\delta}$ exploded torus fibrations.

\begin{defn}
A $ C^{k,\delta}$ exploded function $f\in \mathcal E^{k,\delta,\times}(\mathbb R^{n}\times\et mA)$ is a function of the form 
\[f(x,\tilde z):= g(x,\tilde z)\tilde z^{\alpha}\e a\text{ where }g\in C^{k,\delta}\lrb{\mathbb R^{n}\times\et mA,\mathbb C^{*}},\ \alpha\in\mathbb Z^{m},\ \e a\in\e{\mathbb R}\]

A sequence of exploded functions $g^{i}\tilde z^{\alpha}\e {a_{i}}$ converge in $C^{\infty,\delta}$ if the sequence of functions $g^{i}$ does, and the sequence $a_{i}$ is eventually constant.

\

A $C^{k,\delta}$ exploded fibration is an  abstract exploded fibration locally modeled on open subsets of $\mathbb R^{n}\times \et mA$ with the sheaf $\mathcal E^{k,\delta,\times}$.
\end{defn}

\

The following observations are easy to prove:
\begin{lemma} $C^{k,\delta}$ is an algebra over $C^{\infty}$ for any $0<\delta<1$.

\end{lemma}

\begin{lemma}\label{linear change}

Given any `linear' map 
\[\alpha:\mathbb R^{n}\times\et mA\longrightarrow \mathbb R^{n'}\times\et {m'}B\]
\[\alpha (x,\tilde z):=\lrb{Mx,\tilde z^{\alpha^{1}},\dotsc,z^{\alpha^{m'}}}\] 
where $M$ is a $n$ by $n'$ matrix and $\alpha^{i}_{j}$ is a $m$ by $m'$ matrix with $\mathbb Z$ entries, $\alpha$ preserves $C^{\infty,\delta}$ in the sense that given any function $f\in C^{\infty,\delta}(\mathbb R^{n'}\times\et {m'}B)$,  
\[\alpha\circ f\in C^{\infty,\delta}\lrb{\mathbb R^{n}\times\et mA}\]
\end{lemma}

More difficult is the following:

\begin{lemma}\label{exponentiate vectorfield}
 Given any $C^{k,\delta}$ section of $T\lrb{\mathbb R^{n}\times \et mA}$, we can define a map of the form 
\[\exp(f)(x,\tilde z_{1},\dotsc,\tilde z_{m}):=\lrb{x+f_{\mathbb R^{n}}(x,\tilde z), e^{f_{1}(x,\tilde z)}\tilde z_{1},\dotsc,e^{f_{m}(x,\tilde z)}\tilde z_{m}}\]
 If $h$ is in $C^{k,\delta}$, then $h\circ\exp f$ is.

\end{lemma}

\pf

We must show that expressions of the following form decay appropriately.
\[\Delta_{I}(h\circ\exp f)\]
Note that 
\[e_{S}(h\circ\exp f)=(e_{S}h)\circ\exp (e_{S}f) \]
A little thought shows that we can rewrite our above expression in the following form
\[\Delta_{I}(h\circ\exp f)=\sum_{I'\cup I''=I, I'\cap I''=\emptyset }\lrb{\lrb{ e_{I'}\Delta_{I''}h}\circ\exp\circ \Delta_{I'}}f  \]

As an example for interpreting the notation above, we write
\[\lrb{h\circ\exp\circ \Delta_{S}}f:=h\circ\exp f-h\circ\exp (e_{S}f)\]
as opposed to 
\[h\circ\exp\lrb{\Delta_{S}f}:=h\circ\exp(f-e_{S}f)\]

We have that because $f$ is bounded on compact sets, 
\[(\Delta_{I}h)\circ f\text{ decays appropriately.}\]
We can bound expressions of the form
\[\lrb{\lrb{\Delta_{I'}h}\circ \exp \circ \Delta_{I''}}f\]
by the size of the first $\abs{I''}$ derivatives of $\Delta_{I'} h$ times the sum of all products of the form $\prod \abs {\Delta_{I_{i}}f}$ where $I''$ is the disjoint union of $\{I_{i}\}$. We therefore get that all of the terms in the above expression decay appropriately. The appropriate decay of the derivatives of $h\circ \exp f$ follows from the fact that $C^{\infty,\delta}$ is closed under multiplication, and a similar argument to the above one.

\stop

Consider a map for which the pull back of exploded coordinate functions is in $\mathcal E^{k,\delta,\times}$ and the pullback of real coordinate functions is $C^{k,\delta}$. Any map of this form factors as a composition of maps in the form of Lemma \ref{linear change} and Lemma \ref{exponentiate vectorfield}; (first a map in the form of Lemma \ref{linear change} to the product of the domain and target, then a map of the form of Lemma \ref{exponentiate vectorfield}, then a projection to the target, which is of the form of Lemma \ref{linear change}.) We therefore have the following: 

\begin{cor}
A morphism is $ C^{k,\delta}$ if and only if the  pull back of exploded coordinate functions are $\mathcal E^{k,\delta,\times}$ functions, and the pull back of real coordinate functions are $ C^{k,\delta}$ functions.

\end{cor}

We shall need the following definition of convergence.

\begin{defn}
A sequence of $C^{k,\delta}$ exploded maps 
$f^i:\ex A\longrightarrow \ex B$ converges to $f:\ex A\longrightarrow \ex B$ in $C^{k,\delta}$ if  the pullback  under $f^{i}$ of any local coordinate function on $\ex B$ converges in $C^{k,\delta}$ to the pullback under $f$.

\end{defn}

\subsection{Almost complex structures}

\begin{defn}
An almost complex structure $J$ on a smooth exploded $\mathbb T$ fibration $\ex B$ is an endomorphism of $T\ex B$ given by a smooth section of $ T\ex B\otimes T^*\ex B$ which squares to become multiplication by $-1$, so that for any coordinate chart  modeled on an open subset of $\mathbb R^{n}\times\ex T^{m}$   $J$ of the real  part of $\tilde z_{i}\frac \partial {\partial\tilde z_{i}}$ is the imaginary part. 

An almost complex structure $J$ is a complex structure if there exist local coordinates $\tilde z\in\et mA$ so that for all vectorfields $v$, $iv(\tilde z_{j})=(Jv)(\tilde z_{j})$. 
\end{defn}

This differs from the usual definition of an almost complex structure only in the extra requirement that in standard coordinates, $(Jv)\tilde z_{j}=iv\tilde z_{j}$ which is only required to hold when $\tilde z_{j}$ is locally a coordinate on $\ex T$. For instance, on $\et 11$, this requirement is only valid when $\totb{\tilde z}<\e 0$. This extra requirement makes holomorphic curves $C^{\infty,\delta}$ exploded $\mathbb T$ morphisms, and makes the tropical part of holomorphic curves piecewise linear one complexes. If it did not hold, then we would need to use a different version of the exploded category using $\mathbb R^{*}$ instead of $\mathbb C^{*}$ to explore holomorphic curve theory. The analysis involved would be significantly more difficult. 

The following assumption allows us to use standard (pseudo)holomorphic curve results.

\begin{defn}\label{civilized}
An almost complex structure $J$ on $\ex B$ is civilized if it induces a smooth almost complex structure on the smooth part of $\ex B$. This means  that given any local coordinate chart modeled on $\et mA$, if we consider $\et mA$ as a subset of $\et {m'}n$ defined by setting some monomials equal to $1$, there exists an almost complex structure $\hat J$ on $\et {m'}n$ which induces a  smooth almost complex structure on $\totl{\et {m'}n}=\mathbb C^{n}$ so that the subset corresponding to $\et mA$ is holomorphic, and the restriction of $\hat J$ is $J$.
\end{defn}

The word civilized should suggest that our almost complex structure is well behaved in a slightly unnatural way. Any complex structure is automatically civilized, and there are no obstructions to modifying an almost complex structure to civilize it.  If we assume our almost complex structure is civilized, then holomorphic curves are smooth maps. Otherwise, they will just be $C^{\infty,\delta}$ for any $\delta<1$. 

\begin{example}\label{expl}\end{example}

Suppose that we have a complex manifold $M$ along with a collection of immersed complex submanifolds $N_i$ so that $N_i$ intersect themselves and each other transversely. Then there is a smooth complex exploded fibration $\expl (M)$ called the explosion of $M$. (There is actually a functor $\expl$ from complex manifolds with a certain type of log structure which is defined by the submanifolds $N_i$ to the category of exploded fibrations.) 

We define $\expl(M)$ as follows: choose holomorphic coordinate charts on  $M$ which are equal to balls inside $\mathbb C^n$, so that the image of the submanifolds $N_i$ are equal to the submanifolds $\{z_i=0\}$. Then replace a coordinate chart $U\subset \mathbb C^n$ by a coordinate chart $\tilde U$ in $\et nn$ with coordinates $\tilde z_i$ so that 
\[\tilde U:=\{\tilde z\text{ so that }\totl{\tilde z}\in U\}\]
 Define transition functions as follows: the old transition functions are all of the form $f_i(z)=z_jg(z)$ where $g$ is holomorphic  and non vanishing. Replace this with $\tilde f_i(\tilde z)=\tilde z_jg(\totl{\tilde z})$, which is then a smooth exploded function.  This defines transition functions which define $\expl (M)$ as a smooth complex exploded fibration.

 \psfrag{ETF4manifold}{$\totl{\expl M}=M$}
 \psfrag{ETF4totb}{$\totb{\expl M}$}
\includegraphics{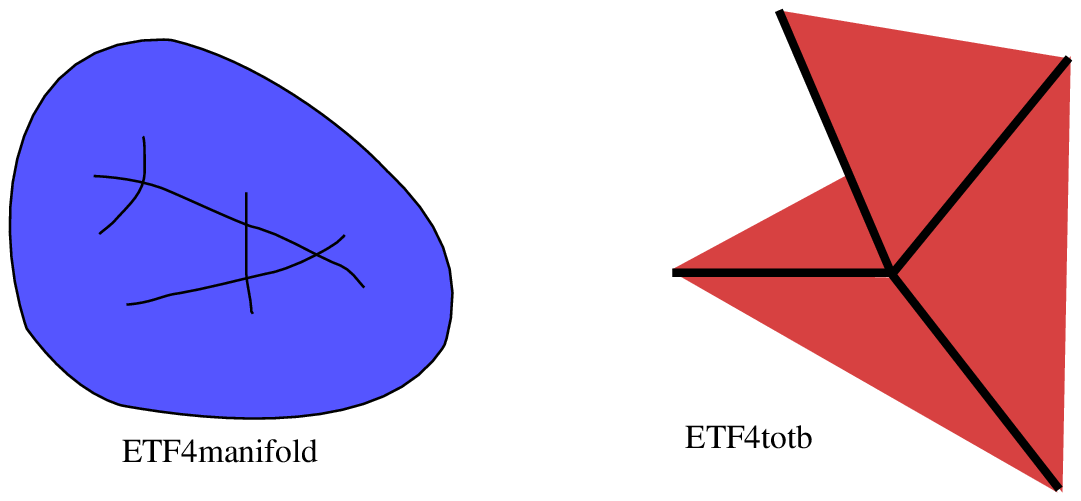}

\begin{defn}
An \exploded  curve is a $2$ real dimensional, complete \exploded fibration with a complex structure $j$. 

A holomorphic curve is a holomorphic map of an \exploded curve to an almost complex exploded fibration.

A smooth exploded curve is a smooth map of an \exploded curve to a exploded $\mathbb T$ fibration. 
\end{defn}

By a smooth component of a holomorphic curve $\ex C$, we shall mean a connected component of the set of points in $\ex C$ which have a neighborhood modeled on an open subset of $\mathbb C$. By an edge of $\ex C$ we shall mean connected component in $\tot{\ex C}$ minus the image of all smooth components. We shall call edges that have only one end attached to a smooth component `punctures'. Each smooth component is a punctured Riemann surface with punctures corresponding to where it is connected to edges. The information in a holomorphic curve $\ex C$ is equal to the information of a nodal Riemann surface plus gluing information for each node parametrized by $\mathbb C^{*}\e{(0,\infty)}$ (for more details of this gluing information see example \ref{local family model} on page \pageref{local family model}).

\psfrag{ETF7v1}{$v^{1}$}\psfrag{ETF7v2}{$v^{2}$}\psfrag{ETF7v3}{$v^{3}$}\psfrag{ETF7v4}{$v^{4}$}
\psfrag{ETF7tot}{$\tot{\ex C}$}
\psfrag{ETF7totb}{$\totb{f(\ex C)}\subset\totb{\ex T^{n}}$}
\psfrag{ETF7eqn}{$\sum df(v^{i})=0$}
\includegraphics{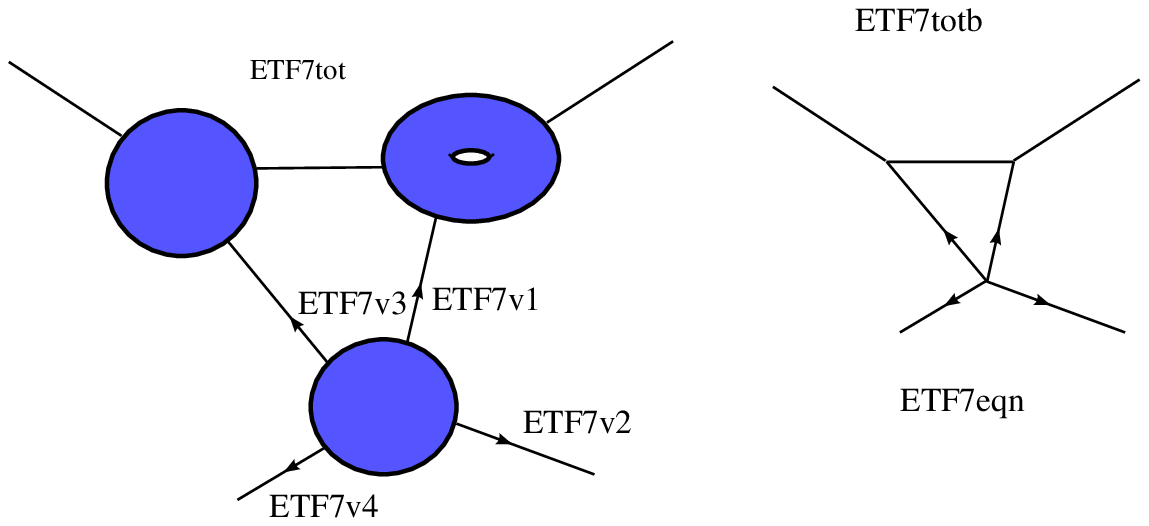}

\begin{example}\label{curve in tropical plane}
\end{example}
Consider a smooth curve $f:\ex C\longrightarrow \ex T^{n}$. This is given by $n$ exploded functions $f^{*}(\tilde z_{1}),\dotsc,f^{*}(\tilde z_{n})\in\mathcal E^{\times}(\ex C)$. Each smooth component of $\ex C$ is sent to the $\lrb{\mathbb C^{*}}^{n}$ worth of points over a particular point in the tropical part $\totb{\ex T^{n}}$. We can choose local holomorphic coordinates $\tilde w\in\et 11$ on $\ex C$. In these coordinates, 
\[f(\tilde w)=(g_{1}(\totl{\tilde w})\e {a_{1}}\tilde w^{\alpha_{1}},\dotsc,g_{n}(\totl{\tilde w})\e{a_{n}}\tilde w^{\alpha_{n}})\ \ g_{i}\in C^{\infty}(\mathbb C,\mathbb C^{*}),\ \alpha\in\mathbb Z^{n}\]
 
 (Our smooth curve is holomorphic if and only if the functions $g_{i}$ in our local coordinate representation above are holomorphic.)
In particular, $f$ restricted to a smooth component gives a smooth map of the corresponding punctured Riemann surface to $\lrb{\mathbb C^{*}}^{n}$,  and the homology class in $H_{1}\lrb{\lrb{\mathbb C^{*}}^{n}}$ of  a loop around the puncture corresponding in the above coordinates to $\totl{\tilde w}=0$ is given by $\alpha$. Of course, the sum of all such homology classes from punctures of a smooth component is zero. As $\alpha$ also determines the derivative of $f$ on edges, this can be viewed as some kind of conservation of momentum condition for the tropical part of our curve, $\totb{f}$. In tropical geometry, this is called the balancing condition.

\begin{example}\label{tropical curve}
\end{example}
One way to consider the image of some holomorphic curves  in $\ex T^{n}$ is as the `locus of noninvertablity' of some set of polynomials 
\[P_{i}(\tilde z):=\sum \tilde c_{i,\alpha}\tilde z^{\alpha}\ i=1,\dotsc,n-1\]
We can consider the set
\[ Z_{\{P_{i}\}}:=\left\{\tilde z\text{ so that }P_{i}(\tilde z)\in 0\e{\mathbb R}\ \forall i\right\}\]
Suppose that for all points $p\longrightarrow Z_{\{P_{i}\}}$, the differentials $\{dP_{i}\}$ at $p$ are linearly independent. Then $Z_{\{P_{i}\}}$ is the image of some holomorphic curve.

Let us examine the set $Z_{\{P_{i}\}}$ more closely. For any point $\tilde z_{0}$, denote by $S_{i,\tilde z_{0}}$ the set of exponents $\alpha$ so that $\totb{P_{i}(\tilde z_{0})}=\totb{\tilde c_{i,\alpha}\tilde z^{\alpha}}$. Then there exists some neighborhood of $\totb{\tilde z_{0}}$ in $\totb{\ex T^{n}}$ so that 
\[P_{i}=\sum_{\alpha\in S_{i,\tilde z_{0}}}\tilde c_{i,\alpha}\tilde z^{\alpha}\]
The points inside $Z_{\{P_{i}\}}$ over $\totb{\tilde z_{0}}$ are then given by solutions of the equations 
\[\sum_{\alpha\in S_{i,\tilde z_{0}}} c_{i,\alpha}z^{\alpha}=0\text{ where } \tilde c_{i,\alpha}=c_{i,\alpha}\totb{\tilde c_{i,\alpha}}\text{ and }\tilde z=z\totb{\tilde z}\]
Note that the above equation has solutions for $z\in (\mathbb C^{*})^{n}$ if and only if $S_{i,\tilde z_{0}}$ has more than $1$ element. This corresponds to the tropical function $\totb{P_{i}}$ (which is continuous, piecewise integral affine, and convex) not being smooth at $\totb{\tilde z_{0}}$. We therefore have that $\totb{Z_{\{P_{i}\}}}$ is contained in the intersection of the non-smooth locus of the tropical polynomials $\totb{P_{i}}$.

\subsection{Fiber product, refinements}

\begin{defn}
Two smooth (or $C^{k,\delta}$) exploded morphisms 
\[\ex A\xrightarrow{f}\ex C\xleftarrow{g}\ex B\]
are transverse if for every pair of points $p_1\longrightarrow\ex A$ and $p_2\longrightarrow \ex B$ so that $f(p_1)=g(p_2)$,  $df(T_{p_1}\ex A)$ and $dg(T_{p_2}\ex B)$ span $T_{f(p_1)}\ex C$.
\end{defn}

\begin{defn}
If $f$ and $g$ are transverse smooth (or $C^{k,\delta}$) exploded morphisms,
\[\ex A\xrightarrow{f}\ex C\xleftarrow{g}\ex B\]
The fiber product $\ex A\fp fg\ex B$ is the unique smooth (or $C^{k,\delta}$) exploded $\mathbb T$ fibration   with maps to $\ex A$ and $\ex B$ so that the following diagram commutes
\begin{displaymath}
\begin{array}{cll}
\ex A\fp fg\ex B & \longrightarrow & \ex A \\ 
 \downarrow & × & \downarrow \\ 
\ex B & \longrightarrow & \ex C
\end{array}
\end{displaymath}
and with the usual universal property that given any commutative diagram
\begin{displaymath}
\begin{array}{lll}
\ex D & \longrightarrow  & \ex A \\ 
\downarrow & × & \downarrow \\ 
\ex B & \longrightarrow & \ex C
\end{array}
\end{displaymath}
there exists a unique morphism $\ex D\longrightarrow \ex A\fp fg\ex B$ so that the following diagram commutes
\begin{displaymath}
\begin{array}{llc}
\ex D & \rightarrow & \ex A \\ 
\downarrow & \searrow & \uparrow \\ 
\ex B & \leftarrow & \ex A\fp fg\ex B
\end{array}
\end{displaymath}
\end{defn}

The existence of the fiber products defined above is left as an exercise. The existence in some special cases was proved in \cite{gokova}.

\begin{example}
 
\end{example}
 Consider the map $f:\et mA\longrightarrow\ex T^n$
given by 
\[f(\tilde z)=(\tilde z^{\alpha^1},\dotsc,\tilde z^{\alpha^n})\]
Denote by $\alpha$ the $m\times n$ matrix with entrys $\alpha^i_j$. The derivative of $f$ is surjective if $\alpha:\mathbb R^m\longrightarrow\mathbb  R^n$ is. Denote by $\abs\alpha\in\mathbb N$ the size of $\alpha^1\wedge\dotsb\wedge\alpha^n\in \bigwedge^n(\mathbb Z^m)$. (In other words, $\abs{\alpha}\in\mathbb N$  is the smallest  number so that the above wedge is at most $\abs\alpha$ times any nonzero element of $\bigwedge^{n}(\mathbb Z^{m})$.) The fiber product of $f$ with the point $(1,\dotsc,1)$ corresponds to the points in $\et mA$ so that $\tilde z^{\alpha^i}=1$ for all $i$. This is then equal to $\abs \alpha$ copies of $\et {m-n}{\e{\ker \alpha}\cap A}$ where we identify $\e{\mathbb R^{m-n}}=\e{\ker\alpha}$.

\

The following should not be difficult to prove.

\begin{conj}
If $f$ and $g$ are complete and transverse
\[\ex A\xrightarrow{f}\ex C\xleftarrow{g}\ex B\]
then the fiber product $\ex A\fp fg\ex B\longrightarrow \ex C$ is complete.
\end{conj}

\begin{defn}\label{family defn}
  A family of exploded $\mathbb T$  fibrations over $\ex F$ is a map $f:\ex C\longrightarrow \ex F$ so that:
 \begin{enumerate}
  \item $f$ is complete
  \item for every point $p\longrightarrow \ex C$, 
  \[df:T_p\ex C\longrightarrow T_{f(p)}\ex F\text{ is surjective }\]
 \[\text{ and }df:{}^{\mathbb Z}T_p\ex C\longrightarrow {}^{\mathbb Z}T_p\ex F\text{ is surjective }\]
 \end{enumerate}

\end{defn}

(The definition of complete is found on page \pageref{complete}. Recall also that integer vectors in ${}^{\mathbb Z}T_p\ex C$ are the vectors $v$ so that for any exploded function, $vf$ is an integer times $f$. For example on $\ex T$ the integer vectors are given by integer multiples of the real part of $\tilde z\frac \partial{\partial \tilde z}$, so the map $\ex T\longrightarrow\ex T$ given by $\tilde z^{2}$ is not a family as it does not obey the last condition above.)

The local normal form for coordinate charts on a smooth family is a map $\et mA\longrightarrow\et {k}B$ given by $(\tilde z,\tilde w)\mapsto \tilde z$,  where  
the cone $B$ is given by the projection of $A$ to the first $k$ coordinates, and the projection of every strata of $ A$ is a strata of $ B$. This may differ from a product because the cone $A$ may not be a product of $B$ with something.

\begin{example}
\end{example}
We can represent the usual compactified moduli space  of stable curves $\bar{\mathcal M}_{g,n}$ as a complex orbifold. There exist local holomorphic coordinates so that the boundary of $\bar{\mathcal M}_{g,n}$ in these coordinates looks like $\{z_{i}=0\}$. As in example \ref{expl} on page \pageref{expl}, we can replace these coordinates $z_{i}$ with $\tilde z_{i}$ to obtain a complex orbifold exploded fibration $\expl \left(\bar{\mathcal M}_{g,n}\right)$. The forgetful map $\pi:\bar{\mathcal M}_{g,n+1}\longrightarrow\bar{\mathcal M}_{g,n}$ induces an map 
\[\pi:\expl\left(\bar{\mathcal M}_{g,n+1}\right)\longrightarrow\expl\left(\bar{\mathcal M}_{g,n}\right)\]
This is a family, and each stable exploded curve with genus $g$ and $n$ marked points corresponds to the fiber over some point $p\longrightarrow\expl\left(\bar{\mathcal M}_{g,n}\right)$. This example is considered in a little more detail in \cite{gokova}. (Obviously, it would take a little more work to prove that $\expl\bar{\mathcal M}_{g,n}$ actually represents the moduli stack of stable exploded curves, but that is not the subject of this paper.)

\begin{example}\label{local family model}
\end{example}
The following example contains all interesting local behavior of the above example. It is not a family only because it fails to be complete (it is not topologically proper). Consider the map 
\[\pi:\et 22\longrightarrow\et 11\text{ given by }\pi^{*}\tilde z=\tilde w_{1}\tilde w_{2}\]
The derivative is surjective, as can be seen by the equation 
\[\pi^{*}(\tilde z^{-1}d\tilde z)=\tilde w_{1}^{-1}d\tilde w_{1}+\tilde w_{2}^{-1}d\tilde w_{2}\]

The fibers of this map over smooth points $\tilde z=c\e 0$ are smooth manifolds equal to $\mathbb C^{*}$ considered just as a smooth manifold with coordinates $w_{1}$ and $w_{2}\in\mathbb C^{*}$ related by $w_{1}w_{2}=c$. (Note that there is no point with $\tilde z=0\e0$.)

In contrast, the fibers of this map over points $\tilde z=c\e x$ with $x>0$ are not smooth manifolds. They can be described using coordinates $\tilde w_{1},\tilde w_{2}\in\et 1{[0,x)}\subset \et 11$, and transition maps give by $\tilde w_{2}\tilde w_{2}=c\e x$ on their intersection, which is $\et 1{(0,x)}$. The picture below shows some of the fibers of this map

\includegraphics{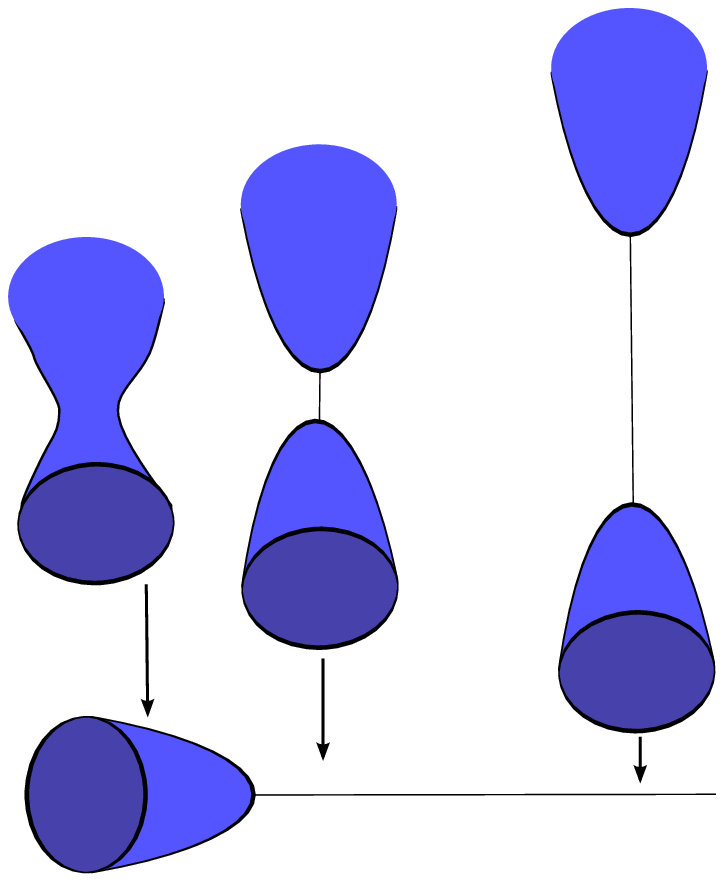}

     \begin{defn}
      A refinement of $\ex B$ is an exploded $\mathbb T$ fibration $\ex B'$ with a map $f:\ex B'\longrightarrow \ex B$  so that
      \begin{enumerate}
 \item $f$ is complete
 \item $f$ gives a bijection between points in $\ex B'$ and $\ex B$
 \item $df$ is surjective
\end{enumerate}

     \end{defn}

\begin{example}\end{example}\label{refinement example}
All refinements are locally of the following form: Suppose that the integral affine cone $A\subset\e{\mathbb R^m}$ is equal to the union of a collection of  cones $\{\hat A_i\}$ so that the intersection of any two of these $\hat A_i$ has codimension at least $1$. Then using the standard coordinates $(x,\tilde z)$ from example \ref{smooth local model} (on page \pageref{smooth local model}) on $\mathbb R^{n}\times \et m{\hat A_{i}}$ for each $\hat A_i$, we can piece the coordinate charts $\mathbb R^n\times \et m{\hat A_i}$ together using the identity map as transition coordinates. In these coordinates, the map down to $\mathbb R^n\times\et mA$ is again the the identity map in these coordinates. (The fact that all refinements are locally of this form can be proved using the observation that the pull back of coordinate functions on $\ex B$ must locally be coordinate functions on $\ex B'$.) So in the case that $\ex B$ is basic, a refinement is simply determined by a subdivision of the tropical part $\totb{\ex B}$. The effect on the smooth part $\totl{\ex B}$ should remind the reader of a toric blowup.

\

The following lemma follows from the above standard local form for refinements. 

\begin{lemma}
 Given a refinement of $ \ex B$, 
 \[f:\ex A\longrightarrow \ex B\]
 Any smooth exploded vector field  $v$ on $\ex B$ lifts uniquely to a smooth exploded vector field $\tilde v$ on $\ex A$ so that $df(\tilde v)=v$.
\end{lemma}

The above lemma tells us that any smooth tensor field (such as an almost complex structure or metric) lifts uniquely to a smooth tensor field on any refinement. The above normal form implies that given any morphism $\ex C\longrightarrow\ex B$  and a refinement $\ex B'\longrightarrow\ex B$, there exists a refinement $\ex C'\longrightarrow\ex C$ and a morphism $\ex C'\longrightarrow\ex B'$ lifting the above map.

\begin{defn} Call a holomorphic curve stable if it has a finite number of automorphisms, and is not a nontrivial refinement of another holomorphic curve. 
\end{defn}

If $\ex B$ has an almost complex structure,  there is a bijection between stable holomorphic curves in $\ex B$ and stable holomorphic curves in any refinement $\ex B'$. In fact, when the moduli space of stable holomorphic curves in $\ex B$ is smooth, the moduli space of stable holomorphic curves in $\ex B'$ is a refinement of this moduli space.


\subsection{Symplectic structures} 

 Our definition of symplectic structures below should be regarded as provisional. Our definitions are enough  to tame holomorphic curves, which is all that is required for this paper, but there is probably a better theory of symplectic structures on exploded fibrations. 
 
\begin{defn}
A symplectic exploded $\mathbb T$ fibration $(\ex B,\omega)$ is an even dimensional exploded $\mathbb T$ fibration $\ex B$ with a two form $\omega$ so that each point has a neighborhood equal to an open set in $\mathbb R^{2n}\times\et mA$ with a two form defined as follows
\[\omega:=\omega_{\mathbb R^{2n}} -i\sum_{j} d{z^{\alpha^j}}\wedge d \bar z^{\alpha^j}\]  
Here $\omega_{\mathbb R^{2n}}$ is the standard form on $\mathbb R^{2n}$ and
$\{\alpha^j\}$ is some finite set generating $A^*$. 

\end{defn}

The second part of this symplectic form can be written as 
\[\sum_i d(\alpha_id\theta_i)\]
\[\alpha_i=\sum_j\alpha_i^j\abs{z^{\alpha^j}}^2\text{ or }\alpha:=\sum_j\abs{z^{\alpha^j}}^2\alpha^j\]

By $d\theta_i$, we mean the imaginary part of $\tilde z^{-1}d\tilde z$. The Hamiltonian torus action on this standard model given by the vectorfields $\frac\partial{\partial \theta_i}$ has moment map given by $\alpha:=(\alpha_1,\dotsc,\alpha_n)\subset \mathbb R^n$. The image of this moment map is the  dual cone to $A$ defined by 
\[\{\alpha\text{ so that }\alpha\cdot a\geq 0 \text{ for all }a\in A\}\] 



%
%


\begin{example}\label{relative example}
\end{example} 
Suppose that we have a compact symplectic manifold $M$ with some collection of embedded symplectic submanifolds $N_i$ which intersect each other orthogonally. (In other words, if $x\in N_i\cap N_j$, the vectors $v\in T_xM$ so that $\omega(v,\cdot)$ vanishes on  $T_xN_i$ are contained inside $T_xN_j$.) Then a standard Moser type argument gives that there exists some collection of neighborhoods $U_i$ of $N_i$ so that:
\begin{enumerate}
\item $U_i$ is identified with some neighborhood of the zero section in a complex line bundle over $N_i$ (with $\mathbb T$ action given by multiplication by $e^{i\theta}$), and a $\mathbb T$ invariant connection one form $\alpha_i$ so that $\omega=\omega_{N_i}+d(r_i^2\alpha_i)$, where $r_i$ indicates some radial coordinate on the line bundle, $N_i$ is identified with the zero section, and $\omega_{N_i}$ indicates the pullback of the form $\omega$ under the projection to $N_i$.
\item These identifications are compatible in the sense that on $\bigcap_{i\in I}U_i$, all the different projections to $N_i$ and the $\mathbb T$ actions commute, and 
\[\omega=\omega_{\bigcap N_i}+\sum d(r_i^2\alpha_i)\] 
where $\omega_{\bigcap N_i}$ is the pullback of $\omega$ under the composition of projection maps to $\bigcap N_i$.
\end{enumerate}
We can associate an exploded $\mathbb T$ fibration to the above as follows: 
\begin{enumerate}
\item Choose coordinate charts $V$ on $M$ which are identified with open balls in $\mathbb C^n$ so that the intersection of $V$ with the submanifolds $N_{i}$ is equal to the union of the planes $z_{k}=0$. If $N_i$ intersects this coordinate chart, it does so where  some coordinate $z_{k_i}=0$. Choose these charts so that the fibers of the complex line bundle over $N_i$ are identified with the slices where all other coordinates are constant, and the $\mathbb T$ action is multiplication of $z_{k_i}$ by $e^{i\theta}$, and $r_i=\abs{ z_{k_i}}$.
\item Replace $V\subset\mathbb C^n$ with the corresponding subset of $\et nn$, where each coordinate $z_i$ is replaced by the standard coordinate  $\tilde z_i$. Replacing $z_i$ with $\tilde z_i$ in transition functions between these coordinates gives transition functions which are smooth exploded isomorphisms, which defines an exploded fibration $\ex M$. A standard Moser type argument then gives that the two form on $\ex M$ corresponding to $\omega$ is in the form defined above.
\end{enumerate}

The following is a picture of the result of this procedure on a toric symplectic manifold, where for our symplectic submanifolds we use the fixed loci of circle actions. We shall represent the smooth part $\totl{\ex M}$ by its moment map:

\psfrag{ETF6M}{$\totl{\ex M}$}
\psfrag{ETF6tot}{$\tot{\ex M}$}
\psfrag{ETF6totb}{$\totb{\ex M}$}
\includegraphics{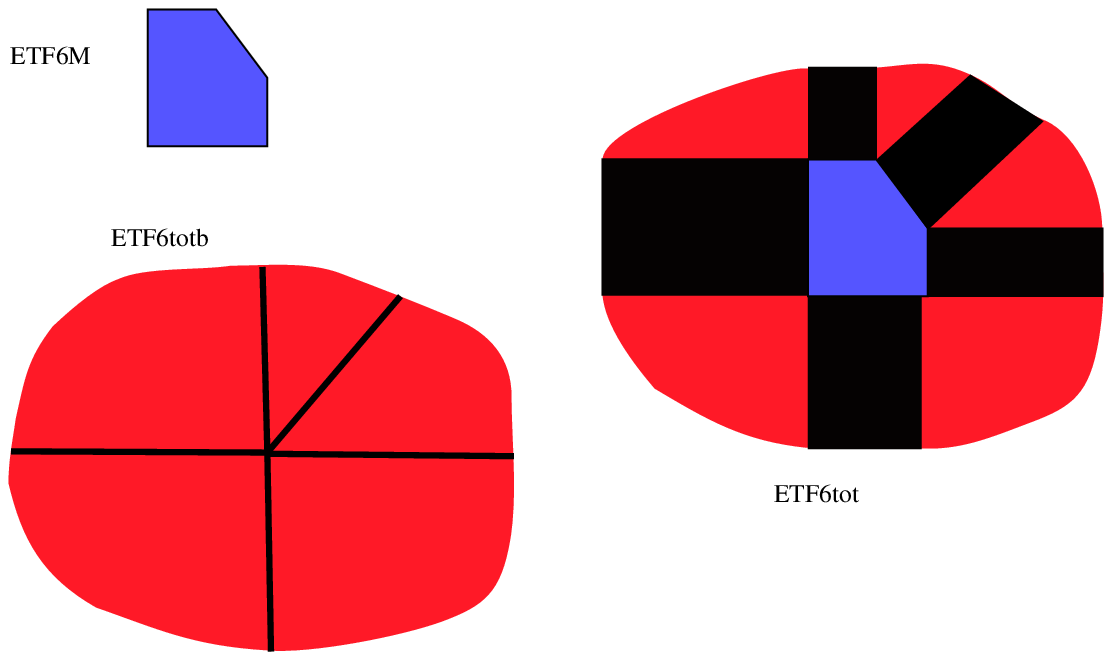}

\begin{defn} An almost complex structure $J$ on a symplectic exploded $\mathbb T$ fibration $(\ex B,\omega)$ is tamed by $\omega$ if $\omega(v,Jv)$ is positive for any vector $v$ so that $vf\neq 0$ for some smooth function $f$.
 \end{defn}

The above definition of taming is clearly inadequate for taming holomorphic curves in $\ex T^{n}$ where it is satisfied trivially by a two form that is identically $0$ because all smooth functions on $\ex T^{n}$ are locally constant. We remedy this as follows: (Note that we shall use freely the fact that any holomorphic curve in $\ex B$ lifts to a holomorphic curve (with a refined domain) in $\ex B'$ for any refinement $\ex B'\longrightarrow\ex B$.)

 \begin{defn}
  A strict taming of a complete almost complex exploded $\mathbb T$ fibration $\ex B$ is a set $\Omega$ of closed two forms on refinements of $\ex B$ containing at least one symplectic form. Each form in $\Omega$ must be nonnegative on  $J$ holomorphic planes,  and there must exist a metric $g$ on $\ex B$ and a positive number $c>0$ smaller than the injectivity radius of $g$ so that
  \begin{enumerate}
   \item Given any smooth exploded $\mathbb T$ curve $f:\ex C\longrightarrow \ex B$, the $\omega$ energy, $\int_{\ex C} f^*\omega$ is a constant independent of $\omega\in \Omega$
   \item\label{local area bound} Given any point $p\longrightarrow \ex B$, there exists a taming form $\omega_p\in\Omega$ so that on the ball of radius $c$ around $p$,
   \[\omega_p(v,Jv)\geq g(v,v)\]
     
     \end{enumerate}

  The $\Omega$ energy of a map $f:\ex C\longrightarrow\ex B$ where $\ex C$ is two dimensional and oriented is given by 
  \[E_{\Omega}(f):=\sup_{\omega\in \Omega}\int f^*\omega\]
  
 \end{defn}

For example, any compact almost complex manifold tamed by a symplectic form is strictly tamed by that symplectic form.
A strict taming is what is needed to get local area bounds on holomorphic curves. 

\

The following example is very important for applications of the theory of holomorphic curves in exploded $\mathbb T$ fibrations.

\begin{lemma}\label{fundamental example}
Suppose that a symplectic exploded fibration $(\ex B,\omega)$ satisfies the following:
\begin{enumerate}
\item $\ex B$ is basic and complete.
\item All strata of the base $\totb{\ex B}$ have the integral affine structure of either a standard simplex or $\e{[0,\infty)^n}$
 \end{enumerate}
If the above hold, then there exists some civilized almost complex structure $J$ and a strict taming $\Omega$ of $(\ex B,J)$ so that $\omega\in\Omega$.

\end{lemma}

\pf

Each strata $\ex B_i$ is in the form of some   $\et n{\totb{\ex B_{i}}}$ bundle over the interior of a  smooth compact symplectic manifold $(M_i,\omega_{M_i})$ with a collection of orthogonally intersecting codimension $2$ embedded symplectic submanifolds $M_j$ like in example \ref{relative example} on page \pageref{relative example}. (There is one such $M_j$ for each strata $\totb{\ex B_j}$ which has $\totb{\ex B_i}$ as a boundary face).  By interior we mean the complement of these symplectic submanifolds. 

As in example \ref{relative example}, there exists a neighborhood of $M_j$,  $U_{i,j}\subset M_i$ which is identified with a complex line bundle over $M_i$ with radial coordinate $r_{i,j}$ and connection one form $\alpha_{i,j}$ so that $\omega_{M_i}=\omega_{M_j}+d(r_{i,j}^2\alpha_{i,j})$. We can make these identifications compatible in the sense that if $M_k$ is the intersection of $M_{j_{1}}$ and $M_{j_{2}}$, then on $U_{i,j_{1}}\cap U_{i,j_{2}}$, \[\omega_{M_i}=\omega_{M_k}+ d(r^2_{i,j_{1}}\alpha_{i,j_{1}})+ d(r^2_{i,j_{2}}\alpha_{i,j_{2}})\]
 the $\mathbb T$ actions and projections from the line bundles commute, and  $r_{i,j_{1}}$ and $\alpha_{i,j_{1}}$ are the pullback of $r_{j_{2},k}$ and $\alpha_{j_{2},k}$ under the projection to $M_{j_{2}}$. (To construct these identifications compatibly, start with the strata $\totb{\ex B_{i}}$ of highest dimension.)
 
 \
 
 \psfrag{ETF19a}{$M_{i}:=\totl{\ex B_{i}}$}
 \psfrag{ETF19b}{$U_{i,j_{1}}$}
 \psfrag{ETF19c}{$U_{i,j_{2}}$}
 \psfrag{ETF19d}{$M_{j_{1}}$}
 \psfrag{ETF19e}{$M_{j_{2}}$}
 \psfrag{ETF19f}{$M_{k}$}
 \psfrag{ETF19g}{$r_{i,j_{1}}$}
 \psfrag{ETF19h}{$r_{i,j_{2}}$}
 \psfrag{ETF19i}{$r_{j_{1},k}$}
 \psfrag{ETF19j}{$r_{j_{2},k}$}
\includegraphics{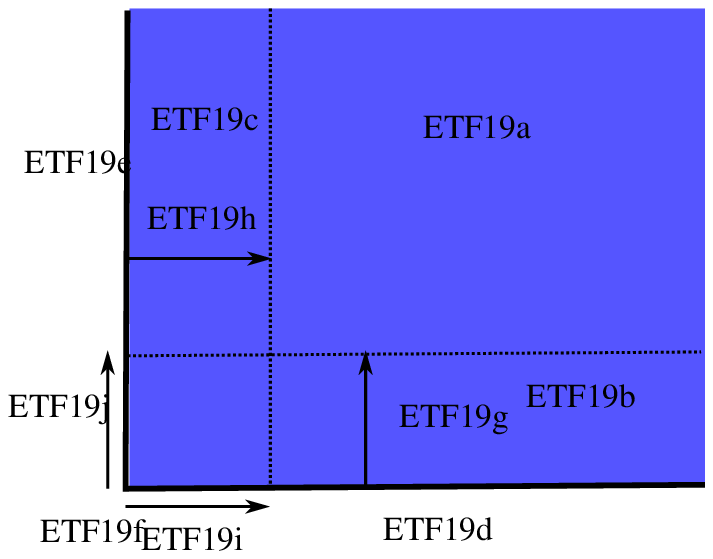}
 
 \
 
  We can arrange this so that in local coordinates in $U_{i,j}$ the fiber coordinate $z$ is actually the smooth part of some coordinate $\tilde z$ from a coordinate chart on $\ex B$. As $\ex B$ is complete, there exists some $r_0>0$ so that the tubular neighborhood around $M_j$ of radius $r_0$ is always contained inside $U_{i,j}$.

We can choose an almost complex structure $J$ tamed by $\omega$ so that on $U_{i,j}$, $J$ is the usual complex structure restricted to the fibers of our line bundle, and is determined by this and the lift (using our connection $\alpha_{i,j}$) of some almost complex structure $J_j$ on $M_j$. (Note that such a $J$ will be civilized in the sense of definition \ref{civilized} on page \pageref{civilized}.)

Denote by ${\ex U}_j$ a connected open subset of $\ex B$ which is the preimage of the union of all $U_{i,j}$ in $\totl{\ex B}$. (This is an open subset that contains the strata $\ex B_j$.) There exist $\mathbb R^*\e{\mathbb R^+}$ valued functions $\tilde r_{i,j}$ defined on ${\ex U}_j$ who's smooth part is equal to  $r_{i,j}$, so that locally there exist standard coordinates so that $\tilde r_{i,j}=\tilde c_k\abs{\tilde z_k}$, and so that $\prod_i\totb{\tilde r_{i,j}}=\e c$ if $\totb{B_j}$ is simplex of finite size. Note that if $\totb{\ex B_j}$ has nonzero dimension, then any adjacent strata $\ex B_k$ with $\totb{\ex B_k}$ of nonzero dimension has the same size, and must obey a similar equation $\prod \totb{\tilde r_{i,k}}=\e c$ with exactly the same constant $\e c$. 

\

\psfrag{ETF20a}{$\totb{\tilde r_{1,j}}$}
\psfrag{ETF20b}{$\totb{\tilde r_{2,j}}$}
\psfrag{ETF20c}{$\totb{\tilde r_{3,j}}$}
\psfrag{ETF20d}{$\totb{\tilde r_{1,j}}\totb{\tilde r_{2,j}}\totb{\tilde r_{3,j}}=\e c$}
\psfrag{ETF20e}{$\totb{\ex B_{j}}$}
\includegraphics{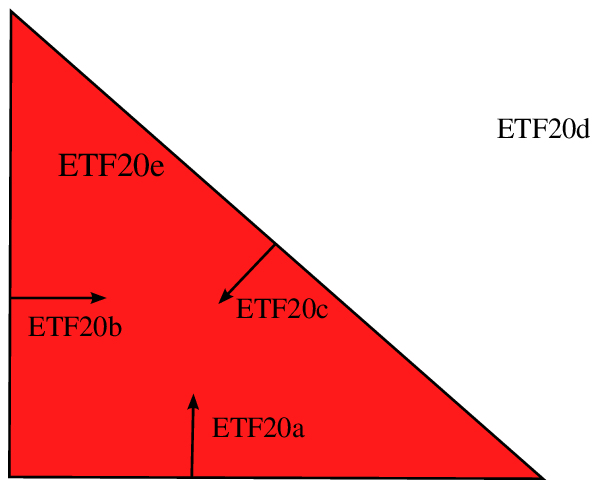}

\

We now provide a strict taming $\Omega$ which for each point $p\longrightarrow \ex B$ includes some symplectic form $\omega_p$ taming $J$ well in a neighborhood of $p$ as in part \ref{local area bound} of the above definition of a strict taming. We do this as follows: Suppose that $p\in \ex U_j$ so that $\tilde r_{i,j}(p)<\frac {r_0}2$ and $\tilde r_{j,i}\geq\frac {r_0}2$ for all $i$ so that this makes sense.
  Consider the set $S_p\subset\mathbb R^*\e{\mathbb R}$ given by the  coordinates of $p$,  $\tilde r_{i,j}(p)$.  We shall define $\omega_p$ on a refinement of $\ex B$ determined by a subdivision of the base $\totb{\ex B}$ as in example \ref{refinement example} on page \pageref{refinement example}. Subdivide the  strata $\totb{\ex B_{j}}$ by the set of all planes of the form 
  $\{\totb{\tilde r_{i,j}}\in\totb{ S_p}\}$. (This is the set of planes through $p$ parallel  to the boundary of the strata and their image under any symmetry of the strata.) 
Similarly subdivide all strata $\ex B_k$ connected to $\ex B_j$ through strata with nonzero dimensional bases by the planes $\{\totb{\tilde r_{i,j}}\in\totb{ S_p}\}$. This subdivision of $\totb{\ex B}$ defines a refinement of $\ex B$ as in example \ref{refinement example}.

\

\psfrag{ETF21a}{$\totb p$}
\includegraphics{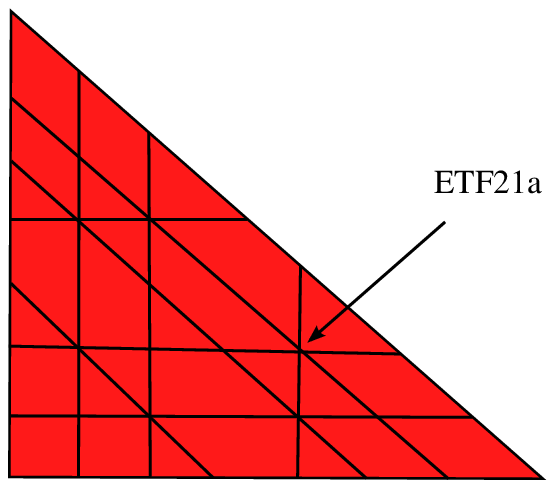}

\

We can now define a symplectic form $\omega_p$ on the refinement of $\ex B$  defined above as follows:  First choose a function 
\[\phi:\mathbb R^*\e{\mathbb R}\longrightarrow\mathbb R \] 
so that 
\[\phi(\tilde r)=0 \ \forall \tilde r>2\e 0\]
\[\phi(\tilde r)+\phi(\tilde r^{-1})=-1\]
and considered as a function $\phi :\mathbb R^*\longrightarrow\mathbb R$, $\phi$ is smooth, increasing and 
\[\frac {d\phi(r)}{dr}>\frac 1 r \text{ for all } \frac 1{\sqrt 2}<r<\sqrt 2\] 
Also choose some smooth monotone function $\rho:\mathbb R\longrightarrow\mathbb R$ so that $\rho(x)=x$ for all $x\geq r^2_0$ and $\rho(x)=\frac {r_0^2}2$ for all $x<\frac {r_0^2}{4}$.

Now define $h:\mathbb R^*\e{\mathbb R^+}\longrightarrow \mathbb R$ by 
\[h(\tilde r):=\rho(\fun{ \tilde r^2})+\frac{r_0^2}{2\abs{S_p}}\sum_{\tilde x\in S_p}\phi( \tilde r\tilde x^{-1}) \]
We shall be replacing $r^2$ with $h(\tilde r)$ in what follows: the relevant properties of $h$ is that it agrees with $r^2$ for $\tilde r$ large enough, is $0$ for $\tilde r$ small enough, and has sufficiently large derivative close to points where $\tilde r\in S_p$, and is monotone increasing.
On $\ex U_j$ we have that $\omega=\omega_{M_j}+\sum d(r^2_{i,j}\alpha_{i,j})$, we shall replace this with 
\[\omega_p:=\omega_{M_j} +\sum d(h(\tilde r_{i,j})\alpha_{i,j}))\]
Now we do the same for all strata $\ex B_k$ connected to $\ex B_j$ through strata with bases of non zero dimension setting on $\ex U_k$
\[\omega_p:=\omega_{M_k} +\sum d(h(\tilde r_{i,k})\alpha_{i,k}))\]
For all other $\ex U_i$, we simply leave $\omega$ unchanged. This defines $\omega_p$ as a two form on our refinement. It is now not difficult to check that the set of forms $\Omega$ given by $\omega$ and $\omega_p$ for all $p$ is a strict taming of $(\ex B,J)$.

\stop

For example, the exploded fibration constructed in example \ref{relative example} from a compact symplectic manifold with orthogonally intersecting codimension $2$ symplectic submanifolds admits an almost complex structure with a strict taming. This is relevant for considering Gromov Witten invariants relative to these submanifolds.

The following  picture of a symplectic form on a refinement in a toric  case may be helpful for understanding the construction in the above proof. (We have drawn $\totl{\ex B}$ in moment map coordiantes.)

\psfrag{ETF8a}{\totl{\ex B}}
\psfrag{ETF8b}{\totb{\ex B}}
\psfrag{ETF8c}{\totl{\ex B'}}
\psfrag{ETF8d}{\totb{\ex B'}}

\includegraphics{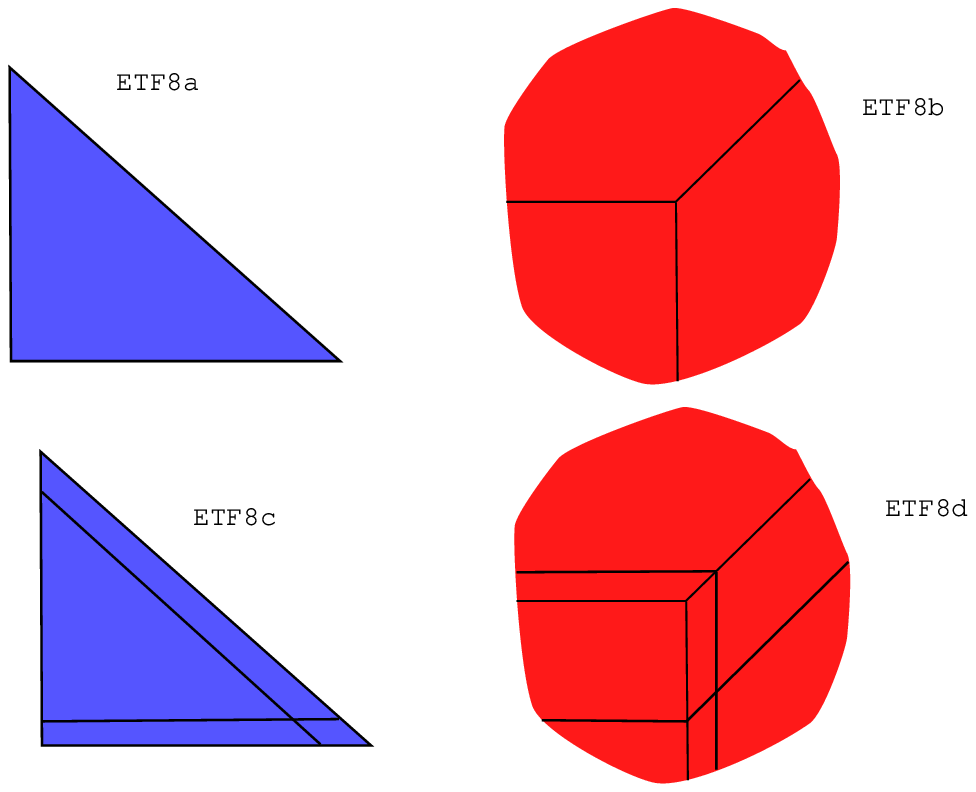}

\begin{example}\label{family example}
\end{example}

Given any symplectic $\ex B$ with a compact tropical part $\totb{\ex B}$ satisfying the conditions of Lemma \ref{fundamental example}, we can construct a family  $\hat{\ex B}\longrightarrow\et 11 $ of symplectic exploded fibrations so that the fiber over $\tilde z=1\e 1$ is $\ex B$, and the fiber over $\tilde z=1$ is a smooth compact symplectic manifold.

\subsection{Moduli stack of smooth exploded $\mathbb T$ curves}\label{perturbation theory}

We shall use the concept of a stack without giving the general definition. (See the article \cite{stacks} for a readable introduction to stacks). The reader unfamiliar with stacks may just think of our use of stacks as a way of keeping all information about families of holomorphic curves in case it is needed for future papers.

 When we say that we shall consider an \exploded fibration $\ex B$ as a stack, we mean that we replace $\ex B$ with a category $\underline{\ex B}$  over the category of \exploded fibrations (in other words a category $\underline{\ex B}$ with a functor to the category of \exploded fibrations) as follows: 
 objects are maps into $\ex B$:
\[\ex A\longrightarrow\ex B\]
and morphisms are commutative diagrams
\[\begin{array}{lll}\ex A & \longrightarrow &\ex B
\\ \downarrow & &\downarrow\id
\\ \ex C &\longrightarrow &\ex B
\end{array}\]
The functor from $\underline{\ex B}$ to the category of \exploded fibrations is given by sending $\ex A\longrightarrow\ex B$ to $\ex A$, and the above morphism to $\ex A\longrightarrow\ex C$.

Note that a maps $\ex B\longrightarrow\ex D$ are equivalent to  functors  $\underline{\ex B}\longrightarrow\underline{\ex D}$ which commute with the functor down to the category of exploded fibrations. Such functors are  morphism of categories over the category of exploded fibrations, and this is the correct notion of maps of stacks. 

For example, a point thought of as a stack is equal to the category of exploded fibrations itself, with the functor down to the category of exploded fibrations the identity. We shall refer to points thought of this way simply as points.

\begin{defn}
 The moduli stack $\mathcal M^{sm}(\ex B)$ of  smooth exploded $\mathbb T$ curves in $\ex B$ is a category over the category of exploded fibrations with objects being families of smooth exploded curves consisting of the following
  
  \begin{enumerate}
 \item 
 A exploded $\mathbb T$ fibration $\ex C$
 \item A pair of smooth exploded $\mathbb T$ morphisms
 
 \[\begin{split}
    &\ex C\longrightarrow \ex B
    \\ \pi &\downarrow
    \\&\ex F
   \end{split}
\]
\item A  section $j$ of $\ker(d\pi)\otimes\left(T^*\ex C/\pi^*(T^*\ex F)\right)$ 
\end{enumerate}
so that
\begin{enumerate}
\item $\pi:\ex C\longrightarrow \ex F$ is a family (definition \ref{family defn} on page \pageref{family defn}).

\item The inverse image of any point $p\longrightarrow \ex F$ is an exploded $\mathbb T$ curve with complex structure $j$.

\end{enumerate}

\

A morphism between families of curves is given by exploded  morphisms $f$ and $c$ making the following diagram commute

\begin{displaymath}
\begin{array}{lllll}
\ex F_1 & \longleftarrow & \ex C_1 & \longrightarrow & \ex B\\ 
\downarrow f &  & \downarrow c & × & \downarrow\id \\ 
\ex F_2 & \longleftarrow & \ex C_2  & \longrightarrow & \ex B
\end{array}
\end{displaymath}
 
so that $c$ is a $j$ preserving isomorphism on fibers. 

The functor down to the category of exploded fibrations is given by taking the base $\ex F$ of a family.
\end{defn}
Note that morphisms are not quite determined by the map  $f:\ex F_{1}\longrightarrow\ex  F_{2}$.  $\ex C_{1}$ is non-canonically isomorphic to the fiber product of $\ex C_{2}$ and $\ex F_{1}$ over $\ex F_{2}$.

\

This is a moduli stack in the sense that a morphism $\underline{\ex F}\longrightarrow \mathcal M^{sm}(\ex B)$ is equivalent to a  family of smooth curves  $\ex F\longleftarrow\ex C\longrightarrow\ex B$ (this is the family which is the image of the identity map $\ex F\longrightarrow \ex F$).

\begin{defn}
A holomorphic curve $\ex C\longrightarrow\ex B$ is stable if it has a finite number of automorphisms, and is not a nontrivial refinement of another holomorphic curve. (If $\ex B$ is basic, this is equivalent to all smooth components of $\ex C$ which are mapped to a point in $\totl{\ex B}$ being  stable as punctured Riemann surfaces.) 
\end{defn}

\begin{defn}
A family of stable holomorphic curves in $\ex B$ is a family of smooth curves so that the map restricted to fibers is holomorphic and stable. The moduli stack of stable holomorphic curves in $\ex B$, is the substack  $\mathcal M(\ex B)\subset\mathcal M^{sm}(\ex B)$ with objects consisting of all families of stable holomorphic curves, and morphisms the same as in $\mathcal M^{sm}$.
\end{defn}

\begin{defn}
 Given a closed $2$ form $\omega$ on the exploded $\mathbb T$ fibration $\ex B$, the smooth exploded $\mathbb T$ curve
 $f:\ex C\longrightarrow\ex B$ 
 is $\omega$-stable if the integral of $f^*\omega$ is nonnegative on any smooth component, and positive on any smooth component which is unstable as a punctured Riemann surface. 

 \
 
Use the notation $\mathcal M^{sm,\omega}(\ex B)$ to indicate the substack of $\mathcal M^{sm}$ consisting of families of $\omega$-stable smooth curves.
\end{defn} 

We have $\mathcal M(\ex B)\subset\mathcal M^{sm,\omega}(\ex B)\subset\mathcal M^{sm}(\ex B)$. The moduli stack of $\omega$-stable curves is a little better behaved than the moduli stack of stable curves because $\omega$-stable curves  have only a finite number of automorphisms. 
 
 We want to make compactness statements about holomorphic curves in families, so we generalize the above definitions for a family $\hat{\ex B}\longrightarrow\ex G$ as follows: 
 
 \begin{defn}
 The moduli stack of smooth curves in a family $\hat{\ex B}\longrightarrow\ex G$, $\ex M^{sm}(\hat{\ex B}\rightarrow\ex G)$ is the substack of $\ex M^{sm}(\hat{\ex B})$ which is the full subcategory which has as objects families which admit commutative diagrams
 \[\begin{array}{cll}(\ex C,j)&\longrightarrow&\hat{\ex B}
 \\ \downarrow& &\downarrow
 \\ \ex F&\longrightarrow & \ex G
 \end{array}\] 
 \end{defn}

The moduli stack of $\omega$-stable smooth curves $\mathcal M^{sm,\omega}(\hat{\ex B}\rightarrow\ex G)$ and the moduli stack of stable holomorphic curves $\mathcal M(\hat{\ex B}\rightarrow\ex G)$ are defined as the appropriate substacks  of $\mathcal M^{sm}(\hat{\ex B}\rightarrow\ex G)$. Note that there is a morphism $\mathcal M^{sm}(\ex B\rightarrow\ex G)\longrightarrow\underline{\ex G}$ which sends the object given by the diagram above to $\ex F\longrightarrow\ex G$. The appropriate compactness theorem for families states that if we restrict to the part of the moduli space with bounded energy and genus, the map $\mathcal M(\hat {\ex B}\rightarrow\ex G)\longrightarrow\underline{\ex G}$ is topologically proper.

 \

The goal of this paper is to prove that the moduli stack of finite energy  stable holomorphic curves with a fixed number of punctures and genus in a complete \exploded fibration with a strict taming is topologically compact. We should say what this means in this context. First, we need a notion of topological convergence. The following definition can be thought of as a notion of what it means for a sequence of points in $\mathcal M^{sm}$ to converge topologically.

\begin{defn}
A sequence of points $\underline p^{i}\longrightarrow \mathcal M(\ex B)$ corresponding to smooth curves  $f^{i}:\ex C^{i}\longrightarrow \ex B$ converges topologically to $ f:\ex C\longrightarrow\ex B$ corresponding to $\underline p\longrightarrow\mathcal M(\ex B)$ in $C^{\infty,\delta}$ if there exist a sequence of families 
  
  \[\ex F\longleftarrow (\ex {\hat C},j_i)\xrightarrow{\hat f^i}\ex B\]
 
  so that this sequence of families converges in $C^{\infty,\delta}$ to 
  
  \[\ex F\longleftarrow (\ex {\hat C},j)\xrightarrow{\hat f}\ex B\]
   
  and a sequence of maps $p^{i}$ of a point into $\ex F$ so that $p^i$ converges to $p$ in $\totl{\ex F}$, $f^{i}$ is the map given by the restriction of $\hat f^{i}$ to the fiber over $p^{i}$, and $f$ is given by the restriction of $f$ to the fiber over $p$. 
  
\end{defn}

(It may seem strange that we use $C^{\infty,\delta}$ instead of smooth convergence for the moduli stack of smooth curves. This is an artifact of our methods of proof. It is highly likely that in the case of integrable complex structures the analogous theorem can be proven with smooth convergence replacing $C^{\infty,\delta}$ convergence. It is also perhaps more natural for our case to consider the moduli stack of $C^{\infty,\delta}$ families. The required results in this setting follow trivially from results in the smooth setting. It is expected that $\mathcal M$ carries a natural type of $C^{\infty,\delta}$ Kuranishi structure.)

\

We say that a sequence of points in $\mathcal M^{sm}(\hat{\ex B}\rightarrow\ex G)$ converge topologically if they converge topologically in $\mathcal M^{sm}(\hat{\ex B})$. Topological convergence in $\mathcal M$ and $\mathcal M^{sm,\omega}$ is defined to simply be topological convergence in $\mathcal M^{sm}$.  We say that a sequence of points in $\underline {\ex G}$ converges topologically if and only if the corresponding sequence of points in $\ex G$ converges topologically. 

\begin{defn}
To say that $\mathcal M$ is topologically compact in $C^{\infty,\delta}$ means that given any sequence of points in $\mathcal M$, there exists a subsequence that converges topologically in $C^{\infty,\delta}$ to some point in $\mathcal M$.

To say that map of stacks $f:\mathcal M\longrightarrow\underline{\ex G}$ is topologically proper in $C^{\infty,\delta}$ is to say that if $\underline p^{i}\longrightarrow\mathcal M$ is some sequence of points so that $f(\underline p^{i})$ converges topologically in $\underline{\ex G}$, then there exists a subsequence of $\{\underline p^{i}\}$ that converges topologically in $C^{\infty,\delta}$ to some point in $\mathcal M$. 
\end{defn}

\

In some sense, the moduli stack of stable holomorphic curves should be complete. (The sense in which this is true involves dealing with transversality issues.) We shall not say much about the other requirement for completeness in this paper, but note that if $\ex U\subset \ex T$ is a connected open subset of $\ex T$, and $\underline{\ex U}\longrightarrow\mathcal M^{sm,\omega}$ is a family of smooth curves so that at least one curve is holomorphic, then the family is holomorphic $\underline{\ex U}\longrightarrow\mathcal M$.

\

The main theorem of the paper, on page \pageref{completeness theorem} can now be understood.

 \section{Estimates}
 
 In this section, we shall prove the analytic estimates required for our compactness theorem. These estimates do not differ much from the standard estimates used to study (pseudo)holomorphic curves in smooth symplectic manifolds. They are also proved using roughly the same methods as in the smooth case. We shall make the following assumptions:
 
  \begin{enumerate}
  \item We'll assume that our target $\ex B$ is complete. This is required to give us the kind of bounded geometry found in the smooth case when the target is compact or has a cylindrical end. We shall also assume that $\ex B$ is basic. The assumption that $\ex B$ is basic  is mainly for convenience in the arguments that follow. It is most useful in establishing  coordinates and notation for the construction of families of smooth curves in the proof of Theorem \ref{completeness theorem}.
 \item We'll assume our almost complex structure $J$ on $\ex B$ is civilized. This artificial (but easy to achieve) assumption means that $J$ induces a smooth almost complex structure on the smooth parts of $\ex B$. This allows us to use standard estimates from the study of holomorphic curves in smooth manifolds with minimal modification. This assumption is needed to have holomorphic curves be smooth morphisms instead of just $C^{\infty,\delta}$ morphisms.
\item We'll assume that $J$ has a strict taming $\Omega$. This is required to get a local area bound for holomorphic curves.
 We shall also work with a particular symplectic form $\omega\in\Omega$. This together with $J$ defines a pseudo metric on $\ex B$ which collapses all tropical directions (so it is like a metric on the smooth part $\totl {\ex B}$.) The strategy of proof for most of our estimates is roughly as follows:  first prove a weak estimate in $\omega$'s pseudo metric, and then improve this to get an estimate in an actual metric on $\ex B$. We shall use this strategy to prove estimates for the the derivative at the center of  holomorphic disks with bounded $\Omega$ energy and small $\omega$ energy, and to get strong estimates on the behavior of holomorphic cylinders with bounded $\Omega$ energy and small $\omega$ energy. 
 \end{enumerate}
 
 This section ends with Proposition \ref{decomposition proposition}, which is a careful statement of the fact that any holomorphic curve with bounded $\Omega$ energy, genus, and number of punctures can be decomposed into a bounded number of components, which either have  `bounded conformal structure  and bounded derivative', or are annuli with small $\omega$ energy, (and so we have strong estimates on their behavior.) 
 
 Throughout this section, we shall use the notation $(\ex B, J, \Omega,\omega,g)$ to indicate a smooth, basic, complete exploded fibration $\ex B$, a civilized almost complex structure $J$, a strict taming $\Omega$, a choice of symplectic taming form $\omega\in\Omega$, and a metric $g$.
  If we say that a constant depends continuously on $(\ex B,J,\Omega,\omega,g)$, we mean that if we have a family $\hat{\ex B}\longrightarrow\ex G$ with such a structure on each fiber, the constant can be chosen a continuous function on $\ex G$. In what follows, we will generally only prove lemmas for a single exploded fibration $\ex B$, and leave the proof of the corresponding statement for families to the diligent reader.
 
 \

We shall use the following lemma which bounds the geometry of a complete \exploded fibration. Note that given a metric $g$ on some complete $\ex B$, the injectivity radius  is a continuous function on $\ex B$, and is therefore bounded below.

\begin{lemma}\label{bounded geometry}
 Given a complete \exploded fibration $\ex B$ with a metric $g$ and any finite number of smooth tensor fields $\theta_i$ (such as an  almost complex structure or a symplectic structure), for any $R>0$ which is smaller than the injectivity radius of $\ex B$, if a sequence of points $p_n\longrightarrow\ex B$ converges topologically to $p\longrightarrow \ex B$, then $(B_R(p_n),g,\theta_i)$ converges  to $(B_R(p),g,\theta_i)$ in the following sense:
 
 Considering $B_R(p_i)$ as a smooth manifold, there exist a sequence of diffeomorphisms $f_n:B_R(p_n)\longrightarrow B_R(p)$ so that $f_n^*( g)$ and $f_n^*(\theta_i)$ converges in $C^\infty$ on compact subsets of $B_R(p_n)$ to $g$ and $\theta_i$.

\end{lemma}

\pf

 This follows from the fact that all smooth tensor fields (and their derivatives) can be locally written as some sum of smooth functions on $\ex B$ times some basis tensor fields.
  
\stop

This tells us that the geometry of a complete \exploded fibration is bounded in the same way that the geometry of a compact manifold is bounded. The analogous result for families is also true.
 
\

 The following is a consequence of our bounded geometry and the standard monotonicity lemma for holomorphic curves first proved in \cite{gromov}. 
 (The proof involves a bound on the covariant derivatives of $J$, the curvature of $g$ and the injectivity radius, all of which vary continuously in families.)
 
 \begin{lemma}\label{strict monotonicity}
  Given  $(\ex B,J,\Omega,g)$, for any $\epsilon>0$ there exists some $E>0$ depending continuously on $(\ex B,J,\Omega,g)$ so that any non constant $J$ holomorphic curve which passes through a point $p$ and which is a complete map when restricted to the $\epsilon$ ball around $p$ has $\Omega$ energy greater than $E$.
 \end{lemma}

\begin{lemma}\label{omega monotonicity}
Given $(\ex B,J,\Omega,\omega)$, define the pseudo-metric 
\[\langle v,w\rangle_\omega:=\frac 12(\omega(v,Jw)+\omega(w,Jv))\]
For any such $(\ex B,J,\Omega,\omega)$ and a choice of $\epsilon>0$, there exists some $E$ depending continuously on $(\ex B,J,\Omega,\omega)$ so that any non constant $J$ holomorphic curve which passes through a point $p$ and which is a complete map when restricted to the $\epsilon$ ball around $p$ in the above pseudo metric has $\omega$ energy greater than $E$.
\end{lemma}
 
 \pf

 The fact that $\ex B$ is basic and complete means in particular that it can be covered by a finite number of coordinate charts, $U\subset\mathbb R^{2n}\times \et mA$. $\mathbb R^{2n}\times \et mA$ with $J$ and $\omega$ can be regarded as subsets of more standard charts $\mathbb R^{2n}\times \et kl$ cut out by setting some monomials equal to $1$. This has the advantage that the smooth part of this is just equal to $\mathbb R^{2n}\times \mathbb C^l$, and we can assume our symplectic and almost complex structures just come from smooth ones here (this uses that $J$ is civilized), so we can use standard holomorphic curve results. Note also that the above pseudo metric is actually a metric on this smooth part. The size of covariant derivatives of $J$ in this metric give continuous functions on $\ex B$, and are therefore bounded because $\ex B$ is complete. To have the bounded geometry required to apply the standard monotonicity lemma for holomorphic curves, we also need some kind of injectivity radius estimate. For a point $q\in\totl{U}$ in the smooth part of a coordinate chart, we can define $i_U(q)$ to be the injectivity radius at $q$ of $\totl U$ with our metric.  We can then define an injectivity radius for a point $p\in\tot{ \ex B}$ to be given by
 \[i(p):=\inf_q\left(d(p,q)+\max_Ui_U(q)\right)\] 
 where $d(p,q)$ indicates distance in our pseudo metric. This defines a continuous function on $\ex B$, and at each point $i(p)>0$, so the fact that $\ex B$ is complete tells us that there is some lower bound $0<c\leq i(p)$.
 
 To prove our lemma we can now take $0<\epsilon<c$. If our holomorphic curve is complete restricted to the $\epsilon$ ball around $p$, it is also complete (which implies that the maps of smooth components are proper maps) when restricted to the $\epsilon$ ball around any $q\in\totl U$ so that $d(p,q)=0$. We can apply the standard monotonicity lemma for holomorphic curves to these balls, so we can choose some $E>0$ so that any holomorphic curve passing through $q$ which is complete restricted to the ball around $q$ is either constant in $\totl U$ or has energy greater than $E$. So either our holomorphic curve must have energy greater than $E$ or it must have some connected component which is constant in the smooth part of every coordinate chart. This connected component must therefore not touch the boundary of our $\epsilon$ ball, and must be a complete map to $\ex B$. The fact that $\Omega$ is a strict taming then allows us to use Lemma \ref{strict monotonicity} to complete our proof.
 
 \stop

 The above proof also implies the following:
 
 \begin{lemma}\label{smooth component monotonicity}
 Given  $(\ex B,J,\Omega,\omega)$, there exists some $\epsilon>0$ which can be chosen to depend continuously on $(\ex B,J,\Omega,\omega)$ so that given any complete holomorphic curve in $\ex B$, the $\omega$ energy of any smooth component is either greater than $\epsilon$, or the image of that smooth component has zero size in the $\omega$ pseudo-metric.

 \end{lemma}

 \begin{lemma}\label{cylinder bound}
  Given $(\ex B,J,\Omega,\omega)$, for any  $\epsilon>0$, there exists some $E>0$ depending continuously on $(\ex B,J,\Omega,\omega)$ so  that any $J$ holomorphic map  $f$ of $e^{-(R+1)}<\abs z<e^{(R+1)}$ with $\omega$-energy less than $E$ is contained inside a  ball of radius $\epsilon$  (in $\omega$'s pseudo metric) on the smaller annulus $e^{-R}\leq\abs z\leq e^R$.
 
 \end{lemma}

 \pf
 
 Suppose to the contrary that this lemma is false.

 Choose $E$ small enough so that  Lemma \ref{omega monotonicity} holds for  $\frac \epsilon 8$ balls. There must be a path in $e^{-(R+1)}<\abs z< e^{R+1}$ joining $z=1$ with an end of this annulus contained entirely inside the ball $B_{\frac \epsilon 8}(f(1))$. (Otherwise, $f$ restricted to some subset would be a proper map to $B_{\frac \epsilon 8}(f(1))$, and we could use Lemma \ref{omega monotonicity} to get a contradiction.) Suppose without loss of generality that it connects $1$ with the circle $\abs z=e^{-(R+1)}$.

 \
 
Suppose for a second that there exists some point $z$ so that $e^{-(R+\frac 34)}<\abs z<e^{-(R+\frac 14)}$ and $f(z)\notin B_{\frac{3\epsilon} 8}(f(1))$.
As above, there must exist some path connecting $z$ with an end which is contained in the ball $B_{\frac \epsilon 8}(f(1))$. There must therefore be some region conformal to $[0,\pi]\times[0,\frac 14]$ so that the   $\dist\left( f(0,t),f(\pi,t)\right)>\frac \epsilon 8$.  Then the Cauchy-Schwartz inequality tells us that 
\[\int_0^\pi\int_0^{\frac 14}\abs{ df}^2_\omega>\frac 14\frac{\left(\frac \epsilon 8\right)^2}{\pi} \]
This in turn gives a uniform lower bound for the energy of our curve. This case can therefore be discarded, and we can assume that all $z$ satisfying 
$e^{-(R+\frac 34)}<\abs z<e^{-(R+\frac 14)}$ are contained in $B_{\frac{3\epsilon} 8}(f(1))$. 

\

As we are assuming this lemma is false, there must be some point $z_2$ so that $e^{-R}\leq\abs {z_2}\leq e^R$ and $f(z_2)\notin B_{\epsilon}(f(1))$. Repeating the above argument (with $z_2$ in place of the point $z=1$) gives us that $z$ must be contained in $B_{\frac {3\epsilon}8}(f(z_2))$ when $e^{R+\frac 14}\leq \abs z\leq e^{R+\frac 34}$. 

We then have $f$ contained in balls  which are at least $\frac \epsilon 4$ apart on the boundary of an annulus, thus we can apply  Lemma \ref{omega monotonicity} to some point in the image under $f$ of the interior at least $\frac \epsilon 8$ from both balls to obtain a lower bound for the energy and a contradiction to the assumption that the lemma was false.

\stop

\

The following gives useful coordinates for the analysis of holomorphic curves:
\begin{lemma}\label{coordinate bound}
Given  $(\ex B,J,\Omega,\omega)$, there exist constants $c_i$ and $\epsilon>0$ depending continuously on $(\ex B,J,\Omega,\omega)$ so that for all points $p\longrightarrow \ex B$, there exists some coordinate neighborhood $U$ of $p$ so that the smooth part, $\totl{U}$ can be identified with a  relatively closed subset of the open $\epsilon$ ball in $\mathbb C^n$ with almost complex structure $\hat J$ and some flat $\hat J$ preserving connection $\nabla$ so that
\begin{enumerate}
\item $p$ is sent to $0$, and the restriction of $\hat J$ to $\totl U$ is $J$.
\item  The metric $\langle\cdot,\cdot\rangle_\omega$ on $\totl U$ is close to the standard flat metric on $\mathbb C^n$ so that 
\[\frac 12 \langle v,v\rangle_\omega\leq \abs v^2\leq 2 \langle v,v\rangle_\omega\]
\item $\hat J$ at the point $z\in\mathbb C^n$ converges to the standard complex structure $J_0$ on $\mathbb C^n$ as $z\rightarrow 0$ in the sense that 
\[\norm{\hat J-J_0}\leq c_0 \abs z\]
 \item The torsion tensor 
 \[\T_\nabla (v,w):=\nabla_vw-\nabla_wv-[v,w]\]
 is bounded by $c_0$, and has its $k$th derivatives bounded by $c_k$.
 \end{enumerate}
\end{lemma}


The only point in this lemma which does not follow from the definition of a civilized almost complex structure and calculation in local coordinates is the fact that the constants involved do not have to depend on the point $p$. This follows from the fact that $\ex B$ is complete.

\begin{lemma}\label{dbar of derivative}
 Given a holomorphic map $f$ of the unit disk to a space with a flat $J$ preserving connection $\nabla$, $f_x$ defines a map to $\mathbb C^n$ defined by parallel transporting $f_x$ back to the tangent space at $f(0)$ which we then identify with $\mathbb C^n$. Such an $f_x$ satisfies the following equation involving the torsion tensor of $\nabla$, $\T_\nabla$.
 
 \[\dbar f_x=\frac 12 J\T_\nabla(f_y,f_x)\]
\end{lemma}

\pf

\[\begin{split}\dbar f_x&=\frac 12(\nabla_{f_x}f_x+J\nabla_{f_y}f_x)
 \\&=\frac 12 \nabla_{f_x}(f_x+J f_y) +\frac 12J(\nabla_{f_y}f_x-\nabla_{f_x}f_y)
 \\&=\frac 12 J\T_\nabla(f_y,f_x)  
  \end{split}\]
\stop

As the $\dbar$ above is the standard linear $\dbar$ operator, this expression is good for applying the following
standard elliptic regularity lemma for the linear $\dbar$ equation. This allows us to get bounds on higher derivatives from bounds on the first derivative of holomorphic functions.

\begin{lemma}\label{elliptic regularity}
For for a given number $k$ and $1<p<\infty$, there exists a constant $c$ so that given any map $f$ from the unit disk $D(1)$,

\[\norm f_{L_{k+1}^p(D(\frac 12))}\leq c\left(\norm{\dbar f}_{L_k^p(D(1))}+\norm f _{L_k^p(D(1))}\right)\]

\end{lemma}

\begin{lemma}\label{derivative bound}
Given $(\ex B,J,\Omega,\omega)$, there  exists some energy $E$, some distance $r>0$, and a constant $c$, each depending continuously on $(\ex B,J,\Omega,\omega)$ so that any holomorphic map $f$ of a disk $\{\abs z\leq 1\}$ into $\ex B$ with $\omega$-energy less than $E$ or contained inside a ball of $\omega$ radius $r$, satisfies 
\[\abs {df}_\omega <c\text{ at }z=0\] 
\end{lemma}

We shall omit the proof of the above Lemma which is a standard bubbling argument, similar, but easier than the proof of the following:

\begin{lemma}\label{strict derivative bound}
Given $(\ex B,J,\Omega,\omega,g)$, for any $E>0$, there exists some $\epsilon>0$, distance $r>0$ and constant $c$, each depending continuously on $(\ex B,J,\Omega,\omega,g)$ so that any holomorphic map $f$ of a disk $\abs z\leq 1$ into $\ex B$ with $\Omega$-energy less than $E$, and either contained in a $\omega$-ball of radius $r$ or having $\omega$-energy less than $\epsilon$ satisfies
\[\abs {df}_g <c\text{ at }z=0\] 

\end{lemma} 

\pf

Suppose that this lemma was false. Then there would exist some sequence of maps $f_i$ satisfying the above conditions with $\abs {df_i(0)}\rightarrow\infty$. 
First, we obtain a sequence of rescaled $J$ holomorphic maps $\tilde f_i
:D(R_i)\longrightarrow \ex B$ from the standard complex disk of radius $R_i$ 
so that
\[\abs{d \tilde f_i}\leq 2\]
\[\abs{d\tilde f_i(0)}=1\]
\[\lim_{i\rightarrow\infty}R_i=\infty\]

We achieve this in the same way as in any standard bubbling off argument such as the proof of lemma 5.11 in \cite{compact}. We can then use Lemma \ref{coordinate bound}, Lemma \ref{dbar of derivative}, and Lemma \ref{elliptic regularity} to get a bound on the higher derivatives of $\tilde f_i$ on $D(R_i-1)$.

As $\ex B$ is complete, we can choose a subsequence so that $f_i(0)$ converges topologically to some $p\longrightarrow \ex B$. Lemma \ref{bounded geometry} tells us that the geometry of $(g,J,\omega)$ around $f_i(0)$ converges to that around $p$.  Considering all our maps as maps sending $0$ to $p$, we can choose a subsequence that converges on compact subsets to a non constant holomorphic map  $f:\mathbb C\longrightarrow \ex B$. Note that $f$ either has $\omega$-energy less than $\epsilon$ or is contained in a ball of radius  $r$ in the $\omega$ pseudo-metric. Because $f$ must have finite $\omega$ energy, we can use Lemma \ref{cylinder bound} to tell us that $f$ must converge in the $\omega$ pseudo-metric as $\abs z\rightarrow\infty$. Then by choosing $\epsilon$ or $r$ small enough, we can prove using Lemma \ref{omega monotonicity} and the standard removable singularity theorem for holomorphic curves, that $f$ must actually have zero size in the $\omega$ pseudo-metric because $f$ either has $\omega$ energy less than $\epsilon$ or is contained in an $\omega$ ball of radius $r$. 

Now that we have that the image of $f$ has zero size in the $\omega$ pseudo metric, we know that $f$ must be contained entirely within some $(\mathbb C^*)^n$ worth of points over a single topological point in $\tot{\ex B}$. This $(\mathbb C^*)^n$ has the standard complex structure, so $f$ gives us $n$ entire holomorphic maps from $\mathbb C\longrightarrow \mathbb C^*$. As $f$ is non constant, at least one of these maps must be non constant. This map must have infinite degree, as it must have dense image in the universal cover of $\mathbb C^*$ which is just the usual complex plane. It therefore has infinite $\Omega$-energy, contradicting the fact that it must have $\Omega$-energy less than $E$.

To prove the family case, it is important to note that (as stated at the start of this section), when we say our constants depend continuously on  $(\ex B,J,\Omega,\omega,g)$, we mean that given any finite dimensional family, $(\hat{\ex B}\longrightarrow\ex G,J,\Omega,\omega,g)$, the constants can be chosen to depend continuously on $\ex G$. With this understood, the proof in the family case is analogous to the above proof.

\stop

 \begin{lemma}\label{omega cylinder bound}
  Given $(\ex B,J,\Omega,\omega)$, for any $0\leq \delta<1$, there exists an energy bound $E>0$, a distance $r>0$, and a constant $c$, each depending continuously on $(\ex B,J,\Omega,\omega)$ so that any holomorphic map $f$ from the annulus $e^{-(R+1)}<\abs z<e^{(R+1)}$ to $\ex B$ with $\omega$ energy less than $E$, or contained inside a ball of $\omega$ radius $r$ satisfies the following inequality
  \[\dist (f(z),f(1))_\omega\leq ce^{-\delta R}(\abs z^\delta+\abs z^{-\delta})
  \text{ for } e^{-R}\leq\abs z\leq e^R\]  
 \end{lemma}

Actually, with a more careful argument, this is true with exponent $\delta=1$, but we will only prove the easier case.

 \pf

For this proof we shall use coordinates $z=e^{t+i\theta}$, and the cylindrical metric in which $\{\frac \partial{\partial t},\frac\partial{\partial \theta}\}$ are an orthonormal frame.

 We choose our energy bound $E>0$ small enough that we can use the derivative bound from Lemma \ref{derivative bound}, and Lemma \ref{cylinder bound} implies that the smaller annulus is contained a small enough ball that we can use the coordinates of Lemma \ref{coordinate bound} (we also choose $r$ small enough that this is true).  We then have the following estimate for the (standard) $\dbar$ of $f$ in these coordinates. 
  
\[\abs{\dbar f}\leq c_1\abs {df}\abs f\]
 
We can run this through the inequality from Lemma \ref{elliptic regularity} to obtain the following estimate on $\dbar f$ restricted to the interior of a disk (which holds on the interior of the cylinder where we have the derivative bound from Lemma \ref{derivative bound}.) 
\begin{equation}\label{dbar bound}\norm{\dbar f}_{L^p(D(\frac 12))}\leq c_1\norm f_{L^p_1D(\frac 12)}\norm f_{L^\infty(D(\frac 12))} \leq c_2\norm f_{L^\infty (D(1))}\norm f_{L^\infty(D(\frac 12))}\end{equation}  

 Here $c_2$ is a constant depending only on $p$ and $(\ex B, J, \omega)$. We can fix $p$ to be something bigger than $2$. By choosing $E$ or $r$ small, we can force $\abs f$ to be as small as we like on the smaller cylinder using Lemma \ref{cylinder bound}. 
 
 Now we can use Cauchy's integral formula 
 
 \[2\pi i f(z_0)=-\int_{\abs z=1}\frac {f(z)}{z-z_0}dz+\int_{\abs z=e^{2l}}\frac {f(z)}{z-z_0}dz+\int_{1\leq\abs z\leq e^{2l}}\frac{\dbar f(z)}{z-z_0}\wedge dz\]
 or in our coordinates, 
 \[\begin{split} f(t_0,\theta_0)&=-\frac 1{2\pi}\int_0^{2\pi} \frac{f(0,\theta)}{1-e^{t_0}e^{i(\theta-\theta_0)}}d\theta+
 \frac 1{2\pi}\int_0^{2\pi} \frac{f(2l,\theta)}{1-e^{t_0-2l}e^{i(\theta-\theta_0)}}d\theta
 \\&+ \frac 1{2\pi}\int_0^{2\pi}\int_0^{2l}\frac {\dbar f(\theta,t)}{1-e^{t_0-t}e^{i(\theta-\theta_0)}}dtd\theta \end{split}\]
 
 Let us now consider each term of this expression for $f(l,\theta_0)$ in the middle of a cylinder under the assumption that the average of $f(2l,\theta)$ is $0$. 
 
 The first term:
 \[\abs{\frac 1{2\pi}\int_0^{2\pi} \frac{f(0,\theta)}{1-e^{l}e^{i(\theta-\theta_0)}}d\theta}\leq\frac 1{e^l-1}\max \abs{f(0,\theta)}\]
 
 The second term using that the average of $f(2l,\theta)$ is $0$,
 \[\begin{split}\abs{\frac 1{2\pi}\int_0^{2\pi} \frac{f(2l,\theta)}{1-e^{l-2l}e^{i(\theta-\theta_0)}}d\theta}
  &=\abs{\frac 1{2\pi}\int_0^{2\pi} f(2l,\theta)\left(\frac 1{1-e^{-l}e^{i(\theta-\theta_0)}}-1\right)d\theta}
   \\&=\abs{\frac 1{2\pi}\int_0^{2\pi} f(2l,\theta)\left(\frac {e^{-l}e^{i(\theta-\theta_0)}}{1-e^{-l}e^{i(\theta-\theta_0)}}\right)d\theta}
   \\&\leq \frac 1{e^l-1}\max\abs{f(2l,\theta)}
   \end{split}\]

   The third term:
 \[\begin{split}\abs{\frac 1{2\pi}\int_0^{2\pi}\int_0^{2l}\frac {\dbar f(\theta,t)}{1-e^{t_0-t}e^{i(\theta-\theta_0)}}dtd\theta}
    &
\leq\norm{\dbar f}_{L^3}\frac 1{2\pi}\norm{\frac 1{1-e^{l-t}e^{i(\theta-\theta_0)}}}_{L^{\frac 32}}
\\&  \leq 
c_3 (l+1)\left(\max_{t\in[-1,2l+ 1]}\abs {f(t,\theta)}\right)
\left(\max_{t\in[0,2l]}\abs {f(t,\theta)}\right)
\end{split}\]  
  
 The constant $c_3$ depends only on $(\ex B,J,\omega)$. It is zero if $J$ is integrable.
 
 To summarize the above, we have the following expression which holds if the average of $f(\theta,2l)$ is $0$. (It also holds if the average of $f(\theta,0)$ is $0$ due to the symmetry of the cylinder.) 
 
 \begin{equation}\label{c estimate}\begin{split}
  \abs {f(l,\theta)}\leq & \frac 1{e^l-1}\left(\max\abs{f(0,\theta)}+\max\abs {f(2l,\theta)}\right)
  \\&+ c_3(l+1)\left(\max_{t\in[-1,2l+ 1]}\abs {f(t,\theta)}\right)
\left(\max_{t\in[0,2l]}\abs {f(t,\theta)}\right)
\end{split}
 \end{equation}

The amount that the average of $f(t_0,\theta)$  changes with $t_0$ inside $[0,2l]$ is determined by the integral of $\dbar f$, which is dominated as above by a term of the form 
\begin{equation}\label{average estimate}\text{change in average }\leq c_4l \left(\max_{t\in[-1,2l+ 1]}\abs {f(t,\theta)}\right)
\left(\max_{t\in[0,2l]}\abs {f(t,\theta)}\right)\end{equation}

 Now let's put these estimates into slightly more invariant terms: Define the variation of $f(t,\theta)$ on $[a,b]\times\mathbb T^{1}$ as follows:
 \[Vf([a,b]):=\max_{t_{1},t_{2}\in[a,b]}\dist(f(t_{1},\theta_{1}),f(t_{2},\theta_{2}))_{\omega}\]
 Now the above two estimates give the following:
 \[Vf([a,a+l])\leq \lrb{\frac 8{e^{l}-1}+c_{5}(l+1)Vf([a-l-1,a+2l+1])}Vf([a-l,a+2l])\]

Now if we choose $l$ large enough, and then make $Vf$ small enough by making our energy bound $E$ small (or directly making $r$ small), we can get the estimate
\[Vf([a,a+l])\leq \frac {e^{-\delta l}}3Vf([a-l,a+2l])
\text{ for }[a-l-1,a+2l+1]\subset[-R,R]\] 

Applying this $3$ times, we get that on the appropriate part of the cylinder,  
\[Vf([a-l,a+2l])\leq e^{-\delta l}Vf([a-2l,a+3l])\]
\[\text{ so }Vf([a,a+l])\leq\frac {e^{-\delta 2l}}3Vf([a-2l,a+3l])\]

Inductively continuing this argument gives that if $[a-nl,a+(n+1)l]\subset[-R+1,R-1]$,
\[Vf([a,a+l])\leq \frac {e^{-\delta n l}}3 Vf([a-nl,a+(n+1)l])\]
The required estimate follows from this.

%
%
%
%
%
%

 \stop

 Lemma \ref{elliptic regularity}, Lemma \ref{dbar of derivative}, and the coordinates from Lemma \ref{coordinate bound} from page \pageref{coordinate bound} can be used with inequality \ref{dbar bound} on page \pageref{dbar bound} in a standard fashion to get a similar estimate on the derivative of $f$ in the $\omega$ pseudo metric. This gives the following corollary:
 
 \begin{cor}\label{cylinder energy decay}
  Given $(\ex B,J,\Omega,\omega)$, there exists some energy $E>0$ depending continuously on $(\ex B,J,\Omega,\omega)$ so that given any $c>0$, there exists some distance $R$ depending continuously on $(\ex B,J,\Omega,\omega,c)$ so that given any holomorphic map $f$ of a cylinder $e^{-l-R}<\abs z<e^{l+R}$ to $\ex B$ with $\omega$ energy less than $E$, the $\omega$ energy of $f$ restricted to $e^{-l}<\abs z<e^{l}$ is less than $ cE $. 
\end{cor}

 \begin{lemma}\label{strong cylinder convergence}
  Given  $(\ex B,J,\Omega,g)$, and any energy bound $E<\infty$ and exponent $0\leq\delta<1$, there exists a covering of $\ex B$ by a finite number of coordinate charts and a constant $c$ so that the following is true:
  
  Given any holomorphic map $f$ of $\Omega$ energy less than $E$ from a cylinder $e^{-(R+1)}\leq \abs z\leq e^{(R+1)}$ to $\ex B$ contained inside a coordinate chart, there exists a map $F$ given in coordinates as 
  \[F(z):=(c_1z^{\alpha_1}\e{a_1},\dotsc,c_kz^{\alpha_k}\e{a_k}, c_{k+1},\dotsc,c_n)\]
 so that 
 \[\dist\left( f(z)- F(z)\right)\leq ce^{-\delta R}\left(\abs z^\delta+\abs{z}^{-\delta}\right)\text{ for }e^{-R}\leq\abs{z}\leq e^{R}\]
  where for $i
  \in [1,k]$, $c_i\in\mathbb C^*$, $a_i\in\mathbb R$,  $\alpha_i\in \mathbb Z$, and for $i\in [k+1,n]$, $c_{i}\in\mathbb R$. $\dist$ indicates distance in the metric $g$.
    
    \
    
  In the case that we a have a family $(\hat{\ex B}\longrightarrow\ex G,J,\Omega,g)$ with the above structure, we can choose $c$ continuous on $\ex G$ and a set of coordinate charts $\hat{\ex U}$ on $\hat{\ex B}$ so that the above coordinate charts on any fiber are the non empty intersections of these with the fiber.   
 \end{lemma}

 \pf 
 
 First, choose any symplectic form $\omega\in\Omega$. Because $\ex B$ is basic and complete, we can then choose a finite number of coordinate charts covering $\ex B$ which are small enough in the $\omega$ pseudo metric so that Lemma \ref{omega cylinder bound} tells us that for any holomorphic map from a cylinder as above contained inside one of these coordinate charts satisfies 
 \begin{equation}\label{omega distance estimate}\dist_\omega(f(z),f(1))\leq ce^{-\delta' R}(\abs z^{\delta'}+\abs z^{-\delta'})\end{equation}
where $\delta'=\frac{\delta+1}2$.  
 
 In our coordinates,
 
 \[f(1)=(c_1\e{a_1},\dotsc,c_k\e{a_k}, c_{k+1},\dotsc,c_n)\]
 
 We can choose $\alpha_i$ so that the winding numbers of the first $k$ coordinates of $f(e^{i\theta})$ are the same as our model map $F$:
 \[F:=(c_1z^{\alpha_1}\e{a_1},\dotsc,c_kz^{\alpha_k}\e{a_k}, c_{k+1},\dotsc,c_n)\] 
 
 The metric $g$ can be compared to the pseudo metric from $\omega$ on the last coordinates $c_{k+1},\dotsc,c_n$, so we only need to prove convergence in the first $k$ coordinates. Let us do this for $f$ restricted to the first coordinate $f_1$.
 
For this, use coordinates $z=e^{t+i\theta}$ on the domain with the usual cylindrical metric and the similar  cylindrical metric on our target. (Choose our coordinate charts so that this metric on the target is comparable to $g$).
Use the notation 
\[\frac{f_1(z)}{c_1z^{\alpha_1}\e {a_1}}=e^{h(z)}\]
where $h$ is a $\mathbb C$ valued function so that $h(0)=0$. The metric we are using on the target is the standard Euclidean metric on $\mathbb C$. Our goal is to prove the appropriate estimate for $h$.

 First note that 
 Lemma \ref{elliptic regularity} and Lemma \ref{dbar of derivative} can be used with the inequality \ref{omega distance estimate} in a standard fashion to get a similar estimate on the derivative of $f$ in the $\omega$ pseudo metric.
 \begin{equation}\abs{df(z)}_\omega\leq ce^{-\delta' R}(\abs z^{\delta'}+\abs z^{-\delta'})\text{ for }e^{-R}\leq \abs z\leq e^{R}\end{equation}
 (The $c$ in the above inequality is some new constant which is independent of f). We can choose our coordinate charts so that if $J_{0}$ indicates the standard complex structure, there exists some constant $M$ so that 
 \[\abs{Jv-J_{0}v}_{g}\leq M{\abs v_{\omega}}\]
 
  This follows from writing down $J$ in coordinates, and the fact that our $\omega$ pseudo metric is a metric on the smooth part of our coordinate chart. This then tells us that 
  $\abs {\dbar f_1}_g$ in our coordinates is controlled by $\abs{df}_\omega$. This then gives that for some new constant $c$ independent of $f$ we get the following inequality:
 
 \begin{equation}\label{dbar h estimate}\abs {\dbar h(z)}\leq c e^{-\delta' R}(\abs z^{\delta'}+\abs z^{-\delta'})\text{ for }e^{-R}\leq \abs z\leq e^{R}\end{equation}  
 
If we choose our coordinate charts small enough in the $\omega$ pseudo-metric, Lemma \ref{strict derivative bound} gives us a bound for $\abs{df}_g$ on the interior of the cylinder, so we have that there exists some $c$ independent of $f$ so that 
 
 \begin{equation}\label{h derivative bound}\abs {dh(z)}< c \text{ for }e^{-R}\leq\abs z\leq e^R\end{equation}
 
 We can now proceed roughly as we did in the proof of Lemma \ref{omega cylinder bound}. In particular, Cauchy's integral theorem tells us that 
  
 \[\begin{split} h(t_0,\theta_0)&=-\frac 1{2\pi}\int_0^{2\pi} \frac{h(t_0-l,\theta)}{1-e^{l}e^{i(\theta-\theta_0)}}d\theta+
 \frac 1{2\pi}\int_0^{2\pi} \frac{h(t_0+l,\theta)}{1-e^{-l}e^{i(\theta-\theta_0)}}d\theta
 \\&+ \frac 1{2\pi}\int_0^{2\pi}\int_{t_0-l}^{t_0+l}\frac {\dbar h(t,\theta)}{1-e^{t_0-t}e^{i(\theta-\theta_0)}}dtd\theta \end{split}\]
 
 We can now bound each term in the above expression as in the proof of Lemma \ref{omega cylinder bound}, except we use the estimate \ref{dbar h estimate} to bound $\abs {\dbar h}$.

 Then we have the following estimate:

 \[\begin{split}
  \abs {h(t_0,\theta_0)-\frac 1{2\pi}\int_0^{2\pi}h(t_0+l,\theta)d\theta }\leq & \frac 1{e^l-1}\max\abs{h(t_0-l,\theta)-\frac 1{2\pi}\int_0^{2\pi}h(t_0+l,\theta)d\theta}
  \\&+\frac 1{e^l-l}\max\abs {h(t+l,\theta)-\frac 1{2\pi}\int_0^{2\pi}h(t_0+l,\theta)d\theta}
  \\&+ c(l+1)e^{-\delta'R}e^{\delta'l}\left( e^{\delta' t_{0}}+e^{-\delta 't_{0}}\right)
\end{split}
 \]

(Of course, this is a new constant $c$, which is independent of $l$ or $h$.)

 Note that the change in the average of $h$ is determined by the integral of  $\dbar h$, which we can bound using estimate \ref{dbar h estimate}. The change in this average can then be absorbed into the last term of the above inequality. Define the variation of $h$ for a particular $t$ as follows:
 
 \[Vh(t):=\max_{\theta}\abs {h(t,\theta)-\frac 1{2\pi}\int_0^{2\pi}h(t,\theta)d\theta }\]
 We then have the following estimate (with a new constant $c$):
 \[Vh(t)\leq \frac 1{e^l-1}(Vh(t+l)+Vh(t-l))+ c(l+1)e^{\delta'(l-R)}\left( e^{\delta' t}+e^{-\delta 't}\right)\]
 
 Recalling that $\delta<\delta'$ and $Vh$ is bounded by equation \ref{h derivative bound}, we can choose $l$ large enough so that for $R$ sufficiently large, using the above estimate recursively tells us that there exists some $c$ independent of $f$ or $R$ so that 
 
 \[Vh(t)\leq ce^{-\delta R}(e^{\delta t}+e^{-\delta t})\text{ for }-R\leq t\leq R\] 
 (This estimate follows automatically from the bound on $Vh$ for $R$ bounded.) The required estimate for $h$ then follows from the fact that the change in the average of $h$ is bounded by the estimate \ref{dbar h estimate}, which is stronger than what we need (which is the same equation with $\delta$ in place of $\delta'$). We therefore have
 \[h(z)\leq ce^{-\delta R}(\abs z^{\delta}+\abs z^{-\delta})\text{ for }e^{-R}\leq\abs z\leq e^{R}\]
 which is the required estimate.
 
 \stop

\

To prove compactness results, we shall divide our domain up into annuli with small energy, and other compact pieces with derivative bounds. For this we shall need some facts about annuli.
Recall the following standard definition for the conformal modulus of a Riemann surface which is an annulus:

\begin{defn}
The conformal modulus of an annulus $A$ is defined as follows. Let $S(A)$ denote the set of all continuous functions with $L^2$ integrable derivatives on $A$ which approach $1$ at one boundary of $A$ and $0$ at the other. Then the conformal modulus of $A$ is defined as 
\[R(A):=\sup_{f\in S(A)}\frac {2\pi}{\int_A(df\circ j)\wedge df}=\sup_{f\in S(A)}\frac {2\pi}{\int_A\abs{df}^2}\] 
\end{defn}

We can extend the definition of conformal modulus to include `long' annuli inside exploded curves as follows:

\begin{defn}
Call a (non complete) exploded curve $\ex A$ an exploded annulus of conformal modulus 
$\log x\e {-l}$ if it is connected, and there exists an injective holomorphic map $f:\ex A\longrightarrow \ex T$ with image $\{1<\abs {\tilde z}< x\e {-l}\}$. Call it an exploded annulus with semi infinite conformal modulus if it is equal to (a refinement of) $\{\abs{\tilde z}<1\}\subset\et 11$.   
\end{defn}

(We use $x\e{-l}$ in the above definition because the tropical part of the resulting annulus will have length $l$.) Two exploded annuli with the same conformal modulus may not be isomorphic, but they will have a common refinement.

\

We shall need the following lemma containing some useful properties of the (usual) conformal modulus.

\begin{lemma}\label{annulus lemma}

\begin{enumerate}
\item\label{conformal modulus 1} An open annulus $A$ is conformally equivalent to $\{1<\abs z< e^R\}$ if and only if  the conformal modulus of $A$ is $R$. If the conformal modulus of $A$ is infinite, then $A$ is conformally equivalent to either a punctured disk, or a twice punctured sphere.
\item\label{conformal modulus 2} If $\{A_i\}$ is a set of disjoint annuli $A_i\subset A$ none of which bound a disk in $A$, then 
\[R(A)\geq\sum R(A_i)\]
\item\label{conformal modulus 3}

 If $A_1$ and $A_2$ are annuli contained inside the same Riemann surface which share a boundary so that $R(A_2)<\infty$, and the other boundary of $A_1$ intersects the other boundary of $A_2$ and the circle at the center of $A_2$, then 
\[R(A_2)+\frac{16\pi^2}{R(A_2)}\geq R(A_1)\]
(This is not sharp.)

\psfrag{ETF9int}{$A_{1}\cap A_{2}$}
\psfrag{ETF9cylinder}{$A_{2}$}
\includegraphics{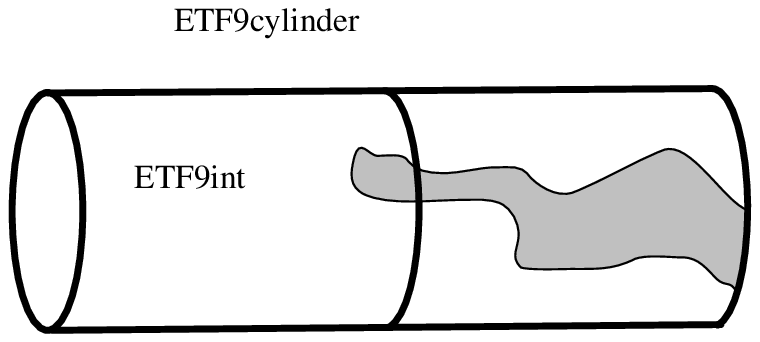}

\item \label{conformal modulus 4} If $A_1$ and $A_2$ are annuli contained inside a Riemann surface so that  every circle homotopic to the boundary inside $A_2$ contains a segment inside $A_1$ that intersects both boundaries of $A_1$, then
\[R(A_1)\leq\frac {4\pi^2}{R(A_2)}\]

\psfrag{ETF10a}{$0$}
\psfrag{ETF10b}{$1$}
\includegraphics{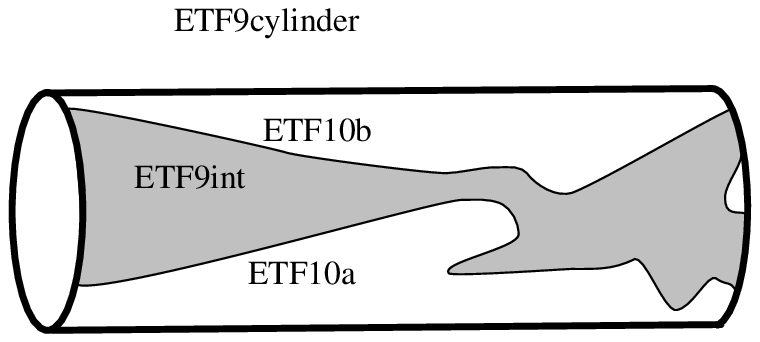}

\end{enumerate}

\end{lemma}

\pf

The first two items are well known. To prove item \ref {conformal modulus 3}, first set $R=R(A_2)$ and put coordinates  $(0,R)\times\mathbb R/2\pi\mathbb Z$ on $A_2$. Consider any function $f\in S(A_1)$. Without losing generality, we can assume that the shared boundary is $(0,\theta)$, and that some segment of the other boundary of $A_1$ is a curve in $A_2$ between $(\frac R2,0)$ and $(R,\theta_0)$ where $\theta_0\geq 0$. Then consider integrating $\abs {df}^2$ along diagonal lines $( Rt, -{4\pi}t+c)$. Each of these lines contains a segment inside $A_1$ on which the integral of $df$ is $1$. The length of each segment is bounded by $((4\pi)^2+R^2)^{\frac 12}$. The integral of $\abs {df}^2$ along this segment is therefore at least $(16\pi^2+R^2)^{-\frac 12}$. This tells us that the integral of $\abs {df}^2$ over $A_1$ is at least $\frac {2\pi R}{((4\pi)^2+R^2)}$, and therefore, 
\[R(A_1)\leq \frac{(16\pi^2+R^2)}{R}=R(A_2)+\frac{16\pi^2}{R(A_2)}\] 

To prove item \ref{conformal modulus 4}, consider a function $f\in S(A_1)$. We shall integrate $\abs{df}^2$ along segments  of the form $(c,t)$ inside $A_2\cap A_1$ traveling from one boundary of $A_1$ to the other. The integral of $df$ along such a segment is $1$, and the length of the segment is at most $2\pi$, so the integral of $\abs{df}^2$ along the segment is at least $\frac 1{2\pi}$, and the integral of $\abs{df}^2$ over $A_1$ is at least $\frac  {R(A_2)}{2\pi}$. Therefore, we have
\[R(A_1)\leq\frac {4\pi^2}{R(A_2)}\] 

\stop

 \begin{prop}\label{decomposition proposition}
 Given $(\ex B,J,\Omega,\omega,g)$, an energy bound $E$ and a number $N$, and small enough $\epsilon>0$, there exists a number bound $M$ depending lower semicontinuously on $(\ex B,J,\Omega,\omega,g,E,N,\epsilon)$ and for any large enough distance $R$, a derivative bound $c$  and conformal bound $\hat R$ depending continuously on $(\ex B,J,\Omega,\omega,g,E,N,\epsilon, R)$ so that the following is true:
 
 Given any complete, stable holomorphic curve $f:\ex C\longrightarrow \ex B$ with energy at most $E$, and with genus and number of punctures at most $N$, there exists some collection of at most $M$ exploded annuli $\ex A_i\subset\ex C$, so that:
 \begin{enumerate}
 \item\label{first inductive condition} Each $\ex A_i$ has conformal modulus larger than $2R$
 
 \item\label{non intersection} Put the standard cylindrical metric on $\ex A_i$, and use the notation $\ex A_{l,i}$ to denote the annulus consisting of all points in $\ex A_i$ with distance to the boundary at least $l$.  \[\ex A_{\frac R 2,i}\cap \ex A_{\frac R 2,j}=\emptyset\text{ if } i\neq j\]

\item  $f$ restricted to $ \ex A_i$ has energy less than $\epsilon$. 

\item\label{conformal bound} Each component of $\ex C-\bigcup \ex A_{R,i}$ is a smooth Riemann surface with bounded conformal geometry in the sense that any annulus inside one of these smooth components with conformal modulus greater than 
$\hat R$ must bound a smooth disk inside that component.

\item\label{proposition derivative bound} The following metrics can be put on $\ex A_i$ and each smooth component of $\ex C-\bigcup \ex A_{R,i}$:
\begin{enumerate}
\item On any component which is a smooth torus, use the unique flat metric in the conformal class of the complex structure so that the area of the torus is $1$. 

\item \label{disk condition}If the component is equal to a disk,  an identification with the standard unit disk can be chosen so that if $\ex A_{R,i}$ is the bounding annulus, $0$ is in the complement of $\ex A_{i}$. Give components such as this the standard Euclidean metric.
\item If the component is equal to some annulus, give it the standard cylindrical metric on $\mathbb R/\mathbb Z\times(0,l)$. Give each $\ex A_i$ the analogous  standard metric.
\item Any component not equal to a torus, annulus, or  disk will admit a unique metric in the correct conformal class with curvature $-1$ so that boundary components are geodesic (there will be no components which are smooth spheres). Give these components this metric.

\end{enumerate}

On any component of $\ex C-\bigcup \ex A_{\frac {9R}{10},i}$, the derivative in this metric is bounded by $c$ 

\[\abs {df}_g<c\]

Moreover, on $\ex A_{\frac {6R}{10},i}-\ex A_{\frac{9R}{10},i}$, the ratio between the two metrics defined above is less than $c$.

\item \label{decomposition stability}


There exists some lower energy bound $\epsilon_0>0$ depending only on $\epsilon$, $R$, and $E$ so that $f$ restricted to any  component of $\ex C-\bigcup \ex A_{{R},i}$ which is a disk, annulus or torus has $\omega$ energy greater than $\epsilon_0$
 \end{enumerate}
 \end{prop}
 
 \psfrag{ETF11d}{{\footnotesize derivative bounded here}}
 \psfrag{ETF11c}{{\footnotesize components with $\omega$ energy bounded below}}
 \psfrag{ETF11b}{$\ex A$}
 \psfrag{ETF11a}{$\ex A_{R}${\footnotesize :well behaved annulus}}
 \includegraphics{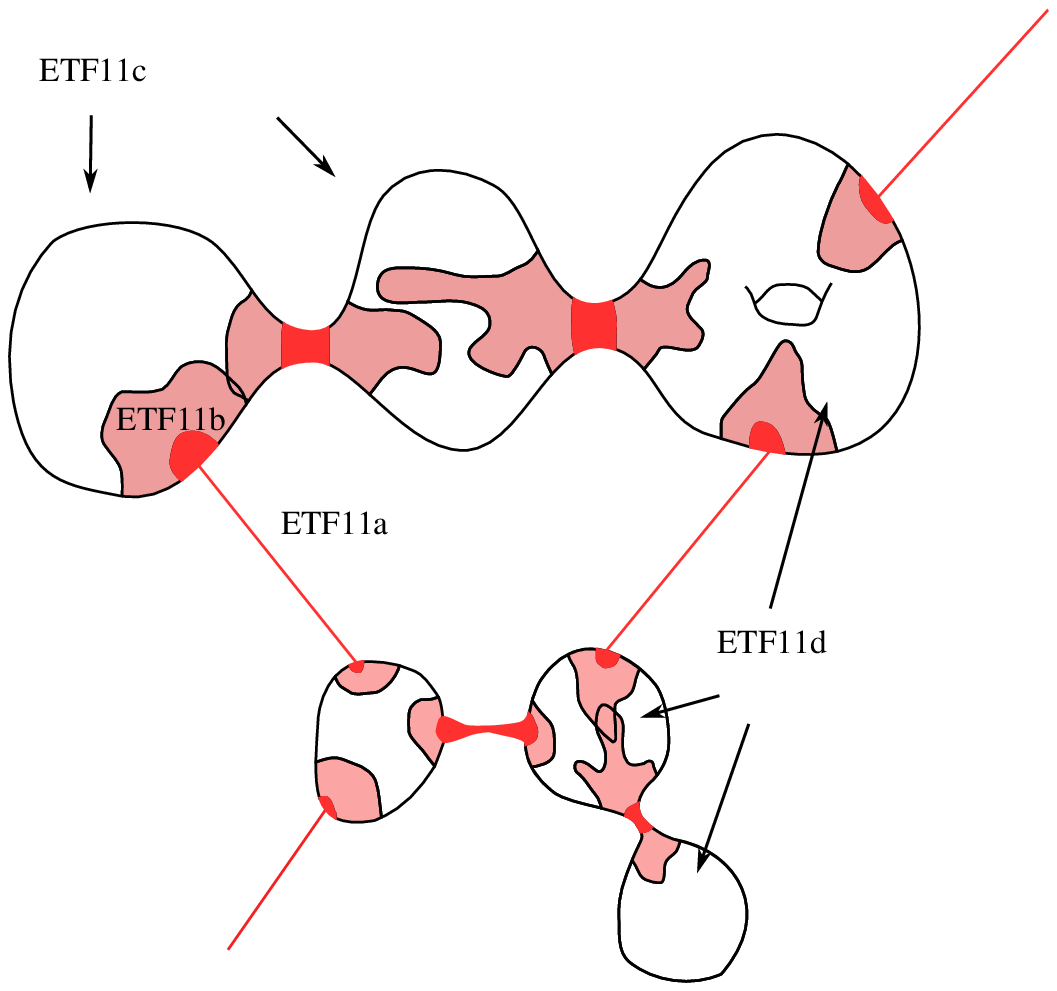}
 
  \pf

Lemma \ref{strict derivative bound} tells us that for $\epsilon$ small enough, any holomorphic map of the unit disk with $\omega$-energy less than $\epsilon$ and $\Omega$-energy less than $E$ must have derivative at $0$ bounded by $c_0$. We shall prove our theorem for $\epsilon$ small enough so that this is true, and small enough that Corollary \ref{cylinder energy decay} also holds with an energy bound of $\epsilon$. We shall also choose our distance $R$ greater than $4\pi$ (for use with Lemma \ref{annulus lemma}), and large enough that Corollary \ref{cylinder energy decay} tells us that if $f$ restricted to some smooth annulus $\ex A_i$ has $\omega$-energy less than $\epsilon$, then the $\omega$-energy of $\ex A_{\frac R 2,i}$ is less than $\frac \epsilon 5$. 

Let us now begin to construct our annuli. First, note that any edge of $\ex C$ has zero $\omega$-energy, so we can choose an exploded annulus $\ex A_i$ containing each edge with $\omega$-energy less than $ \frac \epsilon 5$. (Note for use with item \ref{decomposition stability} that in the case that this bounds a disk, the $\omega$-energy of the resulting disk will be at least $\frac{4\epsilon} 5$ ). We can do this so that these $\ex A_i$ are mutually disjoint. Note that the complement of these $\ex A_i$ is a smooth Riemann surface with boundary. (Note also that no component of $\ex C$ can be isomorphic to $\ex T$ as the  $\Omega$-energy of any complete component is equal to the $\omega$-energy which must be zero on $\ex T$. Similarly, there is no complete component of $\ex C$ which is locally modeled everywhere on $\ex T$.)

 For each connected component of $\ex C$ which is a smooth sphere, we shall now remove an annulus. First, note that we can put some round metric in the correct conformal class of the sphere so that there exist $3$ mutually perpendicular geodesics which divide our sphere into $8$ regions each of which has equal $\omega$-energy. As our sphere must have energy at least $\epsilon$ in order to be stable, we can choose two antipodal regions that each have energy $\frac\epsilon 8$. By choosing $\epsilon_0>0$ small enough, we can then get that there exist two disks with $\omega$-energy at least $\epsilon_0$ which intersect each of these regions, and which have radius as small as we like. Choose $\epsilon_0>0$ small enough that the complement of these disks  is an annulus of conformal radius at least ${k (2R+1)}$ for some integer $k>\frac {5E}\epsilon$. How small $\epsilon_0$  is required to achieve this depends only on $R$, $\epsilon$ and $E$. Then we can divide this annulus into $k$ annuli with conformal modulus at least $(2R+1)$, at least one of which has $\omega$-energy at most  $\frac \epsilon 5$. Add this annulus to our collection. Note that it bounds disks which have energy at least $\epsilon_0$.   
 
 \
 
 \psfrag{ETF12c}{$\omega$-energy $<\frac\epsilon 5$}
 \psfrag{ETF12b}{$\ex A$}
 \psfrag{ETF12a}{$\omega$-energy $>\epsilon_{0}$}
 \includegraphics{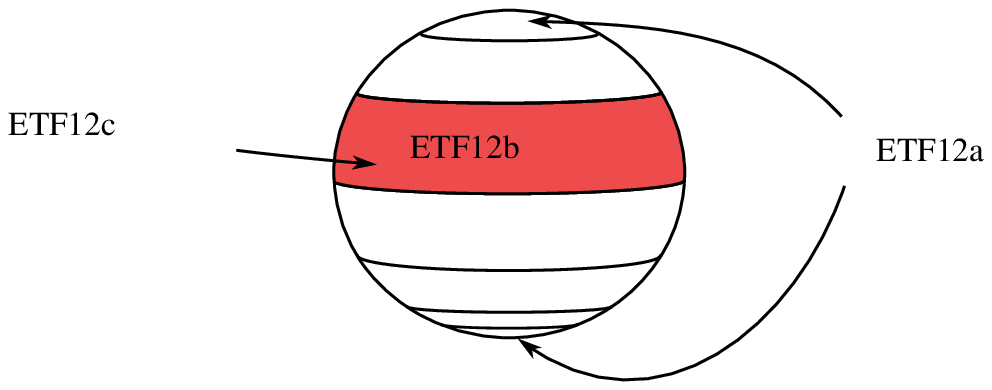}
 
Consider a holomorphic injection of the unit disk $i:D\longrightarrow\ex C$ into the complement of our annuli $\ex A_{R,i}$ constructed up to this point. Suppose that $i(0)$ is in the complement of $\ex A_{\frac R4,i}$. Then Lemma \ref{annulus lemma} tells us that the restriction of $i$ to $\abs z<e^{-\frac {16\pi^2}R}$ is in the complement of $\ex A_{\frac R2,i}$. Now suppose that $\abs {d(f\circ i)}>c_0e^{\frac {16\pi^2}R+k(2R+1)}$ where $k$ is some integer greater than $\frac {5E}\epsilon$. Then the restriction of $f$ to the disk $\abs z< e^{-\frac {16\pi^2}R-k(2R+1)}$ must have energy greater than $\epsilon$ due to Lemma \ref{strict derivative bound}. This disk is surrounded by $k$ disjoint annuli of conformal modulus $(2R+1)$ which are contained in the complement of $\ex A_{\frac R2,i}$. At least one of these must have energy less than $\frac \epsilon 5$. Add this annulus to our collection.  Continue adding annuli in this manner. After a finite (but not universally bounded) number of times, no more annuli can be added in this way. We shall argue this below after adding some more annuli.

\psfrag{ETF13a}{$\ex A_{R}$}

\psfrag{ETF13b}{$\ex A_{\frac R2}$}

\psfrag{ETF13c}{New low energy annulus}

\psfrag{ETF13d}{Disk with $\omega$-energy $>\epsilon_{0}$}

\includegraphics{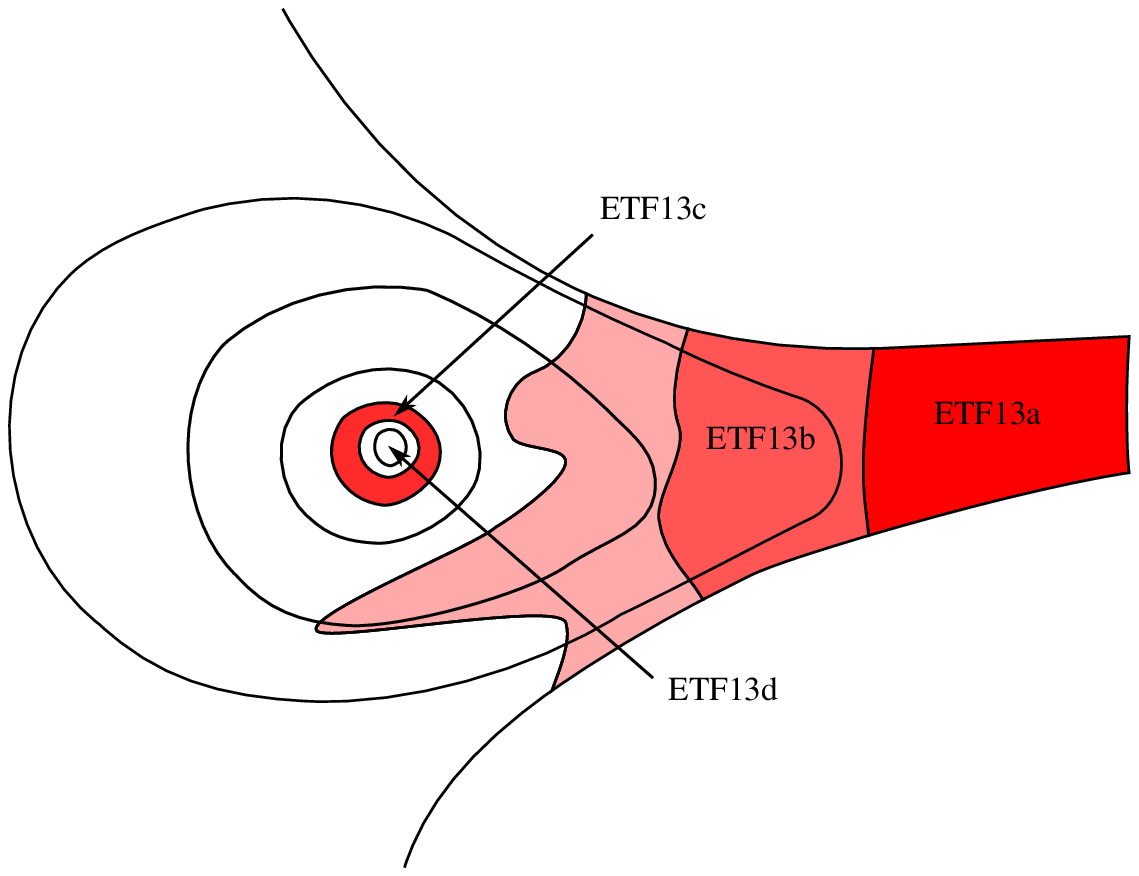}

  We now want to add annuli so that condition \ref{conformal bound} is satisfied. Suppose  some smooth component of $\ex C-\bigcup  \ex A_{R,i}$ contains some annulus of conformal modulus greater than $2R(3+\frac{5E}{\epsilon})$ that does not bound a disk. Note that as we have chosen $R>4\pi$, the size of such an annulus is greater than $2R(1+\frac{5E}{\epsilon})+2R+\frac{32\pi^2}R$. The annulus of size $R$ at its end can't be contained entirely within $\ex A_i$, and then applying Lemma \ref{annulus lemma} part \ref{conformal modulus 4}  with the next annulus  of size $\frac {16\pi^2} R$ in the place of $A_{1}$, we see that our annulus minus annuli at each end of conformal modulus $R+\frac {16\pi^2} R$ must be contained entirely in the complement of $\ex A_{\frac R2,i}$. We then obtain more than $\frac {5E}{\epsilon}$ disjoint annuli in the compliment of $\ex A_{\frac R2,i}$ with conformal modulus greater than $2R$. The restriction of f to at least one of these annuli must have $\omega$-energy less than $\epsilon$.   Add this annulus to our collection.

  \psfrag{ETF14a}{New low energy annulus}
  \psfrag{ETF14b}{$\ex A_{\frac R2}$}
  \psfrag{ETF14c}{$\ex A_{R}$}
  \includegraphics{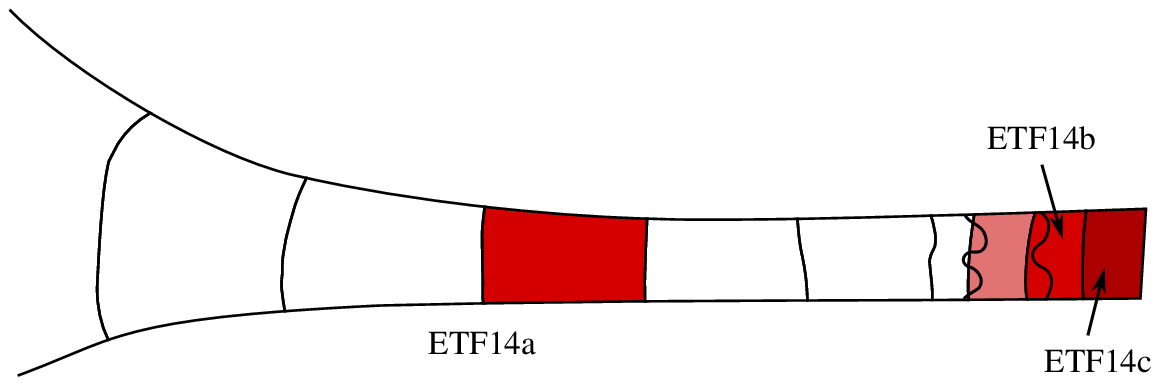}

   We shall now argue that we can only add a finite number of annuli in this manner.
  The number of annuli that bound disks is bounded by $\frac {E}{\epsilon}$ plus the twice the number of connected components of $\ex C$ that are spheres (which is bounded by $\frac E\epsilon$). There are also a finite number of exploded annuli with infinite conformal modulus. The number of these is bounded by our topological bounds, and the fact that there are at most $\frac E\epsilon$ spherical components with fewer than $3$ punctures. Call the complement of all the above annuli $\ex C_0$. There is then a bound on the number of homotopy classes in $\ex C_0$ of our remaining annuli, given by our bounds on the topology of $\ex C$, the fact that we have removed a bounded number of annuli, and the observation that for any Riemann surface with boundary, there are only a finite number of homotopy classes which contain annuli of conformal modulus greater than $4\pi$ (Lemma \ref{annulus lemma} part \ref{conformal modulus 4} can help to prove this). If there were an infinite number of annuli $\ex A_{\frac R2,i}$ in our collection in the same homotopy class, then there would exist an annulus in $\ex C_0$ in that homotopy class which contains all of them. As the $\ex A_{\frac R2,i}$ are disjoint and have conformal modulus at least $R$, this annulus would have to have infinite conformal modulus. Note that there are no nontrivial annuli of infinite conformal modulus in $\ex C_0$ (because we have removed the annuli containing edges and at least one annulus from each spherical component), and if an annulus of infinite conformal modulus surrounds a disk with energy greater than $\epsilon_{0}$, then the derivative of $f$ must have been unbounded there. This is not possible, so we must only have a finite number of annuli in our collection. Note that everything apart from our bound on the number of annuli in a fixed homotopy class is bounded independent of $f$. 

Now we shall merge some annuli so that there exist no annular component of $\ex C-\bigcup \ex A_{R,i}$ with $\omega$-energy less than $\frac\epsilon 5$. Then we will have a bound on the number of annuli which is independent of $f$. Suppose that we have some collection $\{\ex A_1,\dotsc,\ex A_n\}$ of our annuli so that $\ex A_{R,i}$ and $\ex A_{R,i+1}$ bound an annulus which has $\omega$-energy less than $\frac \epsilon 5$. Use the notation $\ex A_{[m, n]}$ to denote the annulus that consists of $\ex A_m$, $\ex A_n$ and everything in between. Then, as we've chosen each of our $\ex A_i$ to have $\omega$-energy less than $\frac\epsilon 5$, $\ex A_{[i,i+2]}$ has $\omega$-energy less than $\epsilon$, so we can apply Corollary \ref{cylinder energy decay} to show that far enough into the interior of this cylinder, there is very little $\omega$-energy. (Note that as the parts of exploded annuli which are locally modeled on $\ex T$ always have no $\omega$-energy, so there is no difficulty in applying this lemma in the seemingly more general setting of exploded annuli.) Applying Lemma \ref{annulus lemma} part \ref{conformal modulus 3}, and noting that as we have chosen $R>4\pi$, we have $2R>R+\frac{16\pi^2}R$, we see that $\ex A_{[i,i+2]}-\ex A_i-\ex A_{i+2}\subset \ex A_{\frac R2, [i,i+2]}$. We have chosen $R$ large enough that Corollary \ref{cylinder energy decay} tells us that the energy of $f$ restricted to $\ex A_{\frac R2,[i,i+2]}$ is less than $\frac \epsilon 5$. We can now repeat this argument inductively to show that the energy of $f$ restricted to $\ex A_{[1,n]}$ is less than $\epsilon$ and $\ex A_{\frac R2, [1,n]}$ is less than $\frac\epsilon 5$. 

\psfrag{ETF15a}{$\ex A_{1}$}
\psfrag{ETF15b}{$\ex A_{2}$}
\psfrag{ETF15c}{$\ex A_{3}$}
\psfrag{ETF15d}{$\ex A_{4}$}
\psfrag{ETF15e}{$\ex A_{\frac R2,[1,3]}$ has low $\omega$-energy}

\includegraphics{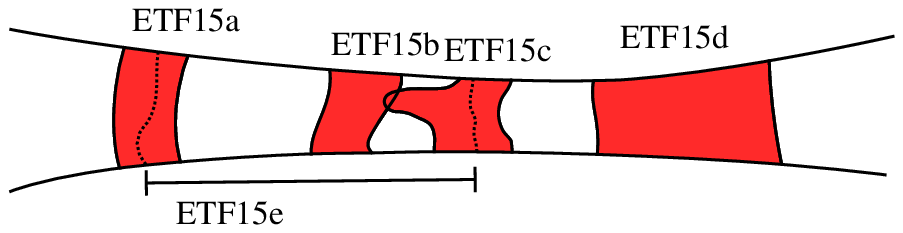}

\

 Now, replace all of our sets of annuli $\{\ex A_1,\dotsc, \ex A_n\}$ of maximal size obeying the above conditions with the annulus given by $\ex A_{\frac R2,1}$, $\ex A_{\frac R2,n}$ and everything in between. (Note that if we had some collection of annuli which bounded small energy annuli in a cyclic fashion, we would obtain a connected component of $\ex C$ which had energy less than $\epsilon$, and was not stable.) 
 
 We shall now check that all the conditions we require are satisfied by this new set of annuli. First, the number is bounded independent of $f$. The size of each annulus is greater than $2R$.
 Our resulting set of annuli also obey the non-intersection condition \ref{non intersection} because the original set did. We showed above that the energy of each new annulus is less than $\epsilon$. 
Because the complement of our original set of annuli obeyed condition \ref{conformal bound} with a bound of $2R(3+\frac{5E}\epsilon)$, by increasing our conformal bound appropriately and using Lemma \ref{annulus lemma}, we can achieve condition \ref{conformal bound}. The lower bound on the $\omega$-energy of unstable components is also satisfied by construction.

All that remains is condition \ref{proposition derivative bound}. Note first that the  complement of our annuli $\ex A_{R,i}$ admits metrics as described in \ref{proposition derivative bound}. The fact that there exists some $c$ so that the metrics we choose on different components differ by a factor of less than $c$ on $\ex A_{\frac {6R}{10},i}-\ex A_{\frac{9R}{10},i}$ follows from conditions \ref{conformal bound} and \ref{non intersection}, and the bound on topology as follows: The complement of $\ex A_{R,i}$ is a Riemann surface with bounded genus, a bounded number of boundary components and bounded conformal geometry from condition \ref{conformal bound}.  An easy compactness argument tells us that there exists a constant $c'$ so that any injective holomorphic map of the unit disk into either $\ex A_{i}$  or the complement of $\ex A_{R,i}$ 
has derivative at $0$ bounded by $c'$. We can  use this (remembering that $R>4\pi>10$ so we can get a disk of unit size centered on any point inside $\ex A_{\frac {6R}{10}}-\ex A_{\frac {9R}{10}}$) to get the ratio of the metric on the complement of $\ex A_{R,i}$ divided by the metric on  $\ex A_{i}$ is  bounded by $c'$ here.

 Then, we can use Lemma \ref{annulus lemma} and our conformal bound from condition \ref{conformal bound} to get some lower bound $r_{0}>0$ so that the boundary of $\ex A_{\frac {9R}{10}}$ and $\ex A_{\frac {6R}{10}}$ is further than $r_{0}$ from the boundary of the complement of $\ex A_{R,i}$. 
 A second easy compactness argument then tells us that if we fix any distance $r_{0}>0$, there exists some constant  $c_{r_{0}}$ (only depending on the bounds in the above paragraph) so that for every point further than $r_{0}$ from the boundary there exists an injective holomorphic map of the unit disk sending $0$ to that point, with derivative at $0$ greater than $c_{r_{0}}$.
  We then get that the above ratio of metrics above is bounded below by $\frac {c_{r_{0}}}{c'}$.


We must now prove our derivative bound on the complement of $\ex A_{\frac{9R}{10}}$.  
   We shall do this by proving that given any holomorphic injection of the unit disk $i:D\longrightarrow\ex C-\bigcup \ex A_{R,i}$ so that $i(0)\notin\ex A_{\frac {9R}{10},i}$, then $\abs{d(f\circ i)}$ is bounded. Then, our conformal bound from condition \ref{conformal bound} will tell us that we have a bound on $\abs {df}$ in the metric that we have chosen.

We do this in two cases. First, suppose that $i(0)$ is in the complement of $\ex A_{\frac R4,i}$ for all our old annuli. Then, Lemma \ref{annulus lemma} tells us that the image of $i$ restricted to the disk of radius $e^{-16\frac{\pi^2}{R}}$ is contained in the complement $\ex A_{\frac R2,i}$ for all our old annuli, as argued above. If the derivative of $f\circ i$ on this disk was large enough, we would be able to add another annulus to our old collection in the manner described above. As this process terminated, $\abs {d(f\circ i)}$ is bounded. Second, suppose that $i(0)$ is contained inside $\ex A_{\frac R4,i}$ for one of our old annuli. Then Lemma \ref{annulus lemma} tells us that the restriction of $i$ to the disk of radius $e^{-16\frac{\pi^2}R}$ must be contained inside $\ex A_i$, and therefore have energy less than $\epsilon$. This means that $\abs{d(f\circ i)}\leq c_0e^{16\frac{\pi^2}R}$, and our derivative is bounded as required.

   \stop

 \section{Compactness}

 \begin{thm}\label{completeness theorem}
  Given a basic, complete \exploded fibration $\ex B$ with a civilized almost complex structure $J$ and strict taming $\Omega$, the moduli stack $\mathcal M_{g,n,E}(\ex B)$ of stable holomorphic curves with a fixed genus $g$ and number of punctures $n$, and with $\Omega$ energy less than $E$ is topologically compact in $C^{\infty,\delta}$ for any $0<\delta<1$. 
  
  More generally, if $(\hat{\ex B}\longrightarrow\ex G,J,\Omega)$ is a family of such $(\ex B,J,\Omega)$, the map $\mathcal M_{g,n,E}(\hat{\ex B}\rightarrow\ex G)\longrightarrow\underline{\ex G}$ is $C^{\infty,\delta}$ topologically proper.
  \end{thm}
  
  In particular, this means that   given any sequence of the above holomorphic curves in the fibers over a topologically convergent sequence in $\ex G$,  there exists a subsequence $f^i$ which converges to a holomorphic curve $f$ as follows:  There exists a
  sequence of  families of smooth curves,
  
  \[\begin{array}{ccc} (\ex {\hat C},j_i) & \xrightarrow{\hat f^i} &\hat{\ex B}
  \\ \downarrow & &\downarrow
  \\ \ex F& \longrightarrow & \ex G\end{array}\]
 
  so that this sequence of families converges in $C^{\infty,\delta}$ to the smooth family
  
  \[\begin{array}{ccc} (\ex {\hat C},j) & \xrightarrow{\hat f} &\hat{\ex B}
  \\ \downarrow & &\downarrow
  \\ \ex F& \longrightarrow & \ex G\end{array}\]

  and a sequence of points $p^i\rightarrow\ex F$ so that $\totl{p^{i}}\rightarrow\totl{ p}$, $f^i$ is the map given by the restriction of $\hat f^i$ to the fiber over $p^i$, and $f$ is given by the restriction of $f$ to the fiber over $p$.  Of course, the case where we just have a single $(\ex B,J,\omega)$ is the same as the family case when $\ex G$ is a point.
  
  \
  
  The definition of $C^{\infty,\delta}$ convergence can be found in section \ref{regularity} starting on page \pageref{regularity} and the definition of the moduli stack is contained in section \ref{perturbation theory}, which starts on page \pageref{perturbation theory}. The notion of topological convergence is introduced on page \pageref{topological convergence}.
  
  \
  
  The proof of this theorem uses  Lemma \ref{strong cylinder convergence} on page \pageref{strong cylinder convergence} and  Proposition \ref{decomposition proposition} on page \pageref{decomposition proposition}. Together, these allow us to decompose holomorphic curves into pieces with bounded behavior. A standard Arzela-Ascoli type argument gets a type of convergence of these pieces. If we were working in the category of smooth manifolds, this would be sufficient to prove the compactness theorem because the type of convergence involved would have a unique limit, and therefore the limiting pieces would glue together. The extra problem that must be dealt with in this case is that we are not dealing with a type of convergence that has a unique limit, so we must work much harder to show that our resulting `limiting family' for each piece glues together to the limiting family $\ex F$. 
  
  \begin{example}
  \end{example}
  The following are examples of the types of non uniqueness we have to deal with:
  \begin{enumerate}
  \item \label{nonunique1}In $\et 11$, let $\tilde z(p^{i})=\frac 1i\e 0$. Then $\totl{p^{i}}\rightarrow \totl{p}$ where $p\rightarrow\et 11$ is any point so that $\tilde z(p)=c\e a$ where $a>0$. Note that there is no `point' $p\rightarrow\et 11$ so that $\tilde z(p)=0\e 0$. (The author did try modifying the definitions of exploded fibrations so that such limits exist, however if this is tried, the theory of holomorphic curves becomes much more complicated, and the whole setup is a lot less natural.) This is one reason that non unique limits need to be considered.  
  \item Expanding example \ref{nonunique1}, consider the maps $f^{i}:\mathbb CP^{1}\longrightarrow \mathbb CP^{1}\times \et 11$ given by $f^{i}(u)=(u,p_{i})$. This sequence of maps does not have a unique limit. We instead consider the family 
  \[\begin{split}
  \mathbb CP^{1}&\times \et 11\xrightarrow{\hat f:=\id}\mathbb CP^{1}\times \et 11
  \\ &\downarrow
  \\ &\et 11
 \end{split}\]
 Then the sequence of maps $ f^{i}$ are given by the restriction of the inverse image of $\hat f$ to the inverse image of $p^{i}$. These converge topologically to $\hat f$ restricted to the inverse image of $p$, where $p\rightarrow \et 11$ is any point so that $\totl {p^{i}}\rightarrow\totl p$.
 \item \label{nonunique3} Expanding example \ref{nonunique1} to see different behavior, consider the map $\pi:\et 22\longrightarrow\et 11$ so that $\pi^{*}\tilde z=\tilde w_{1}\tilde w_{2}$. This was studied in example \ref{local family model} on page \pageref{local family model}. This can be considered to be a family of annuli by restricting to the subset where $\abs{\tilde w_{1}}<1$ and $\abs{\tilde w_{2}}<1$. Then $\pi^{-1}(p^{i})$ is an annulus $1>\abs{\tilde w_{1}}>\frac 1 i$, or $\frac 1 i<\abs{\tilde w_{2}}< 1$. These domains converge topologically to $\pi^{-1}(p)$ where $p=c\e a$ for any $c\in \mathbb C^{*}$ and $a>0$. $\pi^{-1}(c\e a)$ has two coordinate charts which are subsets of $\et 11$, with coordinate $\tilde w_{i}$ so that 
 $\abs{\tilde w_{i}}<1 $, and $\totb{\tilde w_{i}}>\e a$. The transition map between these coordinate charts is given by 
 $w_{1}= c\e a w_{2}^{-1}$. 
 
 \item Let us  expand on example \ref{nonunique3} by adding in the information of maps $f^{i}$ from our domains $\pi^{-1}(p^{i})$ to $\ex T^{2}$. Suppose that 
 \[f^{i}(\tilde w_{1})=(\tilde c_{1}^{i}\tilde w_{1},\tilde c_{2}^{i}\tilde w_{1})\]
 so of course, 
 \[f^{i}(\tilde w_{2})=(\tilde c_{1}^{i}\frac 1i \tilde w_{2}^{-1},\tilde c_{2}^{i}\frac 1 i\tilde w_{2}^{-1})\]
  
  If $(c^{i}_{1},c_{2}^{i})\in \ex T^{2}$ are generically chosen, we will not have a unique limit even if we restrict to the bounded domain $1>\abs{\tilde w_{1}}>\frac 12$. To put all our $f^{i}$ into an individual family, consider the family
  \[\begin{split}\ex T^{2}&\times \et 11\xrightarrow{\hat f(\tilde c_{1},\tilde c_{2},\tilde w_{1},\tilde w_{2}):=(\tilde c_{1}\tilde w_{1},\tilde c_{2}\tilde w_{1})} \ex T^{2}
\\&  \downarrow\id\times \pi
 \\\ex T^{2}&\times \et 11
  \end{split}\] 
  Our individual maps $f^{i}$ are the restriction of $\hat f$ to the fiber over $(\tilde c^{i}_{1},\tilde c^{i}_{2},p^{i})$.
 These converge topologically to the restriction of $\hat f$ to the fiber over any point $(\tilde c_{1},\tilde c_{2},c\e a)$ where $a>0$. Note that if we restricted to the domains where $1>\abs{\tilde w_{1}}>\abs{\tilde w_{2}}$ or $1>\abs{\tilde w_{2}}>\abs{\tilde w_{1}}$, then we get a $\ex T^{2}$ worth of valid topological limits for each domain. There is a non trivial requirement that needs to be satisfied for these topological limits to be glued together.
 \end{enumerate}
  
  Apart from an easy use of an Arzela Ascoli type argument and elliptic bootstrapping, the main hurdle to proving Theorem \ref{completeness theorem} is dealing with the above types of non uniqueness for limits in showing that the different limiting pieces of our holomorphic curves glue together. This is not extremely difficult, but it requires a lot of notation to keep track of everything.

 Note the assumption that $\ex B$ is basic is mainly for convenience in the following arguments. In the family case, by restricting to a subset of $\ex G$ which contains the image of a subsequence, we can also assume that $\hat{\ex B}$ is basic. 

The following lemma will give us good coordinate charts on $\ex B$ or $\hat{\ex B}$.
Recall that if $\ex B$ is basic, we use the notation $\totb{\ex B_{i}}$ for a strata of the tropical part $\totb{\ex B}$, 
and $\overline{\totb{\ex B_{i}}}$ for the polygon  which after identifying some strata is equal to the  closure  of  $\totb{\ex B_{i}}\subset \totb{\ex B}$. A neighborhood of a strata $\ex B_{i}$ is then equal to an open subset of $\totl{\ex B_{i}}\rtimes\et n{\overline{\totb{\ex B_{i}}}}$ using the construction of example \ref{polygon2} on page \pageref{polygon2}. This has a (sometimes defined) action of $\ex T^{n}$ action corresponding to the (sometimes defined) action of $\ex T^{n}$ on $\et n{\overline{\totb{\ex B_{i}}}}$ given by coordinate wise multiplication. 
We will need this $\ex T^{n}$ action to move parts of our holomorphic curves around. 

\begin{lemma}\label{strata charts}
Given any basic \exploded fibration $\hat {\ex B}$, and a family $\hat{\ex B}\longrightarrow\ex G$,  for each strata $\hat{\ex B}_i\subset\hat{\ex B}$, there exists some $\ex U_i\subset\hat{ \ex B}$ containing $\hat{\ex B}_i$ so that
\begin{enumerate}
\item The image in the smooth part $\totl{\ex U_i}\subset\totl{\ex B}$ is an open neighborhood of $\totl{\ex B_{i}}\subset\totl{\ex B}$.
\item \label{strata chart normal form} If $\totb{\hat{\ex B}_i}$ is $n$ dimensional, then there is an identification of $\ex U_{i}$ with an open subset of $\totl{\ex B_{i}}\rtimes \et n{\overline{\totb{\ex B_{i}}}}$ using the construction on page \pageref{polygon2}. We can make this identification so that the (sometimes defined) free action of $\ex T^{n}$ given by considering $\totl{\ex B_{i}}\rtimes \et n{\overline{\totb{\ex B_{i}}}}\subset \totl{\ex B_{i}}\rtimes \ex T^{n}$ satisfies the following:    
\begin{enumerate}
\item The $\ex T^{n}$ action preserves fibers of the family $\hat{\ex B}\longrightarrow\ex G$, in the sense that if $p_{1}$ and $p_{2}$ have the same image in $\ex G$, then $\tilde z*p_{1}$ and $\tilde z*p_{2}$ also have the same image in $\ex G$.
\item For any point $p\longrightarrow\ex U_i$, the set of values in $\ex T^n$ for which this action is defined is  nonempty and convex in the sense that it is given by a set of inequalities of the form 

\[\abs{\tilde z^{\alpha}}>c\e x\in \mathbb R\e{\mathbb R} \] 
where $c\e x$ depends smoothly on $p\longrightarrow \hat{\ex B}$.

\item If $\totb{\hat{\ex B}_i}$ is in the closure of $\totb {\hat{\ex B}_j}$, then the action of $\ex T^n$ on $\ex U_i\cap\ex U_j$ is equal to the action of a subgroup of the $\ex T^m$ acting on $\ex U_j$. 

\end{enumerate}

\end{enumerate}    
\end{lemma}

 \pf
 
 The only point in the above lemma that requires proof is the compatibility of the actions of $\ex T^n$ with restriction and the family. We can construct this first for the strata $\totb{\hat{\ex B}_i}$ with the highest dimension and then extend it to lower dimensional strata.

 Being able to do this amounts to constructing coordinate charts identified with open subsets of $\mathbb R^k\times\et nB$, where coordinates on the projection of this to $\ex G$ consist of some sub collection of  these coordinates (the normal form for coordinates on a family), and 
 so that transition functions are of the form
 \[(x,\tilde w)\mapsto (\phi (x),f(x)\e a\tilde w)\]
 where $\phi$ is a diffeomorphism, and $f$ is a smooth $\mathbb C^*$ valued function. These transition functions preserve a $\ex T^n$ action of the form
 \[\tilde z*(x,\tilde w)=(x,\tilde z\tilde w)\]
 Because of the assumption that the closure of strata in $\totb{\hat{\ex B}}$ are simply connected, we can first reduce to the case that coordinate charts are in the normal form for a family, and  transition functions are of the form
 \[(x,\tilde w)\mapsto (\phi(x),f(x, w)\e a\tilde w)\]
 where  $f$ is smooth and $(\mathbb C^*)^n$ valued. There is then no obstruction to modifying our coordinate charts one by one so that the transition functions no longer have any dependance on $w$. (This amounts to replacing the coordinates $\tilde w$ with $\frac {f(x,w)}{f(x,0)} \tilde w$. This change is well defined on the intersection with previously corrected coordinate charts, does not affect our charts being in the normal form for families,  and there is no obstruction to extending it to the rest of a chart.) 

  \stop

 \
 
 We shall now start working towards proving Theorem \ref{completeness theorem}. The proof shall rely on Lemma \ref{strong cylinder convergence} on page \pageref{strong cylinder convergence} and  Proposition \ref{decomposition proposition} on page \pageref{decomposition proposition}. We shall be using the notation from Proposition \ref{decomposition proposition}.
 First, choose the following:
 
 \begin{enumerate}
 \item Choose a metric, a taming form $\omega\in\Omega$, and finite collection of coordinate charts on ${\ex B}$ each contained in some $\ex U_i$ from Lemma \ref{strata charts} and small enough to apply Lemma \ref{strong cylinder convergence}. For the case of a family, restrict $\hat{\ex B}\longrightarrow\ex G$ to some small coordinate chart on $\ex G$  in which the image of some subsequence converges, and choose our finite collection of coordinate charts  on (the smaller, renamed,) $\hat {\ex B}$ satisfying the above. 
 \item 
  Choose  exploded annuli $\ex A_n^i\subset \ex C^i$ satisfying the conditions in Proposition \ref{decomposition proposition} with $\omega$-energy bound $\epsilon$ small enough and $R$ large enough that each $\ex A^i_{\frac {6R}{10},n}$ is contained well inside some $\ex U_i$ in the sense that it is of distance greater than $2$ to the boundary of $\ex U_{i}$. (The fact that we can achieve this follows from Lemma \ref{omega cylinder bound} on page \pageref{omega cylinder bound}.) Also choose $\epsilon$ small enough and $R$ large enough so that each connected smooth component of $\ex A^i_{\frac {6R}{10},n}$ is contained in one of the above coordinate charts and can have Lemma \ref{strong cylinder convergence} applied to it.
\end{enumerate}

Our taming form $\omega$ is a smooth two form on some refinement of $\hat{\ex B}$. As convergence in a refinement is stronger, we will simply call this refinement $\hat{\ex B}$.

 Use the notation $C^i_m$ to indicate connected components of $\ex C^i-\bigcup_k \ex A^i_{\frac {8R}{10}, k}$. We can choose a subsequence so that the number of such components is the same for each $\ex C^i$, $C^i_m$ has topology that is independent of $i$, and choosing diffeomorphisms identifying them, $C^i_m$ converges as $i\rightarrow \infty$ to some $C_m$ in the sense that the metric from Proposition \ref{decomposition proposition} and the complex structure converges to one on $C_m$. 
 
 \psfrag{ETF16a}{$C^{i}_{1}=(C_{1},j_{i})$}
 \psfrag{ETF16b}{$(C_{2},j_{i})$}
 \psfrag{ETF16c}{$(C_{3},j_{i})$}
 
 \psfrag{ETF16d}{$\ex A^{i}_{4}$}
 \psfrag{ETF16e}{$\ex A^{i}_{5}$}
 \psfrag{ETF16f}{$\ex A^{i}_{6}$}
 
 \psfrag{ETF16g}{Transition functions bounded}
  \includegraphics{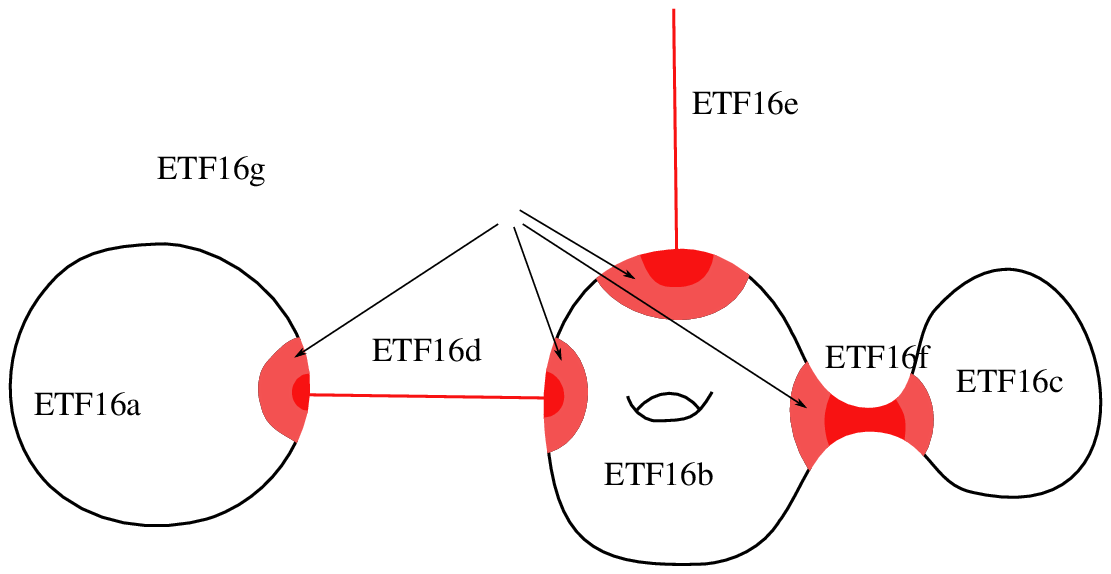}
 
  We shall choose two coordinate charts on $\ex A_n^i$ with coordinates $\tilde z_n^{\pm}$ in a subset of $\et 1{1}$ so that each boundary of $\ex A_n^i$ is identified with $\abs {\tilde z_n^\pm}=1$, and the coordinates are related by $\tilde z_n^+\tilde z_n^-=Q^{i}_{n}\in\mathbb R^{*}\e {\mathbb R^{+}}$. (This means that the conformal modulus of $\ex A_n^i$ is $-\log Q^{i}_{n}$. These $Q_{n}$ will later become coordinates on our family $\ex F$.) Choose a subsequence so that the number of annuli is independent of  $i$, and  $\{Q^{i}_{n}\}$ converges topologically. (In other words, calling any  conformal modulus not a real number infinity, either the conformal modulus of $\ex A_n^i$ converges to a finite number, or converges to infinity.  )
  
  We can define 
 \[\ex A^{\pm}_{n}:=\{\abs {\tilde z^{\pm}_{n}}\leq e^{-\frac{7R}{10}}\}\subset\et 11\] 
 \[\ex A^{i}_{\frac {7R}{10},n}:=\{(\tilde z^{+}_{n},\tilde z^{-}_{n})\text{ so that }\tilde z^{+}_{n}\tilde z^{-}_{n}=Q^{i}_{n}\}\subset\ex A^+_{n}\times\ex A^{-}_{n}\]
 We shall need to keep track of what $\ex A_{n}$ is attached to, so use the notation $n^{\pm}$ to define $C^{i}_{n^{\pm}}$ as the component attached to the end of $\ex A_{n}^{i}$ with the boundary $\abs{\tilde z^{\pm}_{n}}=1$. (Assume, by passing to a subsequence, that $n^{\pm}$ is well defined independent of $i$.)

 We can consider transition maps between $C^{i}_{n^{\pm}}$ and $\ex A^{i}_{\frac {7R}{10},n}$ to give transition maps between $C_{n^{\pm}}$ and $\ex A^{\pm}_{n}$. 
 The bound on the derivative of transition maps from Proposition \ref{decomposition proposition} on the region $\ex A^{i}_{\frac {9R}{10},n}-\ex A^{i}_{\frac {6R}{10},n}$ (and standard elliptic bootstrapping to get bounds on higher derivatives on the smaller region $\ex A^{i}_{\frac {8R}{10},n}-\ex A^{i}_{\frac {7R}{10},n}$) tells us that we can choose a subsequence so that the transition maps between $C^i_{n^{\pm}}$ and $\ex A^\pm_{n}$ converge to some smooth transition map between $C_{n^{\pm}}$ and $\ex A^{\pm}_{n}$. We can modify our diffeomorphisms identifying $C_m^i$ with $C_m$ so that these transition maps all give exactly the same map. Denote the complex structure on $C_m$ induced by this identification by $j_i$. (We shall do something similar for the annuli $\ex A^{i}_{n}$ later on.) We can do this so that $j_i$ and the metrics given by this identification still converge in $C^\infty$ to those on $C_m$.

 \psfrag{ETF17a}{$(C_{1},j)$}
 \psfrag{ETF17a'}{$q_{1}$}
 \psfrag{ETF17b}{$(C_{2},j)$}
 \psfrag{ETF17b'}{$q_{2}$}
 \psfrag{ETF17c}{$(C_{3},j)$}
 \psfrag{ETF17c'}{$q_{3}$}
 \psfrag{ETF17d+}{$\ex A^{+}_{4}$}
\psfrag{ETF17d+a}{$4^{+}:=1$}
 \psfrag{ETF17d-}{$\ex A^{-}_{4}$}
 \psfrag{ETF17d-b}{$4^{-}:=2$}
 \psfrag{ETF17e}{$\ex A^{+}_{5}$}
 \psfrag{ETF17eb}{$5^{+}:=2$}
 \psfrag{ETF17f+}{$\ex A_{6}^{+}$}
 \psfrag{ETF17f-}{$\ex A_{6}^{-}$}
\psfrag{ETF17f+c}{$6^{+}:=3$}
\psfrag{ETF17f-b}{$6^{-}:=2$}
 \includegraphics{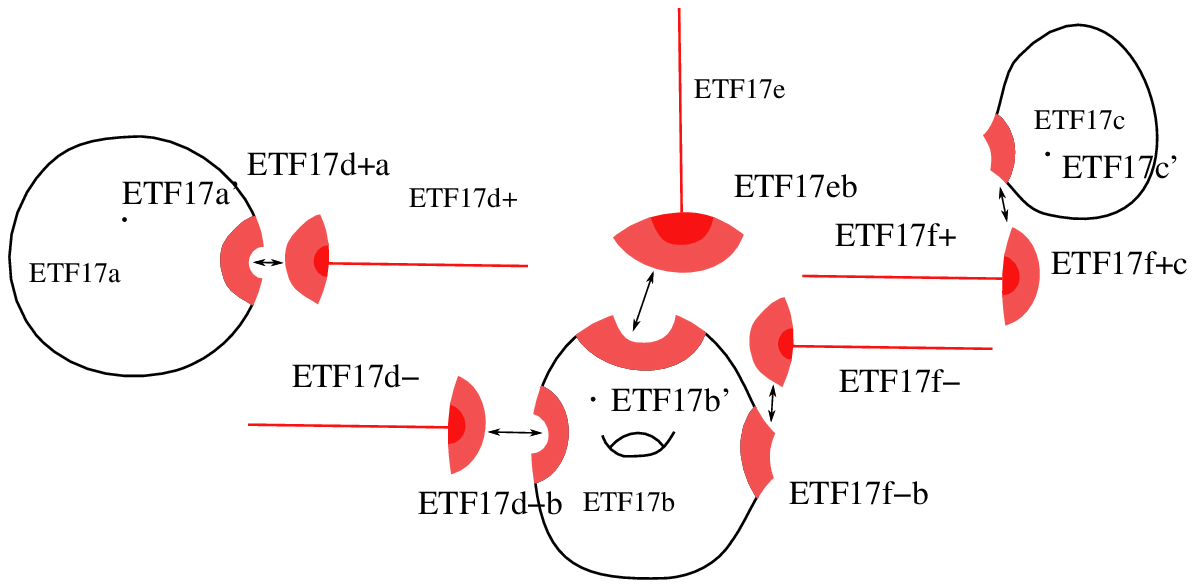}
 
 We shall now start to construct our family $\ex {\hat C}\longrightarrow \ex F$. The first step shall be to get some kind of convergence for each piece of our holomorphic curve.
 
 \

 For each $C_m$ choose some point $q_m\in  C_m$, and consider $Q^i_m:=f^i(q_m)\in \hat{\ex B}$. This sequence has a subsequence that converges in $\totl{\hat{\ex  B}}$ to some point $Q_m\longrightarrow\hat{ \ex B}$. Label the strata of $\hat{\ex B}$ that contains $Q_m$ by $\hat{\ex B}_m$, and consider the chart $\ex U_m$ from Lemma \ref{strata charts} containing $\hat{\ex B}_m$. As the derivative of $f^i$ restricted to $C_m$ is bounded, we can choose a subsequence so that $f^i(C_m)$ is contained well inside $\ex U_m$. (We can choose a subsequence so that the image of $f^{i}(C_{m})$ is contained in the subset of $\ex U_{m}$ which consists of all points some arbitrary distance from the boundary of $\ex U_{m}$.)

 As the derivative of $f^i$ is uniformly bounded on $\ex C^i-\bigcup \ex A^i_{\frac {9R} {10},k}$, we can use lemmas \ref{coordinate bound}, \ref{dbar of derivative}, and \ref{elliptic regularity} to get bounds on the higher derivatives of $f^i$ restricted to $C^i_m$. Lemma \ref{bounded geometry} then tells us that if $f^i(q_m)$ converges topologically to $Q_m\longrightarrow \hat{\ex B}$, the geometry around $f^i(q_m)$ converges to that around $Q_m$, so we can choose a subsequence so that $f^i$ restricted to $ C^i_m$ converges in some sense to a map  
 \[f_{m,Q_m}: C_m\longrightarrow\hat{\ex B}\text{ so that }f_{m,Q_m}(q_m)=Q_m\subset\ex U_{m}\]  
 More specifically, remembering that  everything is contained inside $\ex U_m$, we can use our $\ex T^{n}$ action on $\ex U_{m}$ to say this more precisely.  There exists some sequence $\tilde c_m^i\in\ex T^n$ so that
 \[\tilde c_m^i*f^i:C_m\longrightarrow \ex U_m\]
 converges in $C^\infty$ to $f_{m,Q_m}:C_m\longrightarrow\ex U_m$. 
 
 Of course, there was a choice involved here. We could also have chosen a different topological limit $Q'_{m}$  of $f^{i}(q_{m})$, which would give a translate of $f_{m,Q_{m}}$ by our $\ex T^{n}$ action. These choices will turn up as  parameters on our family. They will need to be `compatible' in the sense that these pieces will need to fit together.

 \
  
We now consider the analogous convergence on annular regions. Because the $\omega$-energy of $f^{i}$ restricted to our annular regions $\ex A^i_n$ is small, we can apply Lemma \ref{strong cylinder convergence} on page \pageref{strong cylinder convergence} to tell us that if the limit of the conformal modulus of $\ex A^i_n$ is infinite, then $f^i$ restricted to the smooth parts of $\ex A^i_n$ converges in some sense (considered in more detail later) to some unique pair of holomorphic maps
  \[f^{\pm}_{n,Q}:\{e^{-\frac {7R}{10}}>\abs {\tilde z_{n}^\pm}>0\}\longrightarrow\hat{\ex B}\] compatible with $f_{n^{\pm},Q_{n^{\pm}}}$ and the transition maps. (Recall that $C_{n^{\pm}}$ is the component attached to $\ex A^{\pm}_{n}$.)  As we have made the assumption that $J$ is civilized, we can use the usual removable singularity theorem for finite energy holomorphic curves to see that these limit maps extend uniquely to  smooth maps on $\{e^{-\frac {7R}{10}}>\abs {\tilde z_{n}^\pm}>\e x\}\subset\ex A^{\pm}_{n}$ for some $x>0$. 
  
  \

  \psfrag{ETF18a}{$C_{4^{+}}:=C_{1}$}
  \psfrag{ETF18a'}{$q_{1}$}
  \psfrag{ETF18aa}{${f_{1,Q_{1}}}: C_{1}\longrightarrow\ex U_{1}$}
  \psfrag{ETF18ab}{$q_{1}\mapsto Q_{1}\in\ex U_{1}$}
  
  \psfrag{ETF18d+}{$\xrightarrow{f^{+}_{4,Q}}\ex U_{4}$}
  \psfrag{ETF18d-}{$\ex U_{4}\xleftarrow{f^{-}_{4,Q}}$}
  \psfrag{ETF18dd}{$\tilde z_{4}^{+}\tilde z_{4}^{-}=Q_{4}$}
 
   \psfrag{ETF18ba}{$C_{2}\xrightarrow{f_{1,Q_{2}}}\ex U_{2}$}
  \psfrag{ETF18bb}{$q_{2}\mapsto Q_{2}\in\ex U_{2}$}

  \psfrag{ETF18b}{$C_{4^{-}}:=C_{2}$}
  \psfrag{ETF18b'}{$q_{2}$}
  
  \includegraphics{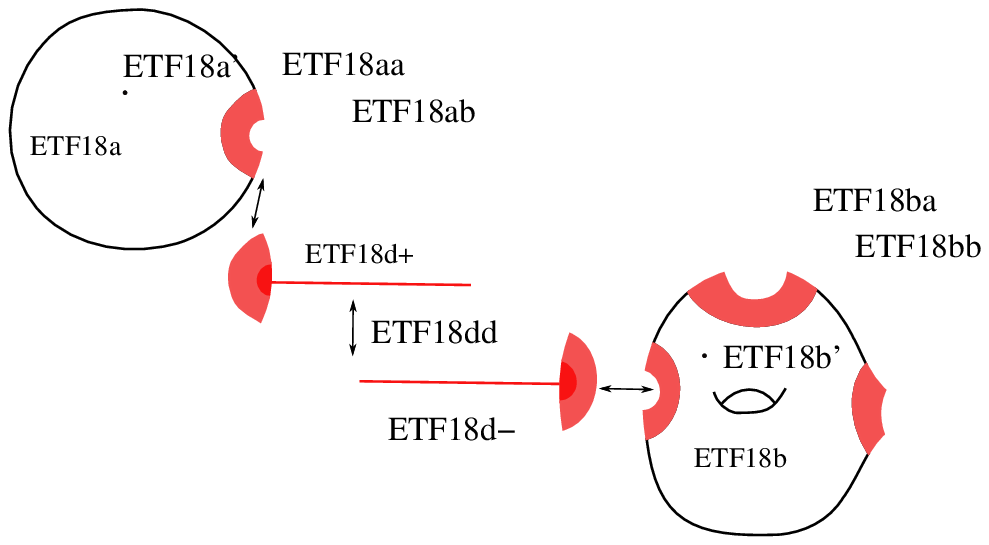}
  
  \
  
   We may assume after passing to a subsequence that $f^i(\ex A^i_{\frac {6R}{10},n})$ are all contained inside a single $\ex U_n$. Recalling our convention that $C_{n^{\pm}}$ intersects $\ex A^i_n$, note that this means  that $\ex U_{n^{\pm}}$ intersects $\ex U_n$ (and $\totb{\ex U_{n^{\pm}}}\subset\totb{\ex U_n}$), so our sequence $\tilde c^i_{n^{\pm}}$ defined earlier also has an action on $\ex U_{n}$, and  
   \[\tilde c^i_{n^{\pm}}* f^i: \ex A_{\frac {7R}{10},n}^i\longrightarrow \ex U_n\] converges in $C^\infty$ on compact subsets adjacent to $C_{n^{\pm}}$ to  $f^\pm_{n,Q}$ . 
   
   \
   
 We now have a type of convergence of the individual pieces we have cut our holomorphic curve into. We shall now define a model for the exploded structure on our family that is too large, as it ignores the requirement that these pieces must fit together.
  
\

 Define $\ex V_m \subset \ex U_m$ to be the exploded fibration consisting of all points $Q_m'\longrightarrow\ex U_m$ so that
 
 \begin{enumerate}
 \item There exists some $\tilde z$ so that $Q_m'=\tilde z*Q_m$ 
 \item  Defining 
 \[f_{m,Q_m'}:=\tilde z*f_{m,Q_m}\]
 The image of $C_{m}$, $f_{m,Q'_{m}}(C_{m})$ is contained well inside $\ex U_{m}$ in the sense that the distance to the boundary of $\ex U_{m}$ is greater than $1$.
 \item If $C_{m}$ is attached to $\ex A^{\pm}_{n}$, then defining
 \[f^{\pm}_{n,Q'}:=\tilde z*f^\pm_{n,Q}\]
The image of the smooth part of $\ex A^{\pm}_{n}$, $\{f^{\pm}_{n,Q'}(\tilde z),\ 0<\abs{\tilde z}<e^{-\frac {7R}{10}}\}$ is also contained well inside $\ex U_{n}$.
  \end{enumerate}

 Note that $\ex V_{m}$ is a smooth exploded fibration which includes all $Q'_m$ which are topological limits for $Q^i_m$. For such $Q'_m$, the map $f_{m,Q'_m}$ will be holomorphic, but for other $Q'_m$, this may not be the case. 

We shall consider the family $\ex F$ as a sub exploded fibration of  

\[\prod_m\ex V_m\prod_n\ex Q_n\subset \prod_m\ex U_m\prod_n\ex Q_n\]

The $\ex Q_n$ above stands for the `gluing parameter' for identifying coordinates $\tilde z_{n}^\pm$ on $\ex A^\pm_n$ via 
\[\tilde z_{n}^+\tilde z_{n}^-=Q_n\in\ex Q_n:=\{\abs{\tilde z}<e^{-2R}\}\subset\et 11\] 
   
   \
   
 We shall consider the following sequence of points \[Q^i:=(f^i(q_1),f^i(q_2),\dotsc,Q^i_{n_1},\dotsc)\in\prod_m \ex U_m\prod_n\ex Q_n\]
  where the conformal modulus of $\ex A^i_{n}$ is equal to $-\log Q^i_{n}$.
 
 \
 
 Note that there is a transitive (sometimes defined) action of $\ex T^k$ on $
 \prod \ex V_m\prod \ex Q_n$ which is the action from Lemma \ref{strata charts} on $\ex V_m\subset \ex U_m$, and multiplication by some coordinate of $\ex T^k$ on $\ex Q_n$. Our family will be given by a complete inclusion
 \[\ex F\longrightarrow \prod_m\ex V_m\prod_n\ex Q_n\subset \prod_m\ex U_m\prod_n\ex Q_n\]
  satisfying the following conditions:
  \begin{enumerate}
  \item \label{F condition 1}The image of $\ex F$ contains some point $Q\longrightarrow \prod_m\ex V_m\prod_n\ex Q_n $ which is a topological limit of $Q^i$, and the image of $\ex F$ is given by the orbit of $Q$ under the action of some subgroup of $\ex T^k$.
  \item\label{F condition 2} The distance in any smooth metric from some subsequence of $Q^i$ to the image of $\ex  F$ converges to $0$.
  \item There is no other inclusion satisfying the above conditions which has smaller dimension than $\ex F$.
  \end{enumerate}
  
  It is clear that such an $\ex F$ must exist. It can be seen from Lemma \ref{strata charts}, item \ref{strata chart normal form} that $\ex F$ is smooth. 

\

Now let us construct $\hat{\ex C}\longrightarrow\ex F$. We shall have charts on $\hat{\ex C}$ given by $\hat C_m:=C_m\times\ex F$ and $\hat{\ex A}_n$. $\hat {\ex A}_n$ has coordinates $(\tilde z_{n}^+,\tilde z_{n}^-,Q)$ where $\tilde z^{\pm}_{n}\in\ex A^{\pm}_{n}$ and $Q_n=\tilde z_{n}^+\tilde z_{n}^-$. Transition maps between $\hat C_m$ and $\hat{\ex A}_n$ are simply given by the transition maps between $C_m$ and $\ex A^\pm_{n}$ times the identity on the $\ex F$ component. As an explicit example, if $\phi$ is a transition map between $C_{m}$ and $\ex A^{+}_{n}$, the corresponding transition map between $\hat C_{m}$ and $\hat{\ex A}_{n}$ is given by 
\[(z,Q)\mapsto\left(\phi(z),\frac{Q_{n}}{\phi(z)},Q\right) \]
  This describes the exploded fibration $\hat{\ex C}$. The map down to $\ex F$ is simply given by the obvious projection to the second component of each of the above charts. This gives a smooth family of exploded curves.   
  Note that condition \ref{decomposition stability} from Proposition \ref{decomposition proposition} ensures that this is a family of stable curves.

\begin{lemma}\label{matching}
 If $Q\longrightarrow \ex F$, and $\totl{Q^{i}_{n}}\rightarrow \totl{Q_{n}}$, then $f^{\pm}_{n,Q}$ can be glued  by the identification $\tilde z_{n}^+\tilde z_{n}^-=Q_n$. This is automatic if the limit of the conformal modulus of $\ex A^{i}_{n}$ is finite.
\end{lemma}

\pf

   In the case that the limit of the conformal modulus of $\ex A^i_n$ is finite, $\lim Q^{i}_{n}=Q_{n}$ for any $Q\longrightarrow \ex F$. If this was not true, we could simply restrict $\ex F$ to the points that satisfy this. This would include any point in $\ex F$ that is the topological limit of $Q^{i}$, it would also satisfy the other conditions required of $\ex F$, but have smaller dimension contradicting the minimality of $\ex F$.
    When we have this, as we have $f^{i}$ converging on $C_{n^{\pm}}$ and $\ex A^{i}_{n}$, and transition maps between these also converging, we can glue the limit.

  In the case that the limit of the conformal modulus of $\ex A^i_n$ is infinite, we want to  glue $f^+_{n,Q}(\tilde z_{n}^+)$ to $f^{-}_{n,Q}(\tilde z_{n}^-)$ over the region where the smooth coordinates $ \totl{\tilde z^+_{n}}=\totl{\tilde z^-_{n}}=0$ via the identification $\tilde z_{n}^+\tilde z_{n}^-=Q_n$. Define the following continuous function 
    \[\phi_n:\ex F\longrightarrow \ex T^l \text{ (the group acting on $\ex U_n$)}\] which detects the failure for $f^\pm_n$ to glue for any $Q\longrightarrow \ex F$: If $\totl{Q_n}=0$,  and $\tilde z_{n}^+\tilde z_{n}^-=Q_n$ so that $\totl{\tilde z^+}=\totl{\tilde z^-}=0$, let $\phi(Q)$ be the element of $\ex T^l$   so that $\phi_n(Q)*f_{n,Q}^+(\tilde z_{n}^+)=f_{n,Q}^-(\tilde z_{n}^-)$. This may only fail to be defined because the left or right hand sides might not be inside $\ex U$.  Note first that if defined, this doesn't depend on our choice of $\tilde z_{n}^\pm$. Also note that this will be defined on some open subset of $\ex F$, and there exists some homomorphism from $\ex T^k$ to $\ex T^l$ so that this map is equivariant with respect to the $\ex T^k$ action on $\ex F$ and the $\ex T^l$ action on $\ex U_m$. We can uniquely extend $\phi_n$ to the rest of $\ex F$ in an equivariant way. 
    
    We want to show that $\phi_n$ is identically $1$ on $\ex F$. For each $Q^i$, there exists a close point $\check{Q}^i$ inside $\ex F$ so that the distance in any smooth metric between $\check Q^i$ and $Q^i$ converges to $0$ as $i\rightarrow\infty$. Then $\phi_n(\check{Q}^i)$ converges to $1$, because these points $Q^i$ come from holomorphic curves that are converging on either end of  $\ex A^i_{\frac {7R}{10},n}$ to $f^{\pm}_n$ in a way given by Lemma \ref{strong cylinder convergence}. This tells us that the sub exploded fibration of $\ex F$ given by $\phi_n^{-1}(1)$ satisfies the conditions \ref{F condition 1} and \ref{F condition 2} above, so by the minimality of the dimension of $\ex F$, it must be all of $\ex F$. This tells us that if $Q_n^{-1}$ is infinite, then we can glue together $f^{\pm}_n$ without any modifications. 
 
 \stop
 
 Note that choosing any point $Q\longrightarrow \ex F$ so that $Q$ is the topological limit of $Q^i$ (such a point must exist by the definition of $\ex F$), the above lemma allows us to glue together $f_{m,Q_m}$ and $f^{\pm}_{n,Q}$ to obtain our limiting holomorphic curve $f$ for Theorem \ref{completeness theorem}. If the above lemma held for every point in $\ex F$, we would have constructed our family. As it is, we need to make some gluing choices.
 
 \
 
 If $Q_n\in\mathbb C^*$ and the limit of the conformal modulus of $\ex A_n^i$ is infinite, we will need to make `gluing' choices.  We do not use the standard gluing and cutting maps as this will not give us strong enough regularity for the resulting family.

\begin{enumerate}

\item
Chose some sequence $\check Q^i\longrightarrow \ex F$ so that the distance between $\check Q^i$ and $Q^i$ converges to $0$. 

We need to take care of the different complex structures obtained by gluing $\ex A^{\pm}_{n}$ by $\check Q^{i}_{n}$ and $Q^{i}_{n}$. Choose an almost complex structure $j_{i}$ on $\ex A^{+}_{n}$ as follows:  Choose smooth isomorphisms
\[  \Phi^{i}_{n}:\ex A^{+}_{n}\longrightarrow \ex A^{+}_{n}\]
so that 
\begin{enumerate}
\item \[\Phi^{i}_{n}\lrb{\tilde z_{n}^{+}}=\tilde z_{n}^{+}\text{ for }\abs {\tilde z_{n}^{+}}\geq e^{-\frac {8R}{10}}\]
\item \[\Phi^{i}_{n}\lrb{\tilde z_{n}^{+}}=\frac{\check Q^{i}_{n}}{Q^{i}_{n}}\tilde z_{n}^{+}\text{ for }\abs{\tilde z_{n}^{+}}\leq e^{-R}\]
\item On the region $ \totl{\tilde z_{n}^{+}}\neq 0$, the sequence of maps $\{\Phi^{i}_{n}\}$ converges in $C^{\infty}$ to the identity.
\end{enumerate}

Now define $j_{i}$ on $\ex A^{+}_{n}$ to be the pullback under $\Phi^{i}_{n}$ of the standard complex structure. Using the standard complex structure on $\ex A^{-}_{n}$, we then get our $j_{i}$ defined on $\hat{\ex A}_{n}$. This is compatible with  the  $j_{i}$ already defined on  $\hat C_{m}$, so we get $j_{i}$ defined on $\hat{\ex C}$. Note that $j_{i}$ restricted to the fiber over $\check Q^{i}$ is the complex structure on $\ex C^{i}$.  From now on, we shall use these new coordinates on $\ex A^{i}_{n}$.

\item We now define the linear gluing map as follows:
\begin{enumerate}
\item Chose some smooth cutoff function 
\[\rho:\mathbb R^{*}\longrightarrow [0,1]\]
so that \[\rho(x)=0\text{ for all }x\geq e^{-\frac {8R}{10} }\]
\[\rho(x)=1\text{ for all }x\leq e^{-\frac{9R}{10}}\] 
Extend this to
\[\rho:\mathbb R^{*}\e{\mathbb R}\longrightarrow [0,1]\]
satisfying all the above conditions. (We defined this first on $\mathbb R^{*}$ so that it was clear what `smooth' meant.)
\item\label{cutnglue}

Given  maps $\phi^+,\phi^-:\et 11\longrightarrow \mathbb C^k$ which vanish at $z=0$ define the gluing map 

\[G_{(\phi^+,\phi^-)}(\tilde z^+,\tilde z^-):=\rho \left( \abs {\tilde z^-}\right)\phi^+(\tilde z^+)+\rho \left( \abs {\tilde z^+}\right)\phi^-(\tilde z^-)\] 

Note that if $\phi^+$ and $\phi^-$ are smooth, $G_{(\phi^+,\phi^-)}:\et 22\longrightarrow \mathbb C^k$ is smooth. Note also that if $\phi^{+}$ and $\phi^{-}$ are small in $C^{\infty,\delta}$, then $G_{(\phi^{+},\phi^{-})}$ is small too.

\end{enumerate}

\item We now define a linear `cutting' map: as follows:
\begin{enumerate}
\item Choose a smooth cutoff function $\beta:\mathbb R^{*}\longrightarrow[0,1]$ satisfying the following:
\[\beta(x)+\beta\lrb{x^{-1}}=1\]
\[\beta(x)=1\text{ for all } x>e^{\frac R{10}}\]

\item
Given a map $\phi:\ex A^{i}_{n}\longrightarrow \mathbb C^k$ and $\check Q^{i}_{n}\in\mathbb C^*$, where 
\[\ex A^{i}_{n}:=\{z^{+}z^{-}=\check Q^{i}_{n}\}\]
we can define the cutting of $\phi$ as follows:
\[\phi^+(z^+):=\beta \left(\frac{\abs{z^+}^2}{\abs {\check Q^{i}_{n}}}\right)\phi(z^+)\]
\[\phi^-(z^-):=\beta \left(\frac{\abs{z^-}^2}{\abs {\check Q^{i}_{n}}}\right)\phi(z^-):=
\beta \left(\frac{\abs{z^-}^2}{\abs {\check Q^{i}_{n}}}\right)\phi\left(z^{+}=\frac {\check Q^{i}_{n}}{z^-}\right)\]

Note that $G_{(\phi^+,\phi^-)}$ restricted to $\tilde z^+\tilde z^{-}=\check Q_{n}$ is equal to $\phi$.
\end{enumerate}

\item Lemma \ref{cylinder bound} tells us that the image of our annuli of finite conformal moduli, $f^i(\ex A^i_{\frac {7R}{10},n})$ is contained in some coordinate chart appropriate for Lemma \ref{strong cylinder convergence}. (We can choose this coordinate chart to be contained inside $\ex U_n$ so that the $\ex T^k$ action on $\ex U_n$ just consists of multiplying coordinate functions by a constant). In fact, we can choose a subsequence so that either $\ex A^i_n$ has infinite conformal modulus for all $i$ or $f^i(\ex A^i_{\frac {7R}{10},n})$ is contained in one of these coordinate charts for all $i$. In that case, Lemma \ref{strong cylinder convergence} together with the conditions on our coordinate change for $\ex A^{i}_{n}$ above tells us that 
\[\begin{split} f^i_n(z)&=(e^{\phi^i_{n,1}(z)}c_{n,1}^i z^{\alpha_1},\dotsc,e^{\phi^i_{n,k}(z)}c_{n,k}^i z^{\alpha_k},\dotsc,c_{n,d-k}+\phi_{n,d-k}^{i})
\\&:=e^{\phi^i_n} c^i_nz^\alpha\end{split}\]

where $\phi^i_n$ is exponentially small on the interior of $\ex A^i_{R,n}$ as required by Lemma \ref{strong cylinder convergence}.  We can regard our local coordinate chart as giving us a local trivialization of $T\ex B$ compatible with our $\ex T^k$ action. We can extend this trivialization using our $\ex T^k$ action to an open subset  $\ex O_n\subset \ex U_n$ which will contain the image of any relevant translation of $f^i_n$. Note that $\ex O_n$ is an open subset of some refinement of $\ex T^k\times\mathbb R^{d-2k}$ (with the above coordinates). We shall define our maps as maps to $\ex T^k\times\mathbb R^{d-2k}$. Note that in particular, if $Q^\infty\longrightarrow\ex F$ is some topological limit of $Q^i$, then $f^{\pm}_{n,Q^\infty}$ in these coordinates are given by 
\[f^+_{n,Q^\infty}(\tilde z^+):=e^{\phi^{+}_n(\tilde z^+)} c_n(\tilde z^+)^\alpha\]
\[:=(e^{\phi^{+}_{n,1}(\tilde z^+)}c_{n,1} (\tilde z^+)^{\alpha_1},\dotsc,e^{\phi^+_{n,k}(\tilde z^+)}c_{n,k} (\tilde z^+)^{\alpha_k},\dotsc,c_{n,d-k}+\phi^{+}_{n,d-k}(\tilde z^{+}) )
\]
\[f^-_{n,Q^\infty}(\tilde z^-):=e^{\phi^-_n(\tilde z^-)}c_n\left(\frac {Q^\infty_n}{\tilde z^-}\right)^\alpha\]
where $\phi^\pm$ are smooth and vanish at $z^\pm=0$.

Similarly, if $\ex A^i_n$ is infinite, then we can consider $f^i_n$ to map to $\ex O_n$ so that 

\[f^i_n(\tilde z^+):=e^{\phi^{i+}_n} c^i_n(\tilde z^+)^\alpha\]
\[f^i_{n}(\tilde z^-):=e^{\phi^{i-}_n(\tilde z^-)}c^i_n\left(\frac {\check Q^i_n}{\tilde z^-}\right)^\alpha\]

\item Use the notation $\frac {p_1}{p_2}$ for two points $p_i\longrightarrow \ex V_n$ to mean the element of $\ex T^k$ so that $p_1=\frac {p_1}{p_2}*p_2$.
Define $F^i_n:\hat{\ex A}_n\longrightarrow \hat{\ex B}$
by
\[F^i_n(\tilde z^+,\tilde z^-,Q):=\frac{Q_{n^{+}}}{\check Q^i_{n^{+}}}*c^i_n(\tilde z^+)^\alpha\]
 and
\[F_n(\tilde z^+,\tilde z^-,Q):=\frac{Q_{n^{+}}}{ Q^\infty_{n^{+}}}*c_n(\tilde z^+)^\alpha\]
 (Recall that $Q_{n^+}$ recorded the position of $C_{n^{+}}$ which is attached to  $\ex A^{+}_n$). Note that even though we defined this as a map to $\ex T^k\times \mathbb R^{d-2k}$, and our chart is some refinement of this,  these $F^i_n$ are still smooth maps to $\hat{\ex B}$. (This can be seen if we restrict this map to local coordinate charts.)


Then we can define $\hat f^i:\hat{\ex A}_n\longrightarrow \hat{\ex B}$ using the cutting and gluing maps above by

\[\hat f^i(\tilde z^+,\tilde z^-,Q):=e^{G_{(\phi^{i+}_n,\phi^{i-}_n)}(\tilde z^+,\tilde z^-)}F^i_n(\tilde z^+,\tilde z^-,Q)\]
We can similarly define $\hat f:\hat{\ex A}_n\longrightarrow \hat{\ex B}$.

Also define
\[\hat f:\hat{ C}_m\longrightarrow \hat{\ex B}\]
by
\[\hat f(z,Q)=f_{m,Q_m}(z)\]
and define $\hat f^i$ on $\hat C_m$ similarly by translating $f^i$ appropriately depending on the difference between $\check Q^i_{m}$ and $Q_{m}$.

This gives  well defined smooth maps
\[\hat f^i:\hat {\ex C}\longrightarrow \hat{\ex B}\]
\[\hat f:\hat{\ex C}\longrightarrow \hat{\ex B}\]
so that $f^i$ is given by the restriction of $\hat f^i$ to the fiber of $\hat{\ex C}\longrightarrow \ex F$ over $\check Q^i\longrightarrow\ex F$. Note that as the vectorfields $\phi^{\pm}$ are vertical with respect to $\hat{\ex B}\longrightarrow \ex G $, $ G_{\phi^{+},\phi^{-}}$ is too, and we have finally constructed our required families.

  \[\begin{array}{ccccccc} (\ex {\hat C},j_i) & \xrightarrow{\hat f^i} &\hat{\ex B} &\ \ \ \ \ \ \ \ \ &
  (\ex {\hat C},j) & \xrightarrow{\hat f} &\hat{\ex B}
  \\ \downarrow & &\downarrow& &\downarrow & &\downarrow
  \\ \ex F& \longrightarrow & \ex G & & \ex F &\longrightarrow &\ex G\end{array}\]

\end{enumerate}

\begin{lemma}
 \[\hat f^i:(\hat{\ex C},j_i)\longrightarrow \hat{\ex B}\] converges in $C^{\infty,\delta}$ to 
 \[\hat f:(\hat{\ex C},j)\longrightarrow \hat{\ex B}\] 
\end{lemma}

\pf

We work in coordinates. First note that $\hat f^i$ converges in $C^\infty$ restricted to $\hat f$ on $\hat C_m$. This also holds for $\hat {\ex A}_n$ if the size of $\ex A^i_n$ stays bounded. Note also that $j_i$ converges in $C^\infty$ to $j$.  

It remains only to show that $\hat f^i$ converges in $C^{\infty,\delta}$ to $\hat f$ on $\hat{\ex A}_n$. Recall that $\hat {\ex A}_n$ has coordinates $(\tilde z^+,\tilde z^-,Q)$ where $Q\longrightarrow \ex F$ and $Q_n=\tilde z^+\tilde z^-$. Note that the maps $F^i_n:\hat{\ex A}_n\longrightarrow\hat{\ex B}$ defined above converge in $C^\infty$ to $F_n:\hat{\ex A}_n\longrightarrow\hat{\ex B}$. We have that

\[\hat f^i(\tilde z^+,\tilde z^-,Q):=e^{G_{(\phi^{i+}_n,\phi^{i-}_n)}(\tilde z^+,\tilde z^-)}F^i_n(\tilde z^+,\tilde z^-,Q)\]
 Similarly, we write
\[\hat f(\tilde z^+,\tilde z^-,Q):= e^{G_{(\phi^+_n,\phi^-_n)}(\tilde z^+,\tilde z^-)}F_n(\tilde z^+,\tilde z^-,Q)\]

We can choose a subsequence so that  $\phi^{i\pm}_n( z^\pm)$ converges in $C^\infty$ on compact subsets of $\{0<\abs{z^{\pm}}\leq e^{-\frac {7R}{10}}\} $ to $\phi^\pm_n(z^\pm)$. Lemma \ref{strong cylinder convergence} and our cutting construction above tells us that for any $\delta<1$, there exists some constant $c$ independent of $n$ so that 
\[\abs{\phi_n^{i\pm}(\tilde z^\pm)}\leq c\abs {\tilde z^\pm}^\delta\]  
 and the same inequality holds for  $\phi^\pm_n$ and all derivatives (with a different constant $c$). This implies that $\phi^{i\pm}_n$ converges to $\phi^\pm_n$ in $C^{\infty,\delta}$ for any $\delta<1$. A quick calculation then shows that $G_{(\phi^{i+}_n,\phi^{i-}_n)}(\tilde z^+,\tilde z^-)$ converges in $C^{\infty,\delta}$ to $G_{(\phi^+_n,\phi^-_n)}(\tilde z^+,\tilde z^-)$. Adding the extra coordinates on which $G_{(\phi^{i+}_n,\phi^{i-}_n)}$ does not depend does not affect $C^{\infty,\delta}$ convergence, so we get that $\hat f^i$ converges to $\hat f$ in $C^{\infty,\delta}$ as required. 

 \stop

\

This completes  the proof of Theorem \ref{completeness theorem}.

\bibliographystyle{plain}
\bibliography{ref.bib}

\end{document}